\documentclass[12pt]{amsart}

\usepackage{amsmath}
\usepackage{amssymb}
\usepackage{graphicx}
\usepackage{pifont}
\usepackage{epsfig}
\usepackage{verbatim}
\usepackage[english]{babel}
\usepackage{fancyhdr, blindtext}
\usepackage{float}
\usepackage{subfig}
\usepackage{xcolor}

\usepackage{graphics}
\usepackage{graphicx,epsfig}
\usepackage{epsfig}
\usepackage{graphicx}

\usepackage{a4wide}
\usepackage{amsmath}
\usepackage{amsfonts}

\usepackage{enumerate}
\usepackage{authblk}
\usepackage{float}
\usepackage{subfig}
\usepackage{amsaddr}

\graphicspath{{./arxiv_figures/}}

\begin{document}

\newcommand\bes{\begin{eqnarray}}
\newcommand\ees{\end{eqnarray}}
\newcommand\bess{\begin{eqnarray*}}
\newcommand\eess{\end{eqnarray*}}
\newcommand{\ve}{\varepsilon}
\newtheorem{definition}{Definition}
\newtheorem{theorem}{Theorem}[section]
\newtheorem{lemma}{Lemma}[section]
\newtheorem{proposition}{Proposition}[section]
\newtheorem{remark}{Remark}[section]
\newtheorem{corollary}{Corollary}[section]
\newtheorem{example}{Example}[section]
\title[Oscillatory patterns in a three-timescale predator-prey model]
{\bf Complex oscillatory patterns in a three-timescale model of a generalist predator and a specialist predator competing for a common prey.}
\author{Susmita Sadhu}
\address{Department of Mathematics,
Georgia College \& State University, Milledgeville, GA 31061, USA.}
\email{susmita.sadhu@gcsu.edu}


\thispagestyle{empty}


\begin{abstract}
\noindent   
In this paper, we develop and analyze a model that studies the interaction between a specialist predator (one that relies exclusively on a single prey species), a generalist predator (one that takes advantage of alternative food sources in addition to consuming the focal prey species), and their common prey in a two-trophic ecosystem featuring three timescales. We assume that the prey operates on a faster timescale, while the specialist  and generalist predators operate on slow and superslow timescales respectively. Treating the predation efficiency of the generalist predator as the primary varying parameter and the proportion of its diet formed by the prey species under study as the secondary parameter, we obtain a host of rich and interesting dynamics, including relaxation oscillations, mixed-mode oscillations (MMOs),  \emph{subcritical elliptic bursting} patterns, \emph{torus canards}, and \emph{mixed-type torus canards}. By grouping the timescales into two classes and using the timescale separation between classes, we apply one-fast/two-slow and two-fast/one-slow analysis techniques  to gain insights about the dynamics. Using the geometric properties and flows of the singular subsystems, in combination with  bifurcation analysis and  numerical continuation of the full system, we classify the oscillatory dynamics and  discuss the transitions from one type of dynamics to the other. The types of oscillatory patterns observed in this model are novel in population models featuring three-timescales; some of which qualitatively resemble natural cycles in small mammals and insects.  Furthermore, oscillatory dynamics displaying torus canards, mixed-type torus canards, and MMOs experiencing a delayed loss of stability near one of the invariant  sheets of the self-intersecting critical manifold before getting attracted to the adjacent attracting sheet of the critical manifold have not been previously reported in three-timescale models.

\end{abstract}

\maketitle


\vspace{0.15in}



\section{Introduction}

Understanding patterns of variations in abundance of species is an enduring endeavor in ecology. In many species,  the temporal patterns of their abundance feature multiple timescales and may be broadly viewed as oscillatory dynamics that constitutes of  small amplitude oscillations representing periods of low densities, interspersed with large amplitude oscillations representing episodes of outbreaks. Mixed-mode oscillations (MMOs) or bursting oscillations \cite{DGKKOW, izh, K11} are one such type of complex oscillatory patterns featuring multiple timescales that can represent population cycles bearing resemblance to data from field studies (see figure 1 in \cite{grotan} and \cite{nelson}, figure 2 in  \cite{stenseth} and figure 3 in \cite{singleton}).  Predator-prey models are building blocks for studying population cycles and are commonly used to understand complex interactions in ecological communities; however, there have been relatively few models \cite{BD1} - \cite{BD2}, \cite{BK, KC, KR, MR, MRfoodchain}, \cite{RMpopbio}-\cite{Sadhunew}  that take multiple timescales into account or analyze dynamics involving evolution on three or more distinct timescales. Moreover, the studies on three-timescale population models \cite{BD1} - \cite{BD2}, \cite{ BK, KR, MRfoodchain} have primarily been on tri-trophic food chains  and much is less known about typical dynamics in other types of food-web models featuring three timescales. To address this subject, in this paper, we develop and analyze a two-trophic predator-prey model governing the interaction between  three species, each operating on a different timescale.

Another aspect that is relatively unexplored in continuous-time predator-prey models is the combined effect of specialist predators and generalist predators on the prey dynamics, where the two predators do not engage in intraguild predation, i.e. the two species of predators do not kill/prey upon each other (see  \cite{KW}).  Most existing work on food web models (see \cite{bazykin, ELS, Hanski et al, Hassell et al, seoetal} and the references therein) treat the predator as a true specialist or a generalist, or as an intraguild specialist or an intraguild generalist \cite{KW}, with some models that have considered a shift in predation pattern of the predator from generalist to specialist according to seasonally varying prey availability (i.e. the predator behaves as a generalist in the seasons when several prey species are available but as a specialist in the seasons when few prey species are present) \cite {BLL, TL}.  The presence of both non-intraguild specialist and generalist predators does not seem to have been modeled thus far particularly in ecosystems that lack strong seasonal variations. In this spirit, we propose a  three-species model composed of specialist and generalist predators and their common prey, where the dynamics of the specialist predator is modeled with Holling type II functional response and that of the generalist predator with Holling type III functional response. Such functional responses are typically associated with the predation behaviors of specialist  and generalist predators. We assume that the generalist predator reproduces with Beverton-Holt function in the absence of the common prey  \cite{ELS} and operates on a slower timescale than the specialist predator. With these assumptions, we study the dynamics of the species in the framework of singularly perturbed system of equations, where the prey evolves on a faster timescale, while the specialist and generalist predators evolve on intermediate and slow timescales respectively. Examples of species modeled by such a system may include small mammals such as rodents preyed on by small mustelids or canids (specialist predators)  \cite{Hanski et al}  and large avian predators (generalist predators), or insects attacked  by  specialist parasitoids, and generalist predators such as insectivorous birds or small mammals or arachnids \cite{Hassell et al}.

\begin{figure}[h!]     
  \centering 
{\includegraphics[width=12.67cm]{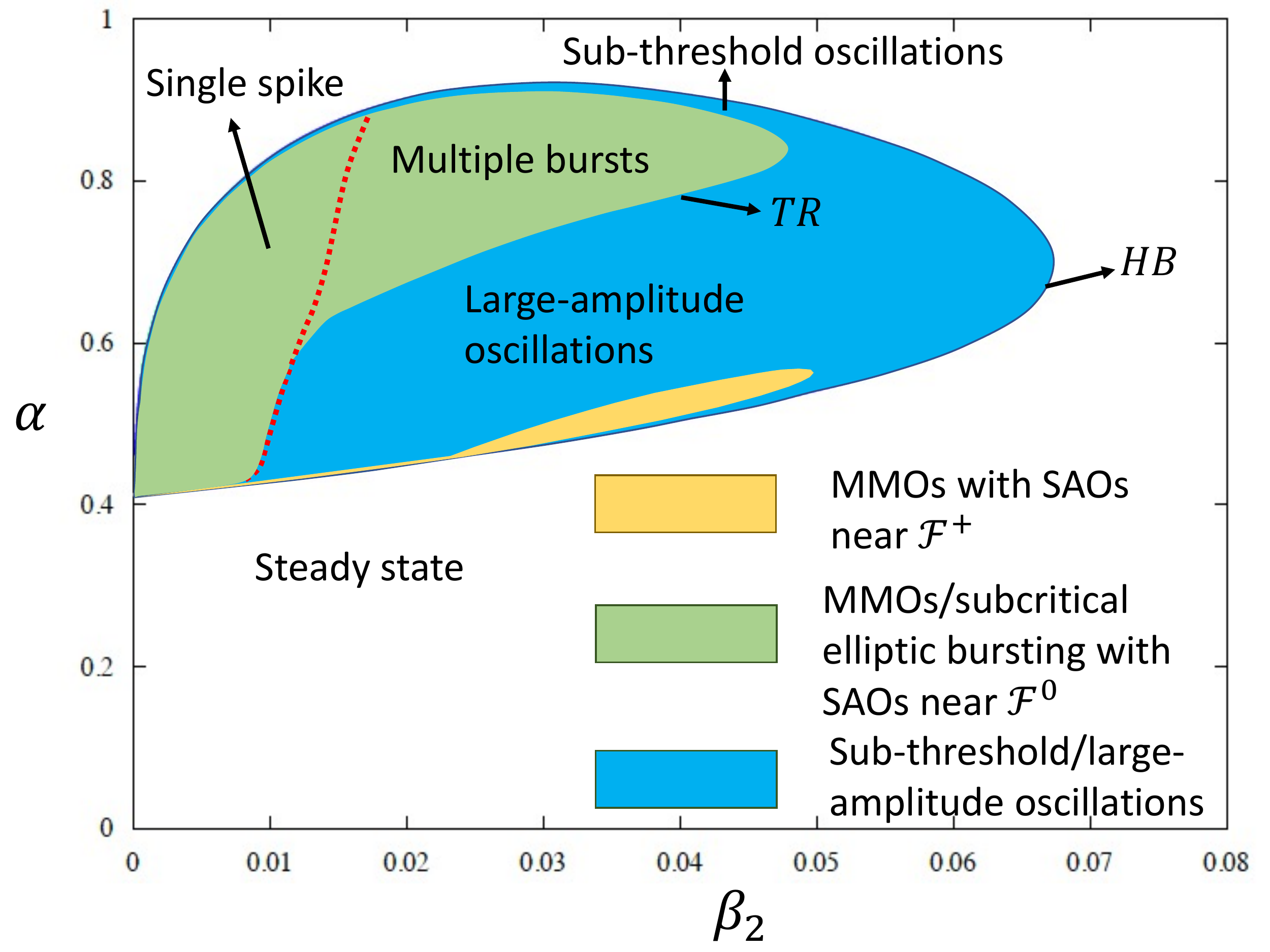}}
\caption{ Two-parameter bifurcation diagram of (\ref{nondim2}) showing transitions between different dynamical regimes for all parameter values as in (\ref{parvalues}). See text for details. HB: Hopf bifurcation, TR: torus bifurcation.}
\label{two_par_bif_full}%
\end{figure}

Treating the predation efficiency of the generalist predator as the primary control parameter and the fraction of generalist predator's  diet that consists of the particular prey species of interest as the secondary parameter, we explore the dynamics of the model, and find a variety of interesting oscillatory patterns such as MMOs, \emph{subcritical elliptic bursting} (subHopf/fold cycle) or \emph{subHopf/subHopf} bursting \cite{izh, TTVWB, VBW}, and relaxation oscillations as shown in figure \ref{two_par_bif_full}. We explain the mechanisms underlying these dynamics using geometrical singular perturbation theory (GSPT) \cite{CT, F, K11} and bifurcation analysis. We get a significant insight about the dynamics by grouping the timescales into two classes and utilizing the timescale separation between the classes using GSPT (see \cite{BD1, kpnew, KPK, LRV} for some examples of three-timescale models where this approach has been used). In one case, we partition the system into slow and fast subsystems in which the fast subsystem consists of a single fast variable and the slow subsystem includes the remaining relatively slower variables and perform the technique of ``one-fast/two-slow analysis". In the other case, the system is divided into a two-dimensional fast subsystem and a one-dimensional slow subsystem and we perform the technique of ``two-fast/one-slow analysis". We study the roles of the critical and superslow manifolds in shaping the dynamics, and explore the bifurcation structures of the two equivalent systems in their singular limits.  An important component of our work includes an exploration of the dynamics of the two-dimensional fast subsystem. Using the geometry of the model and the flows of the lower-dimensional subsystems, we then analyze the different characteristics of the solutions in the three-timescale framework. Finally, we utilize the bifurcation structure of the full system to investigate the parameter dependence of the nature of emergent solutions.

The interplay between the three timescales and coexistence of several mechanisms that are   associated with either two-fast/one-slow systems or one-fast/two-slow systems make the analysis challenging in  three-timescale systems. For instance, the theories of generalized canard phenomenon, singular Hopf bifurcation and delayed Hopf bifurcation  \cite{BKW, DGKKOW, neish1, neish2, Wesc} provide theoretical basis for understanding mechanisms responsible for local oscillatory behavior and bifurcation delay in slow-fast systems with two-timescales. In the present model, it turns out that in a parameter regime where MMOs are observed, a folded node singularity lies in a close vicinity of a delayed Hopf point as well as an equilibrium point of the full system, allowing the mechanisms to interact (c.f. \cite{LRV}). Furthermore, the duration of the quasi-static phase of the subHopf/subHopf type MMO orbits is affected by the relative position of the equilibrium point with respect to a homoclinic bifurcation point of the two-dimensional fast subsystem, adding to more complexity. Interesting dynamics such as torus canards \cite{burketal, Vo} and mixed-type torus canards  \cite{BAD, DBKK}, typically seen in two-timescale systems in neuronal models with at least two fast variables, and not yet been studied in three-timescale systems,  are also observed in this model in vicinities of torus bifurcations. These solutions mark the onset of transition from one kind of oscillatory behavior to another,  and are yet to be fully understood in three timescale systems.

Relaxation oscillation cycles, typically known by boom-and-bust cycles, and chaotic attractors have been extensively studied  in two-timescale and three-timescale ecological models  (see \cite{BD1} - \cite{BD2},  \cite{KR, LXY, MR, MRfoodchain, RMpopbio}  with some recent work on analytical and computational studies on relaxation oscillations \cite{AS, AY},  canard cycles  \cite{PAGK, WZ} and MMOs \cite{KC}, \cite{SCT}-\cite{Sadhunew} in two-timescale predator-prey models. However, to the best of our knowledge, MMOs, torus canards, mixed-type torus canards, and subHopf/fold cycle bursting solutions have not been explored previously in ecological models featuring three timescales.   
Furthermore, torus canards, mixed-type torus, and the MMO dynamics labeled as ``single spike" in figure \ref{two_par_bif_full}, formed by solutions that experience a delayed loss of stability near one of the segments of the self-intersecting critical manifold before they jump toward an adjacent attracting sheet of the other segment of the critical manifold 
 (see figures \ref{crit_mfld}(a), \ref{delayed_hopf}, and \ref{singular_funnel_proj_xz}), are novel in three-timescale settings.  The present work contributes to learning about  different types of  complex oscillatory solutions that can arise in a generic three-timescale predator-prey system, some of which seem to  qualitatively resemble patterns of natural population cycles.

The remainder of the paper is organized as follows. We introduce the model, perform a dimensional analysis, and discuss  the assumptions and physical significance of each parameter in Section 2.  In Section 3, we partition the system into fast and slow subsystems  by grouping the timescales into two classes and perform detailed  fast-slow analyses on these systems using techniques from GSPT.   Combining the two techniques, a GSPT analysis is performed on the full three-timescale model in Section 4. We investigate the bifurcation structure of the model and partition the parameter space into different regions based on the type of oscillatory dynamics. We conclude with a discussion in Section 5.


\section{The Model}

The model studied in this paper reads as follows:
\begin{eqnarray}\label{maineq} \left\{ 
\begin{array}{ll}\label{1}
      \frac{dX}{dT} &= rX\left(1-\frac{X}{K}\right)-\frac{p_1XY}{{H_1}+X}-\frac{\alpha p_2X^2Z}{{H^2_2}+X^2}\\
        \frac{dY}{dT} &= \frac{b_1p_1XY}{{H_1}+X}-d_1Y -m_1Y^2\\
         \frac{dZ}{dT} &= \alpha\left(\frac{b_2p_2X^2Z}{{H^2_2}+X^2}-d_2Z \right)+(1-\alpha)\left(\frac{qZ}{1+m_2Z}-d_3Z\right)
       \end{array} 
\right. \end{eqnarray}
under the initial conditions
\bes \label{maineqic}  X(0)=\tilde{X}\geq 0,\ Y(0)=\tilde{Y}\geq 0, \ Z(0)=\tilde{Z}\geq 0,
\ees
where $X$ represents the population density of the prey and $Y$, $Z$ represent the densities of the two species of  predators.  We assume that $Y$ is a true specialist predator, whereas $Z$ is a generalist predator that does not rely exclusively on $X$ for its food source. The parameters $r$ and $K$ represent  the intrinsic growth rate and the carrying capacity of the prey respectively, $p_1$ is the maximum per-capita predation rate of $Y$, $H_1$ is the semi-saturation constant which represents the prey density at which $Y$ reaches half of its maximum predation rate ($p_1/2$),  $b_1$,  $d_1$ and $m_1$  are respectively the birth-to-consumption ratio,  per-capita natural death rate and density-dependent mortality rate of $Y$. The other parameters $p_2, b_2, d_2$, and $H_2$ are defined analogously for $Z$. In the absence of $X$, we assume that $Z$ reproduces with a Beverton-Holt like function with maximum per-capita reproduction rate $q$ and constant mortality rate $d_3$. We denote the strength of density-dependence in $Z$ by $m_2$. 
 
The net growth rate of $Z$ is considered as a weighted sum of its net growth rates resulting  from consumption of $X$ and other alternative resources with the weight parameter $\alpha$.  
A similar approach was taken in a two-seasons models in \cite{TL}, where the weight parameter was related to the relative length of seasons. In this model, $\alpha$ will be interpreted as the proportion of diet of $Z$ that consists of $X$ and will vary between $0$ and $1$. 
With the following change of variables and parameters:
\begin{eqnarray} {\nonumber} t&=&rT, \  x= \frac{X}{K}, \ y=\frac{p_1Y}{rK},\ z=\frac{p_2Z}{rK}, \ \varepsilon_1 = \frac{b_1p_1}{r},  \ \varepsilon_2 = \frac{b_2p_2}{r}, \  \varepsilon_3 = \frac{q}{r}, \\
{\nonumber}  \beta_1 &=&  \frac{H_1}{K}, \ 
  \beta_2 =\frac{H_2}{K}, \ \delta_1 = \frac{d_1}{b_1p_1},\  \delta_2= \frac{d_2}{b_2p_2}, \  \delta_3= \frac{d_3}{q}, \ 
    \gamma_1= \frac{m_1Y_0}{b_1p_1},  \ \gamma_2 = {m_2Z_0}, 
\end{eqnarray}
where
 \bess   Y_0=\frac{rK}{p_1},\  Z_0=\frac{rK}{p_2}, \eess
system $(\ref{maineq})$ takes the following  dimensionless form:
\begin{eqnarray}\label{nondim1}    \left\{
\begin{array}{ll} \dot{x}&= x\left(1-x-\frac{y}{{\beta_1}+x}-\frac{\alpha xz}{{\beta^2_2}+x^2}\right)\\
 \dot{y}&=\varepsilon_1 y\left(\frac{x}{{\beta_1}+x}-\delta_1-  \gamma_1 y \right)\\
   \dot{z} &= \alpha \varepsilon_2  z\left(\frac{x^2}{{\beta^2_2}+x^2}-\delta_2 \right)+(1-\alpha) \varepsilon_3 z\left(\frac{1}{1+{\gamma}_2 z}-\delta_3 \right),
       \end{array} 
\right. \end{eqnarray}
where the overhead dots denote differentiation with respect to the time variable $t$. 
The quantities $Y_0$ and $Z_0$  will be interpreted as the maximum predation capacities of  $Y$ and $Z$ respectively (see \cite{BD1, Sadhudcds}). We will assume the following conditions on the parameters: 

(A) The maximum per capita birth rate of the prey is much higher than the per capita birth rate of the predators $Y$ and $Z$   i.e. $b_1p_1 \ll r$,  $b_2p_2\ll r$ and $q\ll r$ and that $q =O( b_2p_2)$. For simplicity, we will assume that $q=b_2p_2$. We will further assume that $ b_2p_2 \ll b_1p_1 \ll r$, thus yielding   $0<\varepsilon_2 =\varepsilon_3 \ll  \varepsilon_1 \ll1$.

(B) We will assume that  $Z$ can persist even in the absence of $X$ and that the parameters $\delta_1$, $\delta_2$ and $\delta_3$ satisfy the inequality $0<\delta_1, \delta_2, \delta_3<1$, which implies that the growth rates of the predators are greater than their death rates. This is a default assumption otherwise the predators would die out faster than they could reproduce even at their maximum reproduction rates.

(C) The parameters $\beta_1$ and $\beta_2$ are dimensionless semi-saturation constants measured against the prey's carrying capacity. We will assume that both predators are efficient, and will reach the half of their maximum predation rates before the prey population reaches its carrying capacity, thus yielding $0< \beta_1,\beta_2<1$.

Under the assumptions (A)-(C), system $(\ref{nondim1})$  transforms to a singular perturbed system of equations with three timescales, where the prey exhibits fast dynamics, the specialist predator exhibits intermediate dynamics while the generalist predator exhibits slow dynamics. 

We rewrite system $(\ref{nondim1})$  as
\begin{eqnarray}\label{nondim22}    \left\{
\begin{array}{ll}  \dot{x} &= x\left(1-x-\frac{y}{{\beta_1}+x}-\frac{\alpha xz}{{\beta^2_2}+x^2}\right):=x\phi(x,y,z,\rho)\\
    \dot{y}&= \varepsilon_1 y\left(\frac{x}{{\beta_1}+x}-\delta_1-\gamma_1 y \right):=\varepsilon_1 y\chi(x,y,\rho)\\
    \dot{z}&=  \varepsilon_2 z\left(\alpha \left(\frac{x^2}{{\beta^2_2}+x^2} - \delta_2  \right) +(1-\alpha) \left(\frac{1}{1+\gamma_2 z}-\delta_3 \right)\right):= \varepsilon_2 z\psi(x,z, \rho),
      \end{array} 
      \right. 
\end{eqnarray}
where $\phi=0$, $\chi=0$, and $\psi=0$ are the nontrivial $x$, $y$, and $z$-nullclines respectively and $\rho=(\alpha, \beta_1, \beta_2, \delta_1, \delta_2, \delta_3, \gamma_1, \gamma_2)\in \mathbb{R}^8$ is a vector of parameters. For simplicity, from here onwards we will denote the ratio $\varepsilon_2/\varepsilon_1$ by $\delta$ and replace $\varepsilon_1$ by $\varepsilon$. By assumption (A), we then have that $0< \varepsilon, \delta \ll 1$ and system (\ref{nondim22}) can be  rewritten as 
\begin{eqnarray}\label{nondim2}    \left\{
\begin{array}{ll}  \dot{x} &=x\phi(x,y,z,\rho)\\
    \dot{y}&= \varepsilon y\chi(x,y,\rho)\\
    \dot{z}&=  \varepsilon \delta  z\psi(x,z, \rho).
      \end{array} 
      \right. 
\end{eqnarray}
With respect to the timescale $t$, system (\ref{nondim2}) will be referred to as the fast system. Rescaling time by $\varepsilon$ and letting $s = \varepsilon t$, one obtains the equivalent intermediate system 
\begin{eqnarray}\label{inter}    \left\{
\begin{array}{ll}  \varepsilon  \dot{x}  &=x\phi(x,y,z,\rho)\\
  \dot{y} &= y\chi(x,y,\rho)\\
\dot{z} &=  \delta  z\psi(x,z, \rho), 
      \end{array} 
      \right. 
\end{eqnarray}
where the over dot represents derivative with respect to $s$. Finally, by rescaling time as $\tau= \varepsilon \delta t$ in (\ref{nondim2}), one obtains the equivalent slow system
\begin{eqnarray}\label{slow}    \left\{
\begin{array}{ll}  \varepsilon \delta \dot{x} &=x\phi(x,y,z,\rho)\\
\delta \dot{y} &= y\chi(x,y,\rho)\\
\dot{z} &=  z\psi(x,z, \rho),
      \end{array} 
      \right. 
\end{eqnarray}
where the over dot represents derivative with respect to $\tau$.

 Fixing $\delta$ and treating $\varepsilon$ as the singular parameter, system (\ref{inter}) partitions as a fast-slow system with one fast variable $x$ and two slow variables $y$ and $z$. On the other hand, keeping $\varepsilon$ fixed and treating $\delta$ as the singular perturbation parameter, system (\ref{inter}) partitions into a family of two-dimensional $(x, y)$  fast-subsystems parametrized by $z$. 

  \begin{figure}[h!]     
  \centering 
\subfloat[]{\includegraphics[width=7.7cm]{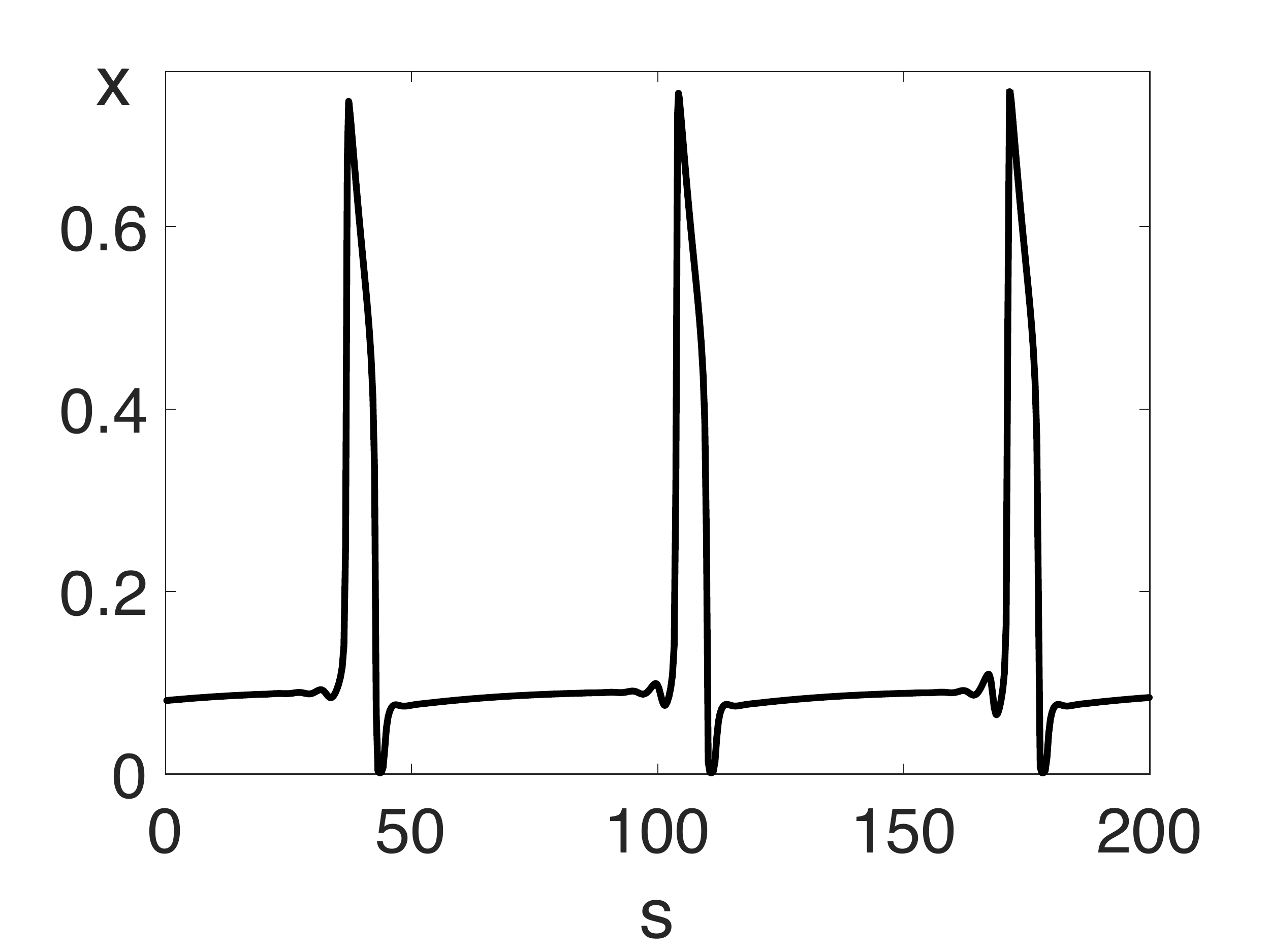}}
\quad
\subfloat[]{\includegraphics[width=7.7cm]{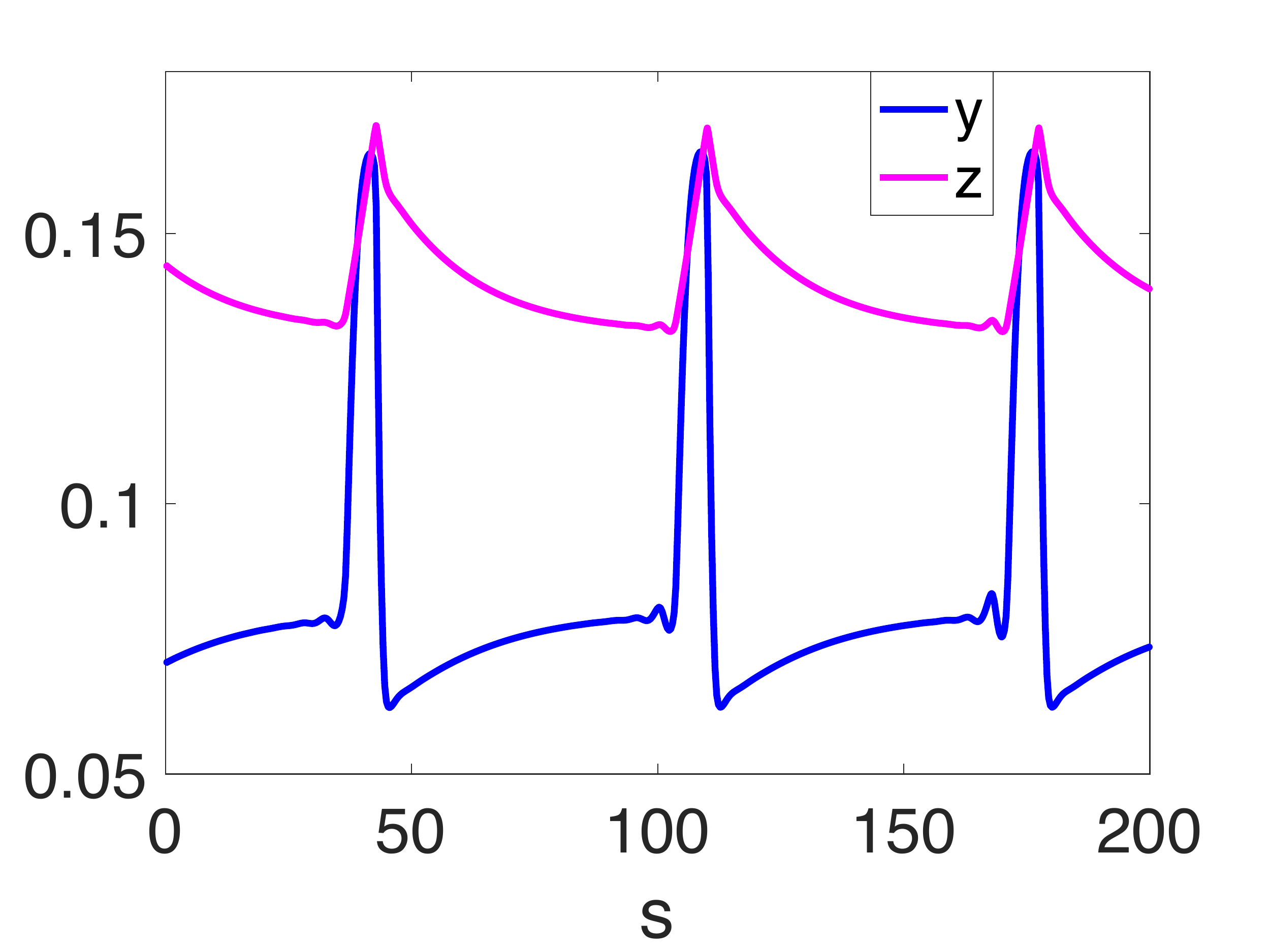}}
\caption{ Mixed-mode time series of (\ref{inter}) in intermediate time for $\beta_2=0.01$, $\alpha=0.75$ and other parameters values as in (\ref{parvalues}). (A) For $x$-component. (B) For $y$ and $z$ components.}
\label{timeseries_example_1}%
\end{figure}

\begin{figure}[h!]     
  \centering 
\subfloat[]{\includegraphics[width=7.7cm]{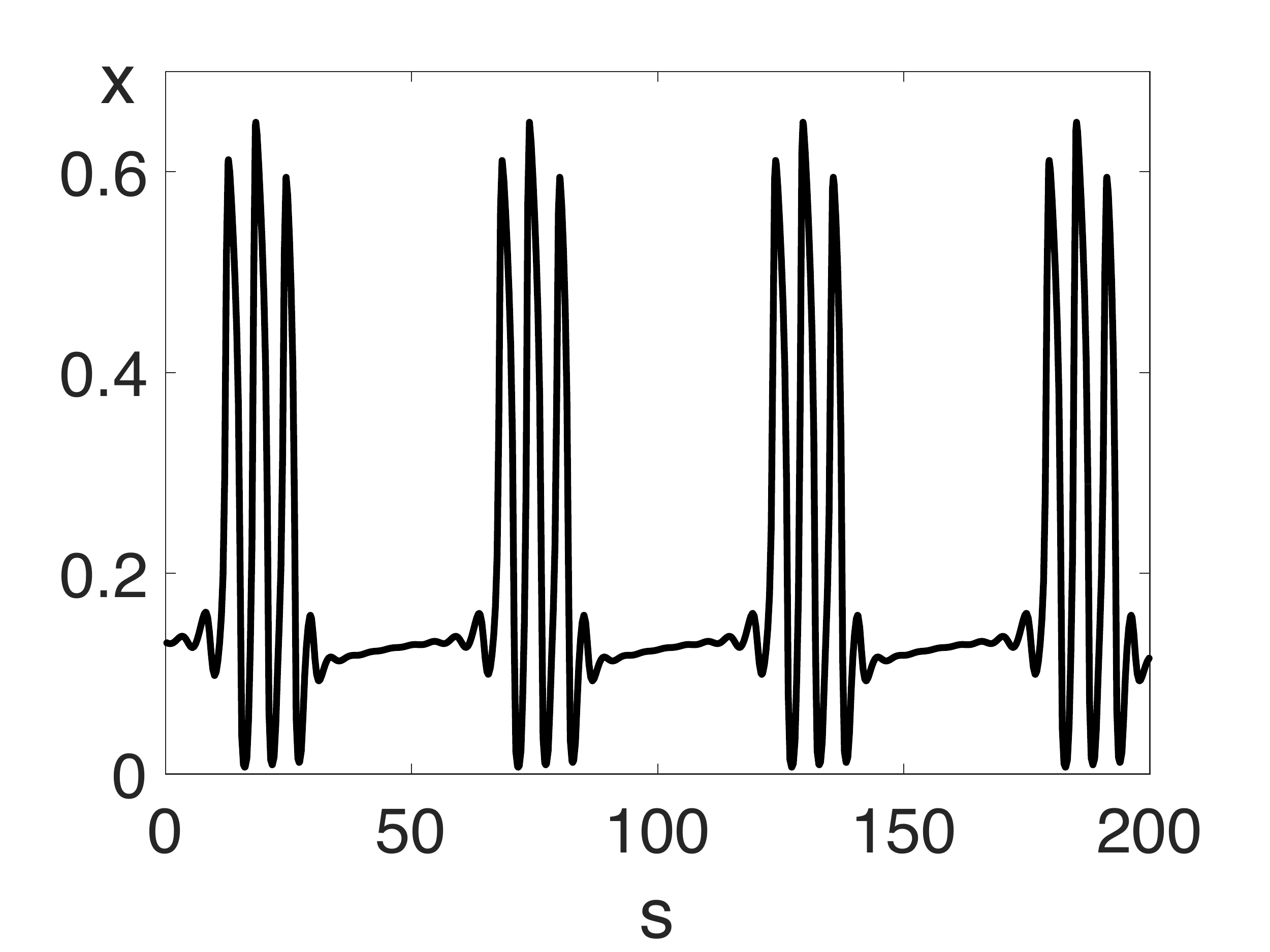}}
\quad
\subfloat[]{\includegraphics[width=7.7cm]{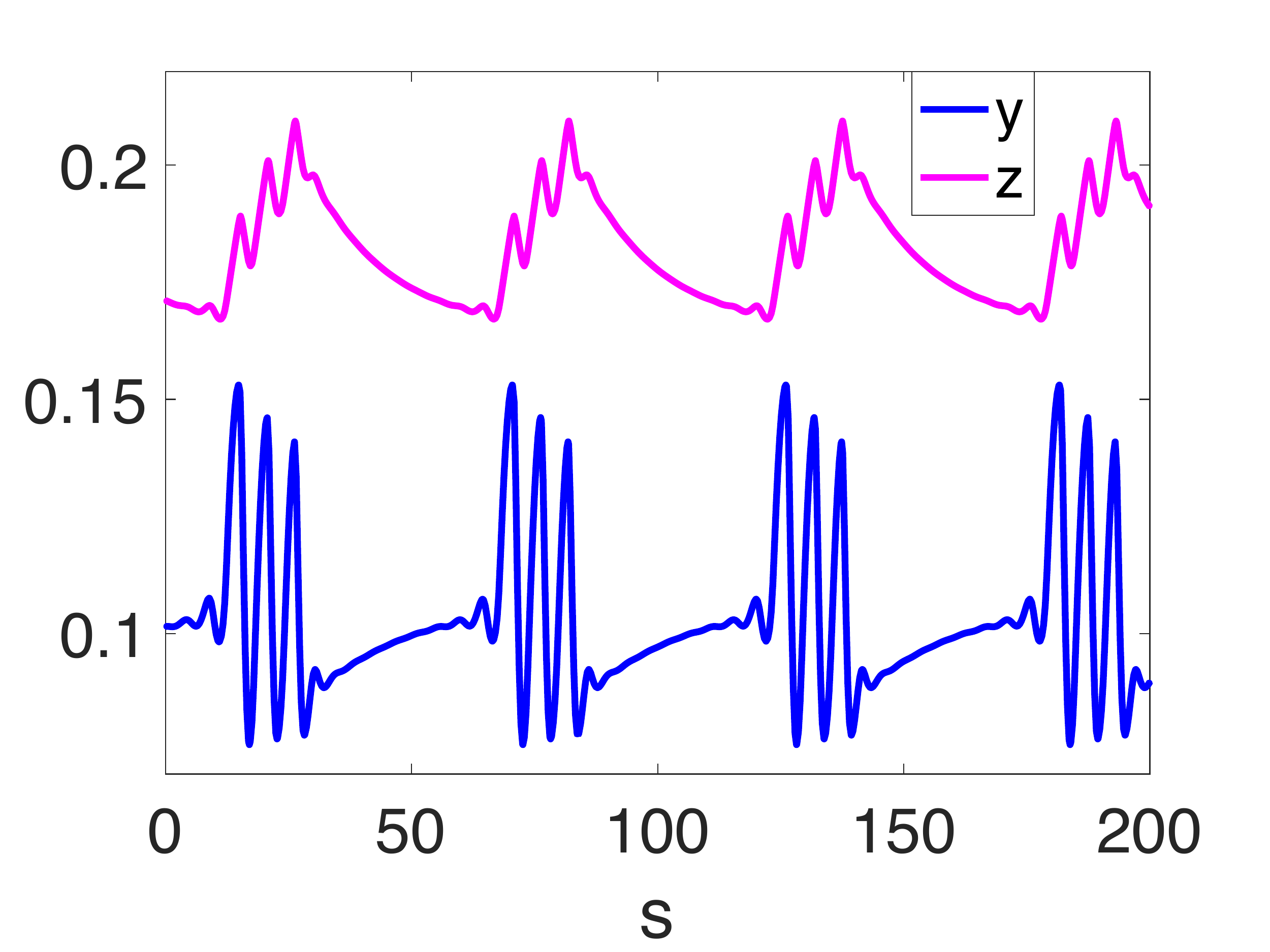}}
\caption{Bursting patterns in system (\ref{inter}) for $\beta_2=0.0245$, $\alpha=0.8$ and other parameters values as in (\ref{parvalues}). (A) For $x$-component. (B) For $y$ and $z$ components.}
\label{timeseries_example_2}%
\end{figure}

Throughout the paper, we will fix the parameter values to
\bes \label{parvalues}
\varepsilon =0.05, \ \delta=0.1, \ \beta_1 =0.1,\ \delta_1=0.15,\ \delta_2=0.35,\ \delta_3=0.65,\  \gamma_1=4.1, \ \gamma_2=15, 
\ees
and vary the predation efficiency $\beta_2$ of the generalist predator as the primary control parameter and the weight $\alpha$ as the secondary parameter.  For the choice of  parameter values in (\ref{parvalues}),  the intersection of the non-trivial nullclines $\phi=0, \chi=0$ and $\psi=0$ produces equilibria in the positive octant. These equilibria will be referred to as the coexistent or non-trivial equilibria and will be denoted by $E^*$. The equilibria lying on the invariant $xy$-plane and the $xz$-plane will be referred to as the boundary equilibria and will be denoted by $E_{xy}$ and $E_{xz}$ respectively. The parameter values chosen here are for illustrative purposes to demonstrate the different types of oscillatory patterns that arises in this system. Representative time profiles of dynamics of system (\ref{inter}) are shown in figures \ref{timeseries_example_1} and \ref{timeseries_example_2}. The phase portraits of the trajectories are shown in figures \ref{singular_funnel_proj_xz} and \ref{crit_mfld}(B) respectively.

\begin{remark} \label{rmknew}  In most ecosystems, it will be reasonable to assume that $\varepsilon \in (0.1, 0.01)$ and $\delta \in (0.2, 0.05)$. For a fixed $\varepsilon$,  and letting $\delta = O(\varepsilon^k)$, we find that the oscillatory patterns in system  (\ref{inter}) are robust for $k>1/2$. Note that  $k\approx 3/4$ in (\ref{parvalues}). It turns out that the three-timescale structure and some of the oscillation patterns are lost if $k$ is chosen to be less than $1/2$ (c.f. \cite{KPK,LRV}). For instance, corresponding to the parameter values in figure  \ref{timeseries_example_1} or figure \ref{timeseries_example_2},  system (\ref{inter})  exhibits relaxation oscillations featuring two-timescales  if $\delta = O(\varepsilon^{1/4})$.
\end{remark}

\section{THE GEOMETRIC SINGULAR PERTURBATION  THEORY APPROACH}

  In this section, we will use  geometric singular perturbation theory and apply Fenichel theory iteratively \cite{CT} to explain the mechanisms underlying the complex dynamics exhibited by system $(\ref{inter})$.  
   \subsection{One-fast/two-slow analysis} 
Fixing $\delta$ and letting $\varepsilon \to 0$, we will decompose system $(\ref{inter})$ into a two-dimensional slow subsystem and a one-dimensional fast subsystem. We will then analyze the planar slow subsystem and study key structures such as the {\emph{critical manifold}}, defined by the equilibria of the fast subsystem. The analysis will be used to study canard-induced MMOs in the full system.
 
       \begin{figure}[h!]     
  \centering 
\subfloat[]{\includegraphics[width=7.7cm]{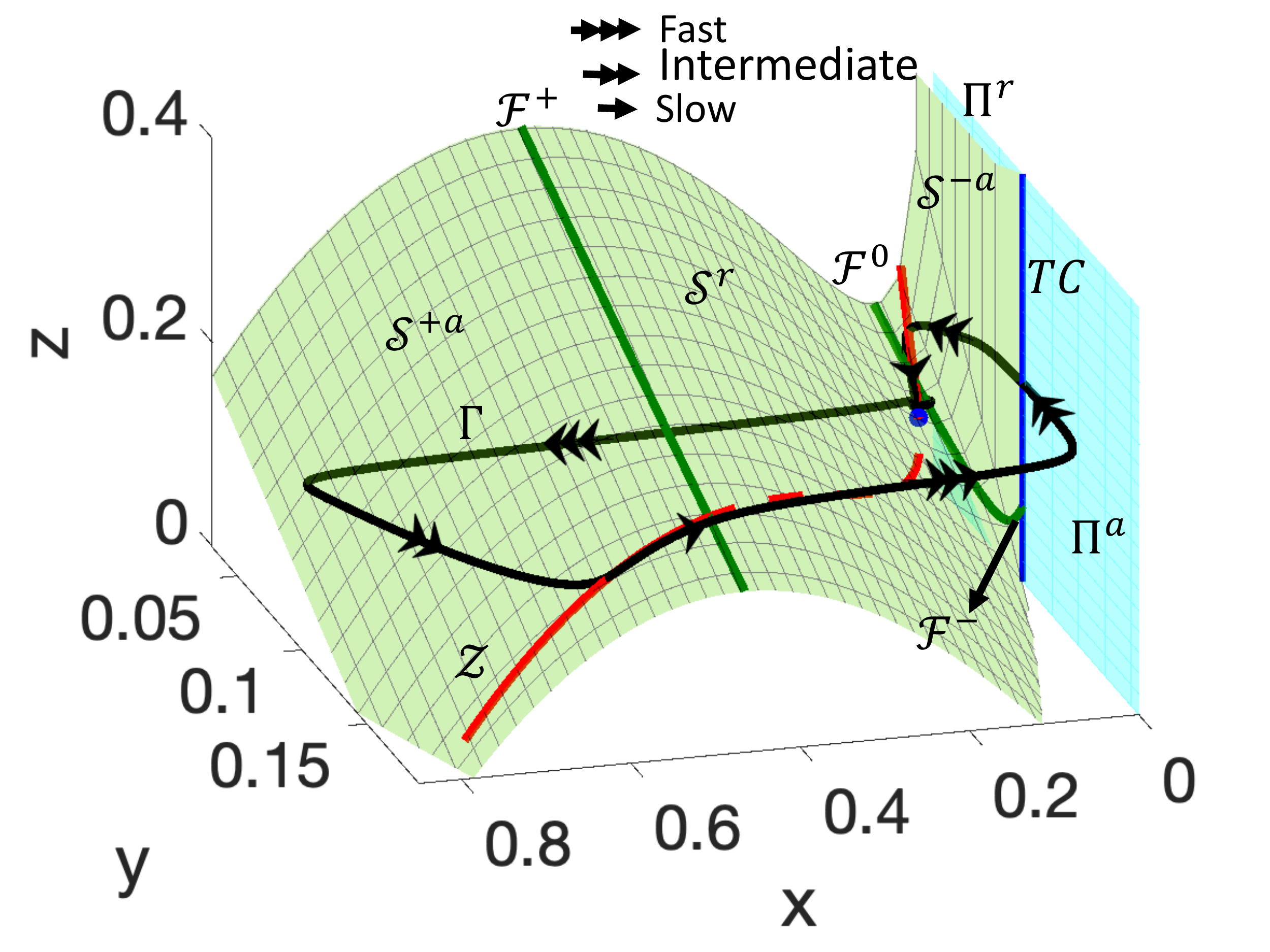}}
\quad
\subfloat[]{\includegraphics[width=7.7cm]{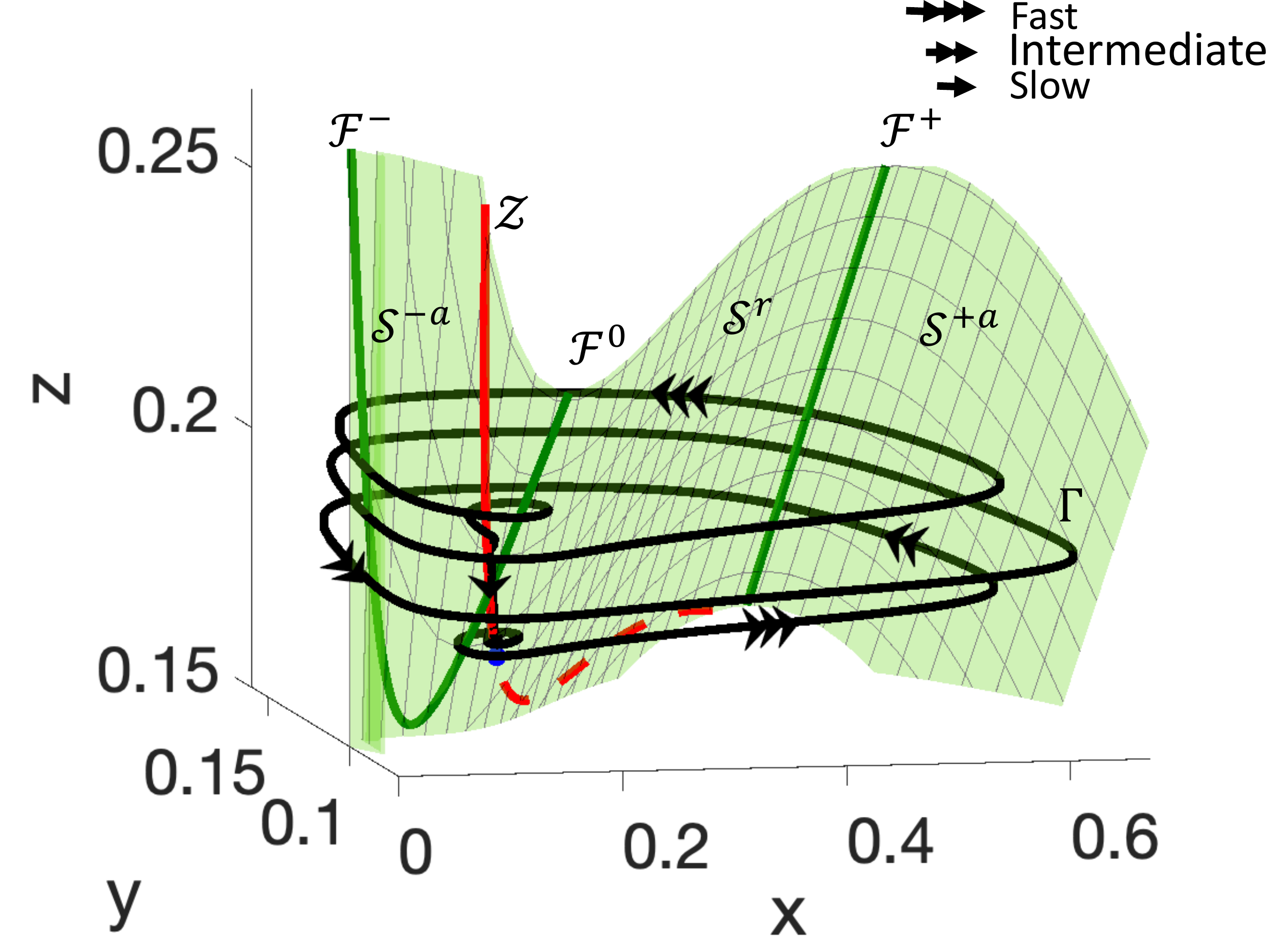}}
\caption{Geometric structure and stable periodic orbits $\Gamma$ of system (\ref{inter}) for different values of $\beta_2$ and $\alpha$ and other parameters as in (\ref{parvalues}). Shown are the critical manifold $\mathcal{M}_1=\Pi \cup S$, the superslow manifold $\mathcal{Z}$ (red) and the coexistence equilibrium state (blue dot). The plane $\Pi$ is divided into an attracting component $\Pi^a$ and a repelling component $\Pi^r$ joined at the transcritical curve $TC$.  The surface $S$ can consist of one or more attracting and repelling sheets joined at the folds $\mathcal{F}^{0}$ and $\mathcal{F}^{\pm}$ (also see table \ref{bif_foldcurve}). Note that $\Gamma$ displays three distinct timescales during the course of its cycle. (A) $\beta_2=0.005$ and $\alpha=0.6$. (B) $\beta_2=0.0245$ and $\alpha=0.8$. } 
\label{crit_mfld}%
\end{figure}

  \subsubsection{The critical manifold $\mathcal{M}_1$}
  By fixing $\delta$ and letting $\varepsilon \to 0$ in system $(\ref{nondim2})$ results into a one-dimensional {\emph{fast-subsystem}}
  \bes
  \label{layer1d}
  \dot{x} = x\phi(x,y,z,\rho),
  \ees where $y$ and $z$ are parameters. The set of equilibria of the fast subsystem defines a two-dimensional manifold called the {\emph{critical manifold}}, $\mathcal{M}_{1,\delta}$, where
\bess \mathcal{M}_{1,\delta}=\left\{(x,y,z): x=0 {\mathrm{~or~}} \phi(x,y,z,\rho)=0\right\}.
\eess
 The geometry of the critical manifold is independent of $\delta$, hence we will suppress the dependence of $\mathcal{M}_{1,\delta}$ on $\delta$ and denote it by $\mathcal{M}_{1}$.  Note that $\mathcal{M}_{1}$  consists of two disjoint components, namely the plane $\Pi =\{ (0,y,z): y, z \geq 0\}$ and the curved surface 
 \bess S=\{(x,y,z) \in {\mathbb{R}^3}^+:\phi(x,y,z,\rho)=0\} = \{(x, y, z) \in {\mathbb{R}^3}^+ : y= F(x, z, \rho)\},
 \eess
 where 
 \[F(x, z, \rho) = (\beta_1+x)\Big(1-x-\frac{\alpha xz}{{\beta^2_2}+x^2}\Big).\]
 The surface $S$ intersects the plane $\Pi$ along the line
$TC=\{(0, \beta_1, z): z\geq 0\}$. System (\ref{layer1d}) undergoes a transcritical bifurcation along $TC$ and $\Pi$  is divided into two normally hyperbolic parts, $\Pi^a =\{(0, y, z): y>\beta_1\}$ and $\Pi^r =\{(0, y, z): y<\beta_1\}$ by $TC$ as shown in figure \ref{crit_mfld}. The surface $S$ 
is folded and can be written as $S= S^a \cup S^r \cup \mathcal{F}$, where 
\bess 
S^a=S\cap \{ \phi_x(x, y, z,\rho)<0\} \ \textnormal {and}  \ S^r=S\cap \{ \phi_x(x, y, z,\rho)>0\}
\eess  
are normally attracting and repelling respectively, and $\mathcal{F}$ is degenerate due to loss of normal hyperbolicity generated by saddle-node bifurcations. Namely,
 \bess \mathcal{F} &=& \{(x, y, z) \in S:  \phi_x(x, y, z,\rho)=0\} \\
& =& 
  \left\{(x, y, z) \in S:   y=\mu(x),\  z= \nu(x), \ x \in [0, 1]\setminus \{x_d\} \right\},
 \eess
 where
 \bes
 \label{foldzcomp}
 \mu(x) = F(x, \nu(x)), \  \nu(x)=\frac{(1-\beta_1-2x)(\beta^2_2 +x^2)^2}{\alpha (\beta_1\beta^2_1 +2\beta^2_2 x-\beta_1x^2)}, \ees
and
\bess
  x_d =  \frac{\beta^2_2}{\beta_1}+\sqrt{ \frac{\beta^4_2}{\beta^2_1}+\beta^2_1}.
 \eess
 
    \begin{figure}[h!]     
  \centering 
{\includegraphics[width=12.7cm]{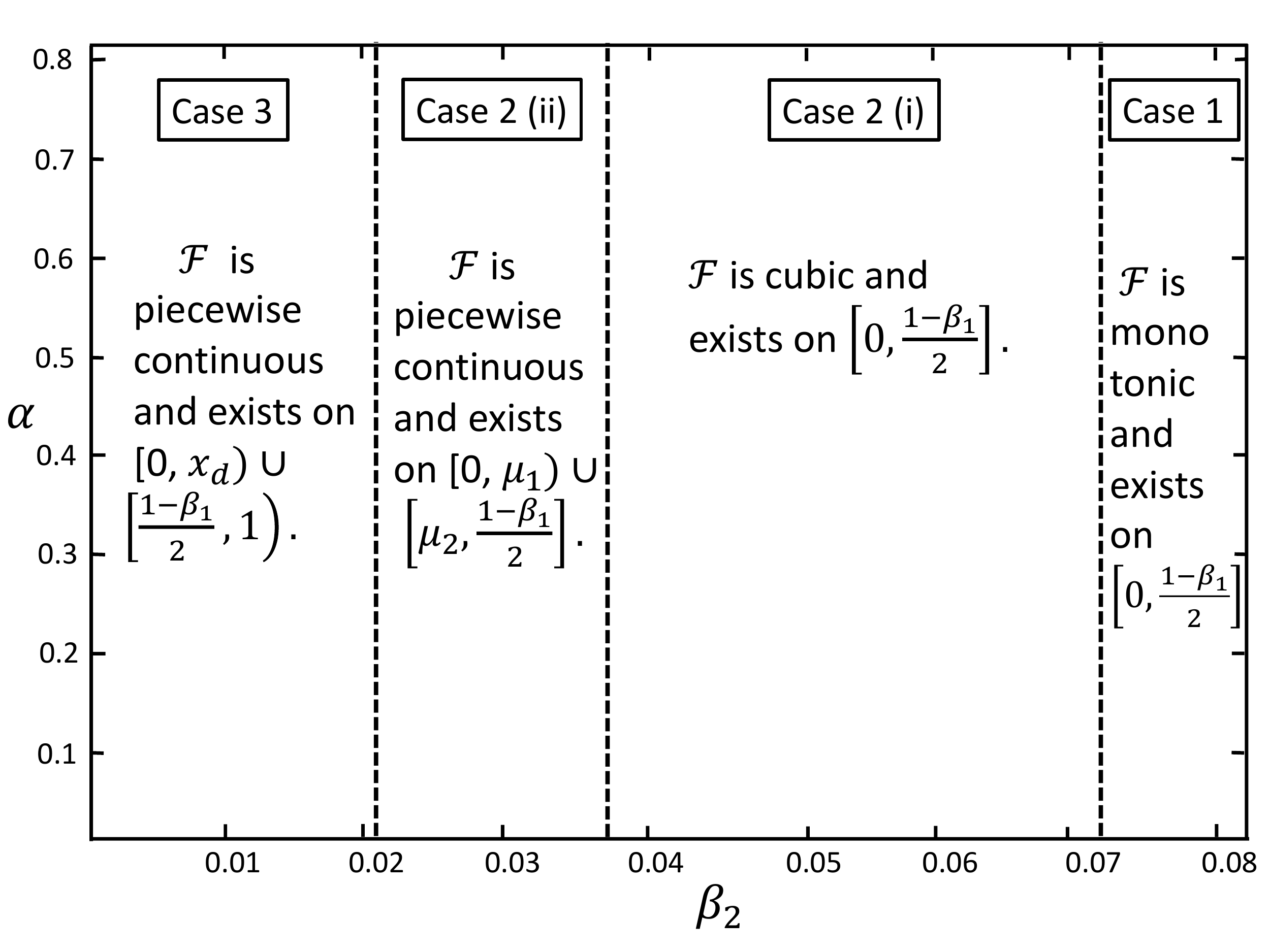}}
\caption{Shape of the fold curve $\mathcal{F}$ with respect to $\alpha$ and $\beta_2$. Here ${\mu}_1,  {\mu}_2$ are the roots of $\mu(x)=0$. See Appendix for the details.} 
\label{bif_foldcurve}%
\end{figure}

\begin{figure}[h!]     
  \centering 
  \subfloat[]{\includegraphics[width=7.7cm]{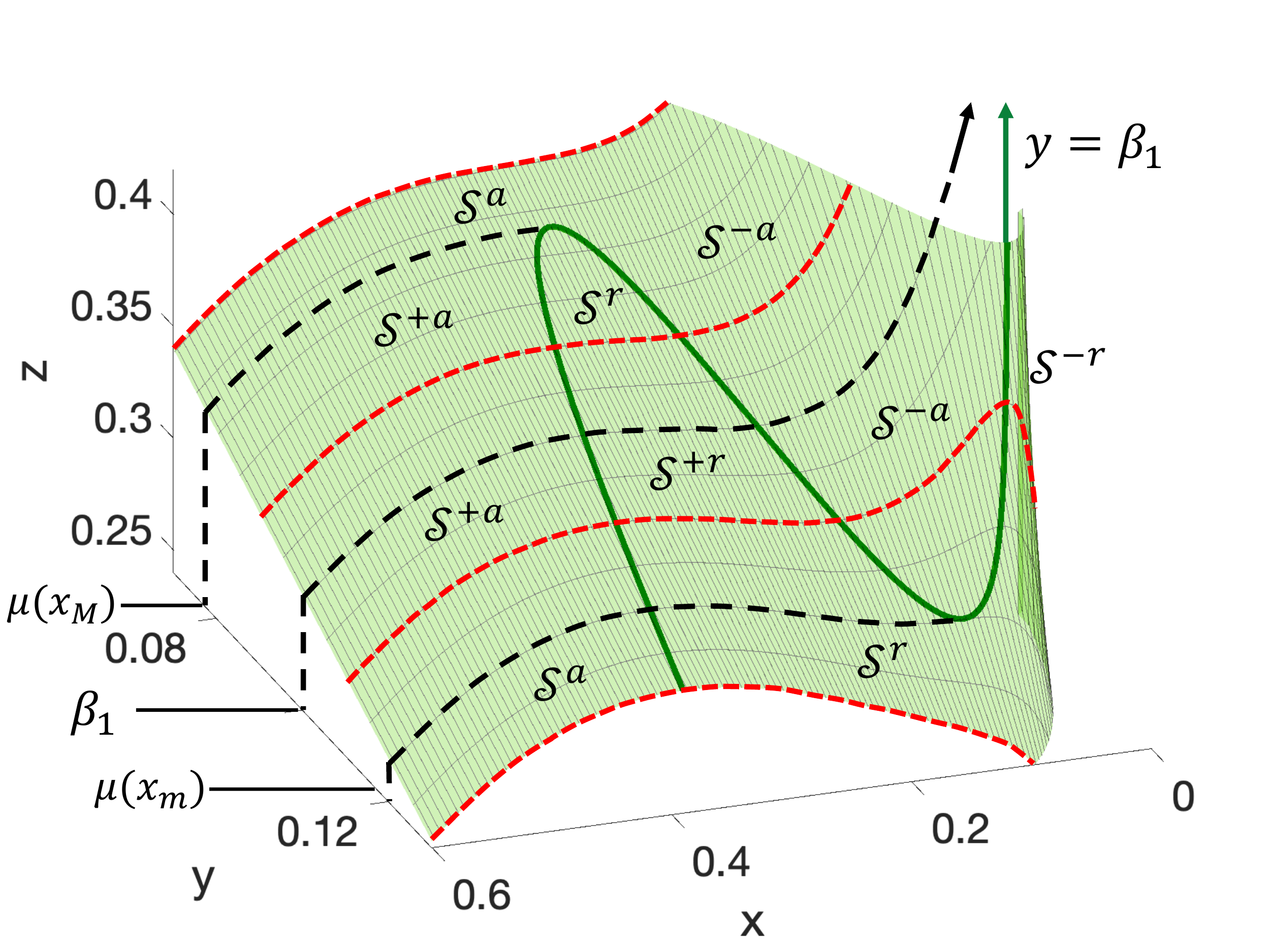}}
  \quad
\subfloat[]{\includegraphics[width=7.7cm]{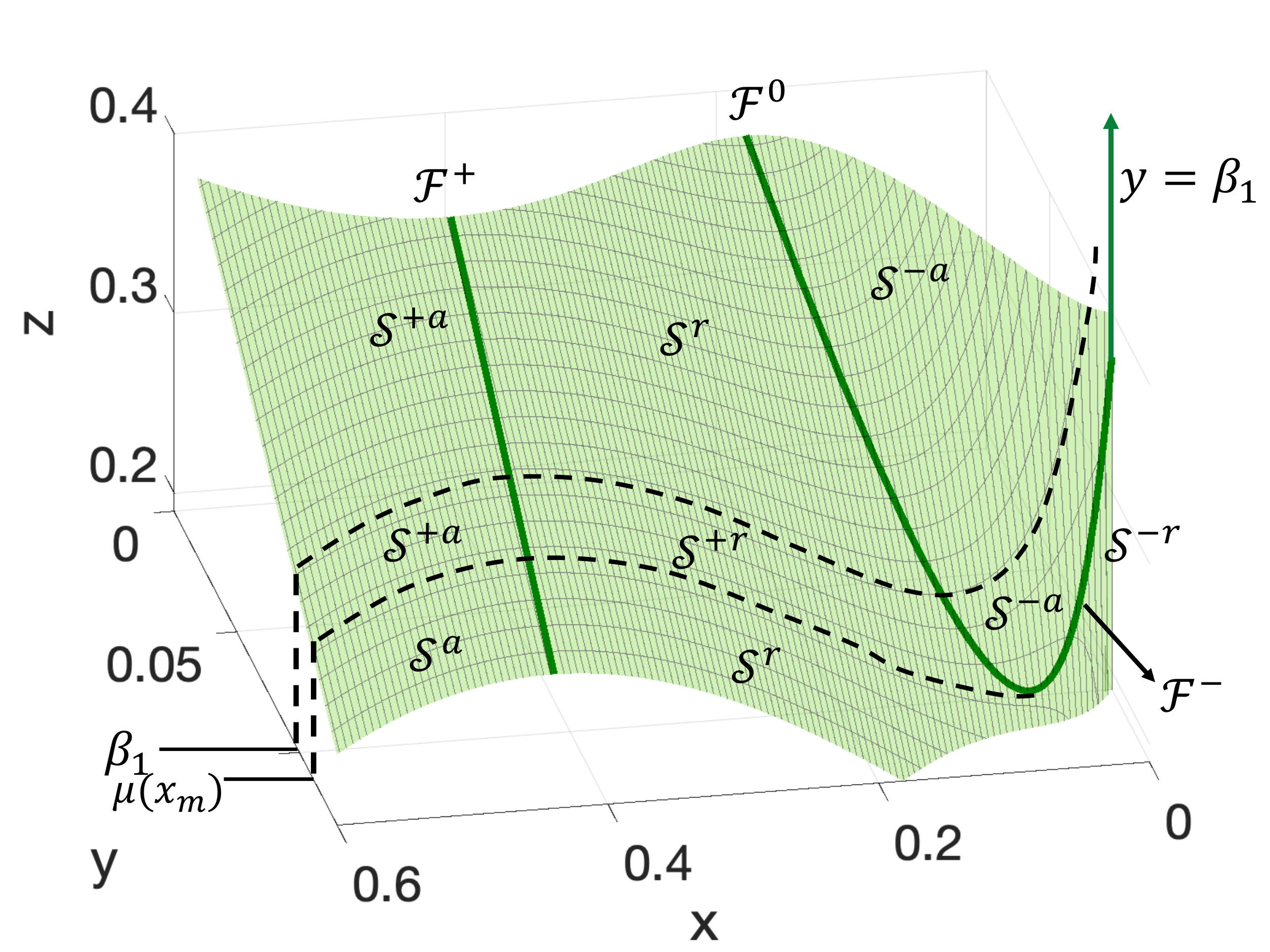}}
\caption{Geometric structure of the surface $S$  of system (\ref{inter}) for different values of $\beta_2$ with $\alpha =0.6$ and other parameters as in (\ref{parvalues}). The black dotted lines divide $S$ into distinct regions, characterized by the number of  attracting and repelling branches that the surface possesses. (A) The surface $S$ is divided into four regions for $\beta_2=0.048$. Also shown are different  cross-sections of $S$ (in red) for constant $y$ values chosen from each region.  (B) The surface $S$ is divided into three regions for $\beta_2=0.025$.}
\label{fold_curves}%
\end{figure}

 For suitable values of $\beta_2$ and $\alpha$, and feasible ranges of $y, z$ such that $\mathcal{F}\subset {\mathbb{R}^3}^+$,  the fold curve may be cubic (i.e. has two folds) (see figure \ref{crit_mfld}(B) and figure \ref{fold_curves}), monotonic or piecewise continuous (see figure \ref{crit_mfld}(A))  determined by the locations of the roots, critical points and the point of discontinuity of $\nu(x)$ (and $\mu(x)$) relative to each other.  The properties of the fold curve as $\alpha$ and $\beta_2$ are varied are illustrated in figure \ref{bif_foldcurve}. Depending on the structure of $S$ characterized by the fold curve, system $(\ref{inter})$ can exhibit MMOs or relaxation oscillations with ``plateaus below" near the $\Pi$ (see figure \ref{timeseries_varying_alpha}) corresponding to Case 3 in figure \ref{bif_foldcurve}  (c.f.\cite{kpnew}),  plateau-less MMOs or relaxation oscillations (see figure \ref{timeseries_varying_alpha_1}) corresponding to Case 2 in figure \ref{bif_foldcurve}  or steady state solutions corresponding to Case 1 in figure \ref{bif_foldcurve}. The details of figure  \ref{bif_foldcurve} are provided in the Appendix.

  \subsubsection{Reduced dynamics on $\mathcal{M}_1$}
 
Taking the singular limit $\varepsilon  \to 0$ in system $(\ref{inter})$  yields the slow subsystem
  \begin{eqnarray}\label{reduced}   \left\{
\begin{array}{ll} 0&=x\phi(x, y, z, \phi) \\
    \dot{y}&=y\chi(x,y,\rho)\\
    \dot{z}&= \delta z\psi(x, \rho).
       \end{array} 
\right. 
\end{eqnarray}
 The flow governed by system (\ref{reduced}) is constrained to the critical manifold $\mathcal{M}_1$ and is referred to as the {\emph{reduced dynamics}}.  On the plane ${\Pi}$, the reduced dynamics solves the system
\begin{eqnarray}\label{plane}   \left\{
\begin{array}{ll} x &=0\\
    \dot{y}&=y\chi(0,y,\rho)\\
    \dot{z}&=  \delta z\psi(0, \rho).
       \end{array} 
\right. 
\end{eqnarray}
The flow descends along $\Pi$ and approaches $(0, 0, 0)$ which is the global attractor of (\ref{plane}). The reduced flow crosses the transcritical line $TC$ from $\Pi^a$ to $\Pi^r$ with finite speed, giving rise to singular canards.  We note that ${\Pi}$ is invariant for all $\varepsilon>0$, hence canards  persist for the full system. As the reduced flow descends along $\Pi^a$ and goes past $TC$, it spends an $O(1)$ amount of time in the intermediate timescale on $\Pi^r$ before it experiences a loss of stability and gets concatenated by a fast fiber to $S^{a}$ as shown in figure \ref{crit_mfld}. This phenomenon of delay is referred to as the Pontragyin's delayed loss of stability and has been studied in a few three-dimensional models (see \cite{kpnew, Sadhudcds}).

We next consider the flow on $S$. As noted earlier, the surface $S$ can be locally expressed as the graph of $y=F(x,z, \rho)$, hence the dynamics of (\ref{reduced}) can be projected  onto the $(x,z)$ coordinate chart. Differentiating $\phi(x,y,z, \rho)=0$ implicitly with respect to time and using the fact that $\dot{\rho}=0$,  gives us the relationship $\phi_x\dot{x}+\phi_y\dot{y}+\phi_z\dot{z}=0$. Thus, the reduced flow (\ref{reduced}) restricted to $S$ reads as
\begin{eqnarray} \label{red2}  \begin{pmatrix}
 -\phi_x\dot{x}  \\
 \dot{z} 
 \end{pmatrix}= \begin{pmatrix}
 \phi_yy \chi+\delta \phi_zz \psi  \\
\delta z\psi
 \end{pmatrix}\bigg|_{y= F(x,z)}.
\end{eqnarray}
 
 System (\ref{red2}) has singularities when $\phi_x=0$ and its solutions blow-up in finite time at $\mathcal{F}$. Hence standard existence and uniqueness results do not hold.  To remove the finite-time blow up of solutions,  we rescale the time $s$ by the factor $-\phi_x$, i.e. $ds = -\phi_x d\tau$ \cite{DGKKOW}, thus transforming system (\ref{red2}) into the {\emph{desingularized system}}
\begin{eqnarray} \label{desing}  \begin{pmatrix}
 \dot{x}  \\
 \dot{z} 
 \end{pmatrix}= \begin{pmatrix}
\phi_yy \chi+\delta \phi_zz \psi  \\
-\delta \phi_x z \psi
 \end{pmatrix}\bigg|_{y=F(x,z)},
\end{eqnarray}
where the overdot denotes $\tau$ derivatives.  We note that system (\ref{desing}) is singularly perturbed with respect to the singular parameter $\delta$. It is topologically equivalent to system (\ref{red2}) on $S^a$. However, the phase-space-dependent time transformation reverses the orientation of the orbits on $S^r$, therefore, the flow of (\ref{red2}) on $S^r$ is obtained by reversing the direction of orbits of (\ref{desing}).   It then follows that the reduced flow on $S$ is either directed towards $\mathcal{F}$ or away from it.
 
By Fenichel's theory \cite{CT, F}, the normally hyperbolic segments of the critical manifold $\mathcal{M}_1$ perturb to locally invariant attracting and repelling slow manifolds ${{\Pi}_{\varepsilon, \delta}}\cup{S_{\varepsilon, \delta}^a}\cup S_{\varepsilon, \delta}^r$  for $\varepsilon> 0$, and the slow flow restricted to these manifolds is an $O(\varepsilon)$ perturbation of the reduced flow on $\mathcal{M}$. However, the theory breaks down in neighborhoods of $\mathcal{F}$.

  \subsubsection{Singular points on $\mathcal{M}_1$}

The set of equilibria $S_E^{\delta}$ of  (\ref{red2}) and (\ref{desing}) that do not lie on the fold curve $\mathcal{F}$ are {\em{ordinary singularities}},  i.e. 
\bess S_E^{\delta} :=\{(x, y, z) \in S\setminus \mathcal{F}: \chi=0 \land z=0\} \cup \{(x, y, z) \in S\setminus \mathcal{F}: y=0 \land \psi=0\}\\
 \cup \{(x, y, z) \in S\setminus \mathcal{F}: \chi=0 \land \psi=0\}.
\eess
 On the other hand, the elements of the set $S_F^{\delta}$ defined by
 \bess
 S_F^{\delta}: =\{(x, y, z) \in \mathcal{F}: \phi_yy \chi+\delta \phi_zz \psi =0\}
 \eess
are called {\em{folded singularities}} or canard points. These points are equilibria  of (\ref{red2}) and (\ref{desing}) and form isolated points of $\mathcal{F}$ (see \cite{DGKKOW} for the classification of folded singularities). Further degeneracies may occur if one of the eigenvalues passes through $0$ and can give rise to \emph{folded saddle node (FSN) bifurcation of types I and II}  \cite{DGKKOW,Wesc}.  
In figure \ref{two_par_desing}, we summarize the variations in bifurcations of  system (\ref{desing}) over a range of values of $\beta_2$ and $\alpha$. In the  bifurcation diagram, the FSN II (a) and FSN II (b) curves correspond to transcritical bifurcations of folded singularities and ordinary singularities $E_{xz}$ (boundary equilibrium state on the $xz$-plane), and folded singularities and ordinary singularities $E^*$ (coexistent equilibrium) respectively. Hopf bifurcation of the full system $(\ref{inter})$ occurs within an $O(\varepsilon)$ of the FSN II curves (c.f.  figure {\ref{two_par_bif_full}).   The FSN I (a) and FSN I (b)  curves represent saddle-node bifurcations  of folded singularities of  system (\ref{desing}). In the region enclosed by the FSN II (b) curve, $E^*$ exists as an ordinary saddle, and as an ordinary node outside this region.  Similarly $E_{xz}$ exists as an ordinary node inside the region bounded by the FSN II (a) curve and as an ordinary saddle  otherwise. There may exist up to four folded singularities in the region to the left of FSN I (a), exactly two folded singularities between FSN I (a) and FSN I (b), and no  folded singularities to the right of FSN I (b).  In region A, there exists a folded node, a folded focus and two folded saddles. Region B also contains four folded singularities, namely, two folded nodes, a folded saddle and a folded focus. Region C has a folded node, folded saddle and a folded focus. In region D, there exist a folded node and a folded focus, and region E contains a folded saddle and a folded focus.  In regions F and H, there exists a folded node/folded focus and a folded saddle, whereas region G contains either two folded nodes or a folded node and a folded focus. We will return to this bifurcation diagram in a later section.

   \begin{figure}[h!]     
  \centering 
{\includegraphics[width=10.7cm]{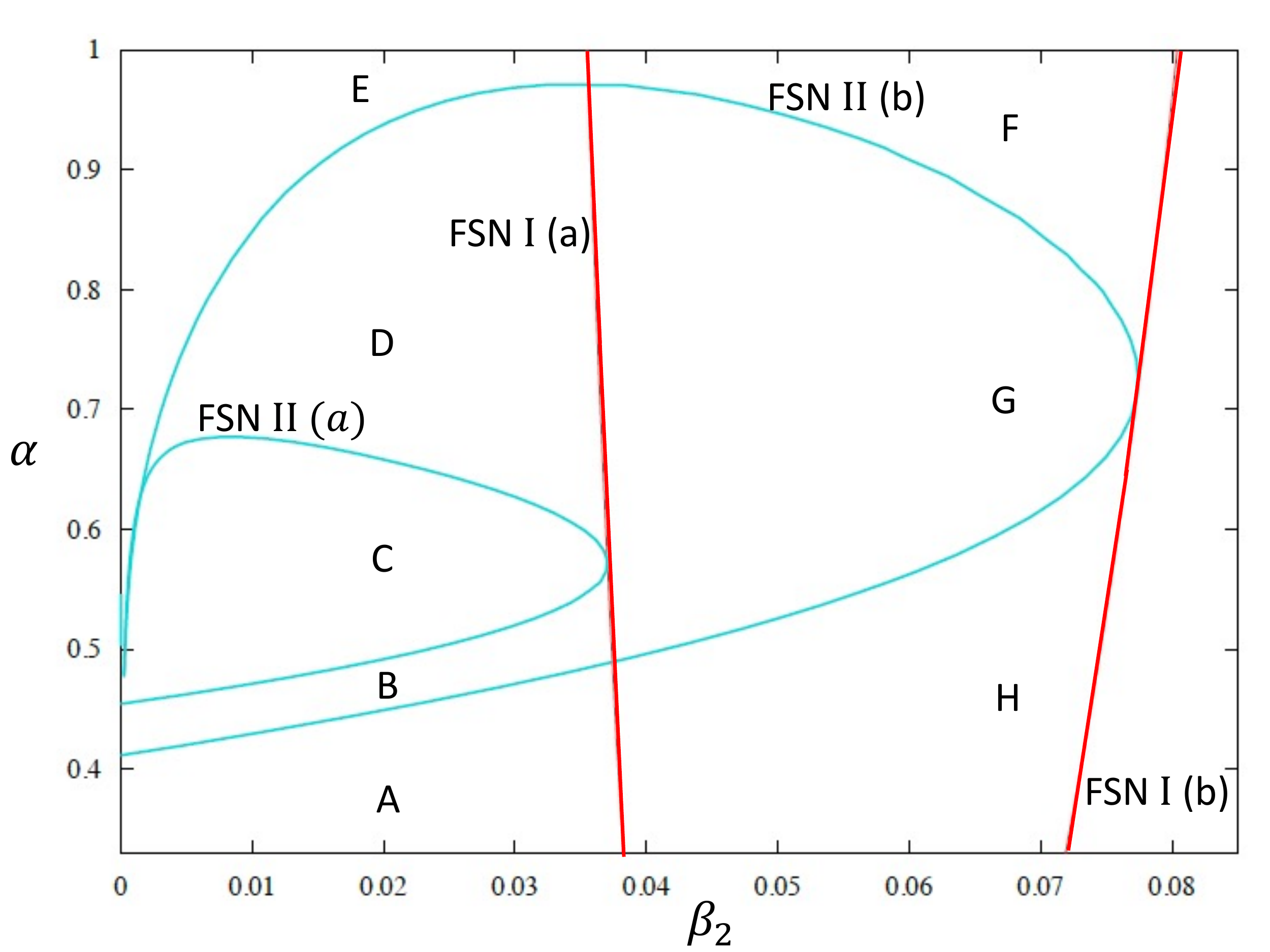}}
\caption{Two-parameter bifurcation structure of the desingularized system (\ref{desing}). The FSN I (a) and FSN I (b) curves represent folded saddle-node bifurcations of type I, while FSN II (a) and FSN II (b) curves represent folded saddle-node bifurcations of type II. See text for details.} 
\label{two_par_desing}%
\end{figure}

In the scenario when a folded node exists, the fold curve $\mathcal{F}$ and the \emph{strong singular canard} $\gamma_{0, \delta}$ form a trapping region ({\emph{singular funnel}}) on $S^a$ such that all solutions in the funnel converge to the folded node \cite{DGKKOW}.
 For certain choices of $\beta_2$ and $\alpha$, system (\ref{red2}) may admit two singular funnels, one located on $\mathcal{S}^{-a}$ and the other on $\mathcal{S}^{+a}$. A trajectory of the full system (\ref{inter}) may experience local oscillations near one or multiple branches of the fold curve $\mathcal{F}$, depending on which singular funnel it enters in the singular limit $\varepsilon \to 0$. Figure \ref{singular_funnel_proj_xz} (also see figure \ref{singular_funnel_proj_xz_1}) represents the scenario when an orbit enters in a close vicinity of the singular funnel on $\mathcal{S}^{-a}$, the  funnel being bounded by the middle branch  $\mathcal{F}^0$ of the fold curve and the strong singular canard  $\gamma_{0, \delta}^-$. The orbit passes close to the folded node singularity and makes rotations as it goes past the equilibrium $E^*$ of the full system (\ref{inter}).


   \subsection{Two-fast/one-slow analysis} 
 In this subsection, we will decompose system $(\ref{inter})$ into a one-dimensional slow subsystem with slow variable ($z$) and a two-dimensional fast subsystem with  fast variables ($x$ and $y$) by keeping $\varepsilon$ fixed and treating $\delta$ as the singular parameter.  We will analyze the bifurcation structure of the fast subsystem with the slow variable treated as a parameter. Key bifurcation structures such as the {\emph{superslow manifold}}, defined by the stationary solutions of the fast subsystem, periodic solutions, Hopf and homoclinic bifurcations of the fast subsystem will be useful in understanding bursting dynamics of the full system.
   
\subsubsection{The superslow manifold $\mathcal{M}_2$} 

Fixing $\varepsilon$ and letting $\delta \to 0$ in system (\ref{nondim2}) yields the fast subsystem

\begin{eqnarray}\label{layer1}    \left\{
\begin{array}{ll}  \dot{x} &=x\phi(x,y,z,\rho)\\
    \dot{y}&= \varepsilon y\chi(x,y,\rho)\\
    \dot{z}&=  0.
      \end{array} 
      \right. 
\end{eqnarray}
The equilibria of (\ref{layer1}) forms a one-dimensional critical manifold $\mathcal{M}_{2,\varepsilon}$, called the {\emph{superslow manifold}}, consists of three components, namely $\mathcal{M}_{2, \varepsilon} := \mathcal{K} \cup  \mathcal{L} \cup \mathcal{Z}$, where
\bess
\mathcal{K} = \{(0, 0, z): z\geq 0 \}, \ \mathcal{L} = \left\{(x, y, z) \in \mathbb{R}^3: y=0, z=H(x), x>0\right\}, \eess
and
\bess
  \mathcal{Z} &=& \{ (x, y, z) \in {\mathbb{R}^3}^+: \phi=0 \land \chi=0\}\\
  &=& \left\{ (x, y, z) \in S: z = G(x) \right\}
\eess
with $H(x)$, $G(x)$ defined by 
\bess
H(x) &=& \frac{1}{x}(1-x)(\beta^2_2 +x^2) \ \textnormal{and} \\
G(x) &=& \frac{1}{\alpha x}\left((1-x)(\beta^2_2 +x^2) -\frac{((1-\delta_1)x - \delta_1\beta_1)(\beta^2_2+x^2)}{\gamma (\beta_1+x)^2}\right).
\eess

Note that both $H(x)$ and $G(x)$ can have at most two relative extreme values on the interval $(0, 1)$. We further note that the geometry of  the superslow manifold $\mathcal{M}_{2,\varepsilon}$ does not depend of $\varepsilon$, though its stability does.  From here on we will suppress its dependence on $\varepsilon$, and denote the superslow manifold by $\mathcal{M}_{2}$. 
The sets $\mathcal{L}$ and $\mathcal{Z}$ are degenerate on $\mathcal{F}_{\mathcal{L}}$ and  $\mathcal{F}_{\mathcal{Z}}$ respectively, where
\bess
\mathcal{F}_{\mathcal{L}} &=& \{(x, y, z)\in  \mathcal{L}: H'(x)=0\} \\
& = & \{(x_{\mathcal{L}},0, H(x_{\mathcal{Z}})): x_{\mathcal{L}} \ \textnormal{is a positive root of}\ H'(x)=0 \}
\eess
and
\bess
\mathcal{F}_{\mathcal{Z}} &= &\{(x, y, z)\in \mathcal{Z}: G'(x)=0\}\\
& =& \{(x_{\mathcal{Z}}, F(x_{\mathcal{Z}}, G(x_{\mathcal{Z}})), G(x_{\mathcal{Z}})): x_{\mathcal{Z}} \ \textnormal{is a positive root of}\ G'(x)=0 \}.
\eess
These sets contain the isolated fold points of  $\mathcal{L}$ and $\mathcal{Z}$, and are referred to as  the ``knees" of these curves. In the scenario when $\mathcal{Z}$ has exactly two folds, $G(x)$ is cubic-shaped consisting of three branches, namely ${\mathcal{Z}}^-$,  $\mathcal{Z}^0$ and ${\mathcal{Z}}^+$ joined at the fold points. Denoting the $x$-coordinates of the folds of $\mathcal{Z}$ by $x^-_{\mathcal{Z}}$ and $x^+_{\mathcal{Z}}$, we then have that
\bess
{\mathcal{Z}}^- =\{(x, y, z)\in \mathcal{Z}: 0<x <x^-_{\mathcal{Z}}\},\  {\mathcal{Z}}^+ =\{(x, y, z)\in \mathcal{Z}: x^+_{\mathcal{Z}}<x <1\}
\eess
and
\bess
{\mathcal{Z}^0} =\{(x, y, z)\in \mathcal{Z}: x^-_{\mathcal{Z}} < x <x^+_{\mathcal{Z}}\}.
\eess
Similarly, the curve  $\mathcal{L}$ consists of three branches ${\mathcal{L}}^-$, ${\mathcal{L}}^0$ and ${\mathcal{L}}^+$ joined at the folds when $H(x)$ is cubic-shaped. The relative position of the folds of $\mathcal{L}$ and $\mathcal{Z}$ with respect to the fold curve $\mathcal{F}$  can play a crucial role in organizing the reduced flow on $\mathcal{M}_1$ and shaping the geometrical structure of mixed-mode oscillatory patterns \cite{kpk}. In this paper, we have considered parameter values such that $\mathcal{L}$ and $\mathcal{Z}$ are cubic-shaped and will discuss the role they play in organizing bursting phenomena in the full system.  

Linearization of (\ref{layer1}) around ${\mathcal{L}}^{\pm}$, ${\mathcal{L}}^0$, ${\mathcal{Z}}^{\pm}$ and $\mathcal{Z}^0$  determines the stability of these branches. Typically the middle branches ${\mathcal{L}}^0$ and $\mathcal{Z}^0$ consist of saddles of (\ref{layer1}). The other branches may lose normal hyperbolicity when system (\ref{layer1}) undergoes Hopf bifurcations. Denoting the set of Hopf points of (\ref{layer1}) by $\mathcal{M}_{DH}^{\varepsilon}$, we have that $\mathcal{M}_{DH}^{\varepsilon}  = {\mathcal{L}}_{DH}^{\varepsilon}  \cup \mathcal{Z}_{DH}^{\varepsilon}$, where
\bess
{\mathcal{L}}_{DH}^{\varepsilon} &=& \{(x, y, z) \in \mathcal{L}: x\phi_x + \varepsilon  \chi =0\}  \ 
\textnormal{and} \\
 \mathcal{Z}_{DH}^{\varepsilon} &=&  \{(x, y, z) \in \mathcal{Z}: x\phi_x + \varepsilon y \chi_y =0 \land  \varepsilon(\phi_x \chi_y - \phi_y \chi_x)>0\}.
\eess
The set ${\mathcal{L}}_{DH}^{\varepsilon}$ divides ${\mathcal{L}}^{\pm}$ into attracting and repelling branches, ${{\mathcal{L}}^{\pm}}^a$ and ${{\mathcal{L}}^{\pm}}^r$  respectively. Similarly, ${\mathcal{Z}}_{DH}^{\varepsilon}$ divides ${\mathcal{Z}}^{\pm}$ into its attracting and repelling branches, ${{\mathcal{Z}}^{\pm}}^a$  and ${{\mathcal{Z}}^{\pm}}^r$ respectively. In a neighborhood of ${\mathcal{Z}}_{DH}^{\varepsilon}$, ${{\mathcal{Z}}^{\pm}}^a$ consists of stable foci of (\ref{layer1}), while  ${{\mathcal{Z}}^{\pm}}^r$ consists of unstable foci of (\ref{layer1}).

Along $\mathcal{Z}$, the $x$-components of the degenerate nodal points of (\ref{layer1}) are roots of $\Lambda_{\varepsilon} (x)=0$, where 
\bess \Lambda_{\varepsilon}(x) &= &( x\phi_x(x, F(x, G(x)), G(x))-\varepsilon \gamma_1F(x, G(x)))^2 \\
 &-& 4\varepsilon x F(x, G(x))\left(-\gamma_1 \phi_x(x, F(x, G(x)), G(x)) +\frac{\beta_1}{(\beta_1+x)^3}\right).
\eess
We will consider roots of $\Lambda_{\varepsilon} (x)=0$ which are located in the interval $[0, 1]$. 
 The eigenvalues of the linearization of (\ref{layer1}) along $\mathcal{Z}$ are 
   \bes \label{eigtwofast} \lambda^{\pm}_{\varepsilon}(x)= \frac{1}{2} \left[ x\phi_x(x, F(x, G(x)), G(x))-\varepsilon \gamma_1F(x, G(x)) \pm \sqrt{\Lambda_{\varepsilon}(x)}\right].
 \ees
 The criticality of the Hopf bifurcation at the Hopf points $\mathcal{Z}_{DH}^{\varepsilon}$ is determined by the sign of $ \Delta |_{\mathcal{Z}_{DH}^{\varepsilon}}$, where
 \bess
 \Delta |_{\mathcal{Z}_{DH}^{\varepsilon}}: &=& \Big(\frac{-y \chi_y\sqrt{-xy \phi_y \chi_x}}{2y\chi_x (2\phi_x+x\phi_{xx})} - \frac{y \chi_y}{2\sqrt{-xy \phi_y \chi_x}} + \frac{(3\phi_{xx} +x \phi_{xxx})\sqrt{-xy \phi_y \chi_x}}{2(2\phi_x+x\phi_{xx})^2} \\
 &+& \frac{y\chi_x(\phi_y +x \phi_{xy})}{2 (2\phi_x+x\phi_{xx})\sqrt{-xy \phi_y \chi_x}} \Big)\Big|_{\mathcal{Z}_{DH}^{\varepsilon}}.
 \eess
 A subcritical (supercritical) Hopf bifurcation occurs when $ \Delta |_{\mathcal{Z}_{DH}^{\varepsilon}}>0 (<0)$ \cite{BE}.

   \begin{figure}[h!]     
  \centering 
{\includegraphics[width=9.5cm]{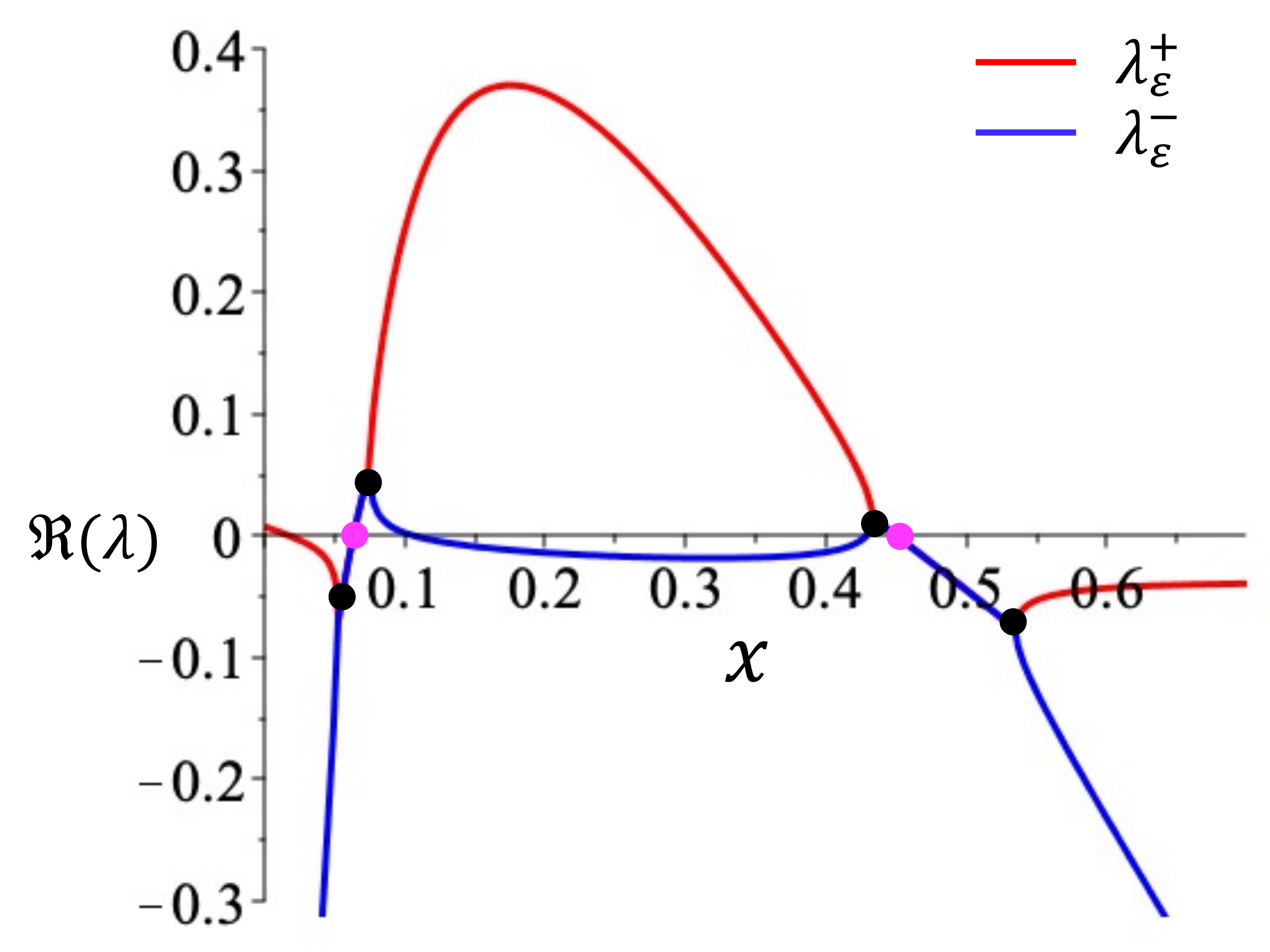}}
\caption{Real parts of eigenvalues of the layer problem (\ref{layer1})  for  $\beta_2=0.005$ and $\alpha=0.6$. The degenerate nodes are marked by black dots and the delayed Hopf bifurcation points by magenta dots.}
\label{realeigen}%
\end{figure}

Figure \ref{realeigen} shows the real parts of $\lambda^{\pm}_{\varepsilon}(x)$ with respect to $x$ for $\beta_2=0.005$ and $\alpha=0.6$. For these parameter values, there exists four degenerate nodal points $x_{DN1}^-< x_{DN2}^- < x_{DN1}^+ < x_{DN2}^+$ in $[0, 1]$. For $0<x<x_{DN1}^-$ and $ x_{DN2}^+<x<1$,  $\lambda^{+}_{\varepsilon}$  and  $\lambda^{-}_{\varepsilon}$ are the weak and strong stable eigenvalues respectively.  A branch switch occurs at $x_{DN2}^- $ with $\lambda^+_{\varepsilon}$ now  being strongly unstable and $\lambda^-_{\varepsilon}$ being weakly unstable or stable for $x_{DN2}^-<x<x_{DN1}^+$. System  (\ref{layer1}) also undergoes Hopf bifurcations twice; the $x$ components of the Hopf points ${DH}^{1,2} \in \mathcal{Z}_{DH}^{\varepsilon} $ are located in the intervals $(x_{DN1}^-, x_{DN2}^-)$ and $(x_{DN1}^+, x_{DN2}^+)$ with $\Delta |_{{DH}^1} = 0.145$ and $\Delta |_{{DH}^2} = 1.33$ indicating that both bifurcations are subcritical. We will refer to this diagram in figure \ref{delayed_hopf}.

The flow on $\mathcal{M}_2$ is governed by the superslow subystem 
\begin{eqnarray}\label{reduced2}    \left\{
\begin{array}{ll}  0 &=x\phi(x,y,z,\rho)\\
    0 &=  y\chi(x,y,\rho)\\
    \dot{z}&=  z \psi(x, z,\rho),
      \end{array} 
      \right. 
\end{eqnarray}
obtained by letting $\delta \to 0$ in system (\ref{slow}). The superslow flow is singular at the folds $\mathcal{F}_{\mathcal{L}} \cup \mathcal{F}_{\mathcal{Z}}$. Note that system (\ref{reduced2}) is singular at the fixed points $(1, 0, 0)$ and the coexistent state, $E^*$, of the full system. Canard solutions may arise when $E^*$ coincides with $\mathcal{F}_{\mathcal{Z}}$ and singular Hopf bifurcation occurs. This aspect will be not be discussed in this paper. The superslow flow occurs along $\mathcal{M}_2$ until it reaches a Hopf point $\mathcal{M}_{DH}^{\varepsilon}$. By Fenichel's theory, the normally hyperbolic segments of $\mathcal{M}_2$ perturb to a locally invariant slow manifold ${\mathcal{M}_2}_{\varepsilon, \delta} := {\mathcal{K}}_{\varepsilon, \delta} \cup {\mathcal{L}}_{\varepsilon, \delta} \cup {\mathcal{Z}}_{\varepsilon, \delta}$ for $\delta> 0$, and the flow restricted to these components are $O(\delta)$ perturbation of the superslow flow on $\mathcal{M}_2$.

The slow flow on  ${\mathcal{Z}}_{\varepsilon, \delta}^{\pm}$  can experience a delay in being repelled from ${\mathcal{Z}}_{\varepsilon, \delta}^{\pm r}$ after it goes past the Hopf bifurcation point $\mathcal{Z}_{DH}^{\varepsilon}$, as the accumulated contraction to  ${\mathcal{Z}}_{\varepsilon, \delta}^{\pm a}$ must get balanced by the total expansion from  ${\mathcal{Z}}_{\varepsilon, \delta}^{\pm r}$. Such a  mechanism of bifurcation delay is referred to as the {\emph{delayed loss of stability}} \cite{neish1, neish2} (also see \cite{KPK,LRV} for details).  

Figures \ref{delayed_hopf} and \ref{bursting} include bifurcation diagrams of the the fast subsystem (\ref{layer1}) for varying $\alpha$. In each case, the bifurcation diagram has the same qualitative features, namely an S-shaped curve of fixed points,  $\mathcal{Z}={\mathcal{Z}}^- \cup \mathcal{Z}^0 \cup {\mathcal{Z}}^+$, and two unstable branches of periodic orbits that emerge from subcritical Hopf bifurcations of (\ref{layer1}) located at ${\mathcal{Z}}^{\pm}$. These branches either terminate in homoclinic bifurcations (HC) with nearby saddle points or make  large excursions in phase plane  before returning to the stable manifold of the saddle.  The former type of homoclinic connection will be referred to as a ``small homoclinic loop" while the latter as a ``big homoclinic loop".  We will revisit the bifurcation structure of (\ref{layer1}) in the next section.

  \begin{figure}[h!]     
  \centering 
{\includegraphics[width=12.67cm]{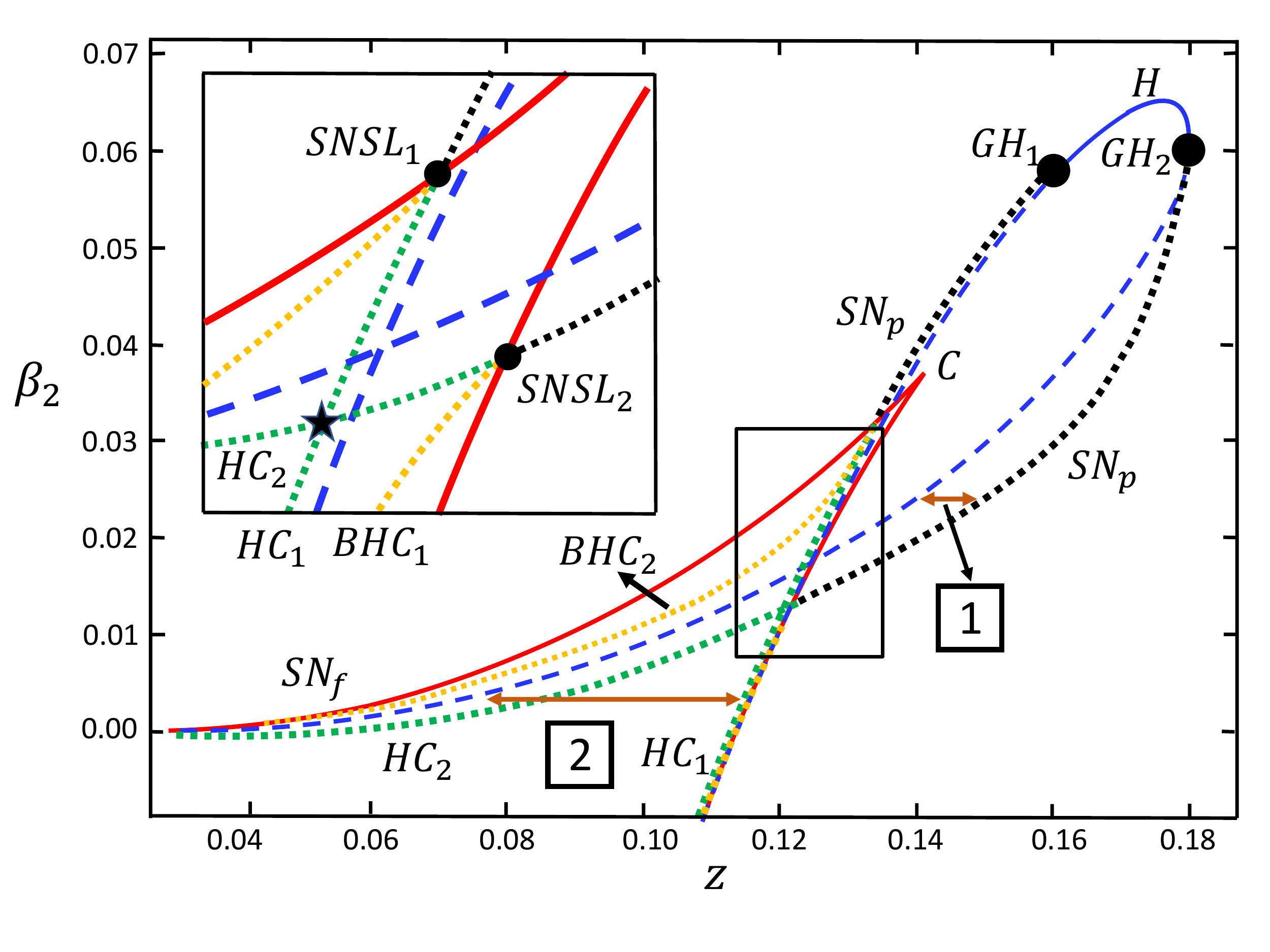}}
\caption{Two parameter bifurcation diagram of (\ref{layer1}) in $(z, \beta_2)$ parameter plane with  $\alpha=1$. The HC curves intersect with each other at the ``gluing bifurcation" denoted by a star and intersect with the $SN_f$ curve at SNSL points denoted by dots. The inset shows a qualitative representation of the region around the gluing bifuraction. $SN_f$ -  saddle-node bifurcation, $SN_p$- saddle-node of periodics, H - Hopf bifurcation, GH - generalized Hopf, HC - small homoclinic bifurcation, BHC - big homoclinic bifurcation, SNSL - saddle-node separatrix loop bifurcation.}
\label{two_par_fast}%
\end{figure}

\subsection{Two-parameter bifurcation of the fast subsystem (\ref{layer1})}
Figure \ref{two_par_fast} shows the changes in the bifurcation structure of the fast subsystem (\ref{layer1}) as $\beta_2$ is varied.  The bifurcation diagram was generated using XPPAUT. The qualitative features of the diagram remain the same (the diagram gets stretched to the right as $\alpha$ is decreased)  for all $\alpha \in (0, 1]$, and therefore we choose $\alpha=1$ as a representative. To this end, we compute loci of the codimension-1 bifurcations from figures \ref{delayed_hopf} and  \ref{bursting}  in the $(z, \beta_2)$ parameter plane of the fast subsystem. The curves of saddle-node bifurcations ($SN_f$) and Hopf bifurcations ($H$) are shown in red and blue respectively. The region enclosed by the $SN_f$ curve contains three equilibrium points, one of which is a saddle. The number of equilibria changes from  three to one upon crossing this curve. Codimension-2  bifurcations such as a cusp bifurcation ($C$) and a pair of Takens-Bogadnov bifurcations (not shown here; one of them occurs at $(z, \beta_2)=(0.023,-0.294)$ and the other at $(z, \beta_2)=(0.034, 0.0003)$) lie on the $SN_f$ curve where $H$ meets with the two branches of the $SN_f$ curve tangentially. Besides $C$, there exist several other noteworthy codimension-2 bifurcations. A pair of degenerate Hopf bifurcations denoted by $GH_1$ and $GH_2$, mark the points at which the Hopf bifurcation changes its criticality from subcritical (shown in dashed blue) to supercritical (in solid blue) and vice-versa. The Hopf curve also has a point of self-intersection where the two equilibria that are not saddle simultaneously undergo subcritical Hopf bifurcation. A pair of saddle-node of periodic orbits ($SN_p$) curves emerge from the degenerate Hopf points and terminate on the $SN_f$ curve at  codimension-2  saddle-node separatrix loop (SNSL) bifurcations. Inside the cusp region, a pair of homoclinic bifurcation curves, $HC_1$ and $HC_2$,  emanate from the Takens-Bogadnov bifurcation points and terminate at $SNSL_1$ and  $SNSL_2$ respectively on the $SN_f$ curve.  These homoclinic curves pertain to the ``small homoclinic loops".  The $HC_1$ and $HC_2$ curves intersect with each other at a codimension-2 bifurcation, referred to as a ``gluing bifurcation" \cite{GTW} where the stable and unstable manifolds of the saddle form a figure-eight.  In addition, there exist two more homoclinic bifurcation curves $BHC_1$ and $BHC_2$  corresponding to the ``big homoclinic loops"  closely following $HC_1$ and $HC_2$ respectively. The ``big homoclinic loops" encircle the two spiral coexistence equilibira, whereas the ``small homoclinic loops" encircle only one of the coexistence equilibrium points (see \cite{ELS} for an illustration of the homoclinic loops).

The inset in figure \ref{two_par_fast} is a qualitative representation of the bifurcation region around the ``gluing bifurcation". It shows the relative position of the small and big homoclinic bifurcation curves along with location of the bifurcation curves of saddle-node of periodic orbits. The precise location of these curves is very difficult to compute. The homoclinic  curves  $BHC_1$ and $HC_1$  occur in a very close vicinity of each other and to the Hopf curve $H$, so we present a qualitative depiction. Each of the curves $BHC_1$  and $BHC_2$ should also terminate at a pair of SNSL bifurcations. The region between these curves and the  Hopf curve is very narrow to locate the precise parameter values of the SNSL bifurcations.

System (\ref{layer1}) exhibits bistability between the equilibria $\mathcal{Z}^-$ and $\mathcal{Z}^+$ in the region enclosed by the subcritical Hopf curves, labeled as \fbox{2} in figure \ref{two_par_fast}. As the fast subsystem transitions from $\mathcal{Z}^-$ to $\mathcal{Z}^+$, system (\ref{nondim2}) exhibits MMOs of subHopf/subHopf type as shown in  figure \ref{delayed_hopf} (these dynamics correspond to MMO orbits with SAOs along $\mathcal{F}^0$ and $\mathcal{F}^+$ in figure \ref{two_par_bif_full}). Similarly, in the regime between the $SN_p$ curve  and the subHopf  curve labeled as \fbox{1}  in  figure \ref{two_par_fast}, a limit cycle attractor and a point attractor coexists. As the fast subsystem (\ref{layer1}) transitions from one attractor to the other, system (\ref{nondim2})   exhibits subHopf/fold cycle bursting as shown in figure \ref{bursting}(d). These dynamics correspond to subcritcal elliptic bursting patterns in figure \ref{two_par_bif_full}.  Other types of bursting dynamics as classified in \cite{izh} are also possible in system (\ref{nondim2}) but are beyond the scope of this paper.  We also remark that figure \ref{two_par_fast} qualitatively resembles the two-parameter bifurcation diagrams of the layer problems of the slow-fast models studied in  \cite{DZL} (Fig. 5) and \cite{GTW} (Fig. 1) to some extent.

\section{Analysis of the full model}

We recall that system (\ref{nondim2}) has three timescales with $x$ being the fast variable, $y$ and $z$ being the intermediate and slow variables respectively. The fastest timescale dominates the evolution of a trajectory unless it is near the critical manifold $\mathcal{M}_1$ or the superslow manifold $\mathcal{M}_2$, where the slower timescales come into effect. Taking the double limit $(\varepsilon, \delta) \to (0, 0)$ in system (\ref{nondim2}) yields the fast subsystem

\begin{eqnarray}\label{layersub}    \left\{
\begin{array}{ll}  \dot{x} &=x\phi(x,y,z,\rho)\\
    \dot{y}&= 0\\
    \dot{z}&=  0,
      \end{array} 
      \right. 
\end{eqnarray}
where the intermediate and slow variables $(y, z)$ are frozen. This system is precisely the layer problem (\ref{layer1d}) we studied in the 1-fast/2-slow approach.  Taking the double limit in the intermediate timescale gives the $1$D intermediate subsystem
\begin{eqnarray}\label{intersub}    \left\{
\begin{array}{ll}  0 &=x\phi(x,y,z,\rho)\\
  \dot{y}&= y \chi(x, y, \rho)\\
    \dot{z}&=  0.
      \end{array} 
      \right. 
\end{eqnarray}
In this case, the flow is governed by the intermediate variable $y$, restricted to the plane $\Pi$ or the surface $S$, and the slow variable $z$ remains the same. The fast variable $x$ immediately responds to changes in state before the intermediate flow takes over. The trajectories of system (\ref{intersub}) are referred to as intermediate fibers. The intermediate flow is not defined on the fold curve $\mathcal{F}$.  Finally, the slow subsystem is obtained by taking the double limit of the slow system (\ref{slow}) and is identical to system (\ref{reduced2}). 

The critical manifold $\mathcal{M}_1$ remains the set of equilibria of the fast subsystem (\ref{layersub}) and serves as the phase space of  the intermediate subsystem (\ref{intersub}). Similar to desingularizing the reduced system (\ref{reduced}), desingularization of (\ref{intersub}) describes the flow on $S$ 
\begin{eqnarray} \label{desingfull}  \begin{pmatrix}
 \dot{x}  \\
 \dot{z} 
 \end{pmatrix}= \begin{pmatrix}
\phi_yy \chi  \\
0
 \end{pmatrix}\bigg|_{y=F(x,z)}.
\end{eqnarray}
The folded singularities of (\ref{desingfull})  is the set $S_{F}^0$ of isolated points that lie at the intersection of the superslow curves $\mathcal{Z}$ or $\mathcal{L}$ and the fold curve $\mathcal{F}$, i.e. 
\[S_{F}^0= \{(x, y, z)\in \mathcal{F}: \chi (x, y)=0\}\cup \{(x, y, z)\in \mathcal{F}: y=0\}. \]  On the other hand, the ordinary singularities in the double singular limit are the set of points
\[S_{E}^0= \{(x, y, z)\in S\setminus \mathcal{F}: \chi (x, y)=0\}\cup \{(x, y, z)\in S \setminus \mathcal{F}: y=0\} \] 
formed by the superslow curves $\mathcal{Z}$ and $\mathcal{L}$ off the fold curve $\mathcal{F}$. Note that the ordinary singularities $S_E^{\delta}$ of (\ref{desing}) do not persist as singularities for system (\ref{desingfull}). This is due to the fact that the constraints $z=0$ or $\psi=0$ are no longer required to hold for existence of ordinary singularities in (\ref{desingfull}). The weak eigenvalue of a folded singularity of (\ref{desingfull}) is zero, which then implies that the folded singularity is a FSN in the double singular limit. As $\varepsilon \to 0$, the set of Hopf points $\mathcal{M}_{DH}^{\varepsilon}$ on the superslow manifold $\mathcal{M}_2$ satisfy $\phi_x=0$, and thus merge with the set of folded singularities $S_{F}^0$. 

The slow subsystem  (\ref{reduced2}) approximates the slow flow of system (\ref{inter}) for sufficiently small $\delta>0$ under the assumption that the variables $x$ and $y$ change much rapidly than $z$ and thereby quickly approach their steady states under small changes in $z$. The superslow manifold $\mathcal{M}_2$, which is the set of equilibria of (\ref{desingfull}) is the phase space of  (\ref{reduced2}). The equilibrium points of the full system $E^*$ and the boundary equilibrium $E_{xy}$ are the only equilibria of (\ref{reduced2}). The knees of $\mathcal{Z}$ are singular points of the slow flow.

By standard GSPT \cite{CT}, the normally hyperbolic portions $S^{\pm a}$ of $\mathcal{M}_1$ and 
${\mathcal{Z}}^{\pm}$ of $\mathcal{M}_2$ perturb to  
$S^{\pm a}_{\varepsilon, \delta}$, and ${\mathcal{Z}}^{\pm}_{\varepsilon, \delta}$ respectively. For a trajectory that starts on $S^{+ a}_{\varepsilon, \delta}$ (say), will follow the intermediate flow on it until it gets attracted to ${\mathcal{Z}}^{+}_{\varepsilon, \delta}$ or reaches a vicinity of $\mathcal{F}^+$. In the former case, it follows the superslow flow on  ${\mathcal{Z}}^{+}_{\varepsilon, \delta}$ and can experience a delay of loss of stability resulting in SAOs, while in the latter case, it jumps to the opposite attracting branch of the slow manifold $S^{- a}_{\varepsilon, \delta}$ resulting in a large amplitude oscillation. 

System (\ref{inter}) exhibits a variety of complex oscillatory patterns, including, but not limited to, mixed-mode oscillations (MMOs) and bursting as seen in figures \ref{timeseries_example_1}-\ref{timeseries_example_2}.  The small-amplitude oscillations in an MMO orbit may occur near one of the three branches of the fold curve $\mathcal{F}$. An MMO orbit can pass very close to a folded node  of (\ref{desing}) as well as a  Hopf bifurcation point of the fast subsystem (\ref{layer1}), and therefore the SAOs in an MMO orbit can be organized by the canard dynamics arising from a folded node singularity, typically referred to as sector type dynamics \cite{DKP, KPK}, or the slow passage effect associated with a Hopf point, referred to as delayed-Hopf type dynamics \cite{DGKKOW}. In most cases (see figures \ref{timeseries_example_1} and \ref{timeseries_varying_alpha}), it turns out that the amplitude of the SAOs in the MMOs orbits of system (\ref{inter}) are exponentially small, and can be  associated with unstable limit cycles born at subcritical Hopf bifurcations of (\ref{layer1}) leading to transient oscillations.

In the next subsection, we will explore the bifurcation structure of system (\ref{nondim2}) by treating the predation efficiency $\beta_2$ of the generalist predator as the primary bifurcation parameter and the fraction $\alpha$ of the generalist predator's diet that consists of $x$ to study the different parameter regimes in which interesting ecological phenomena occur.

\subsection{One-parameter bifurcation}

  \begin{figure}[h!]     
  \centering 
{\includegraphics[width=12.67cm]{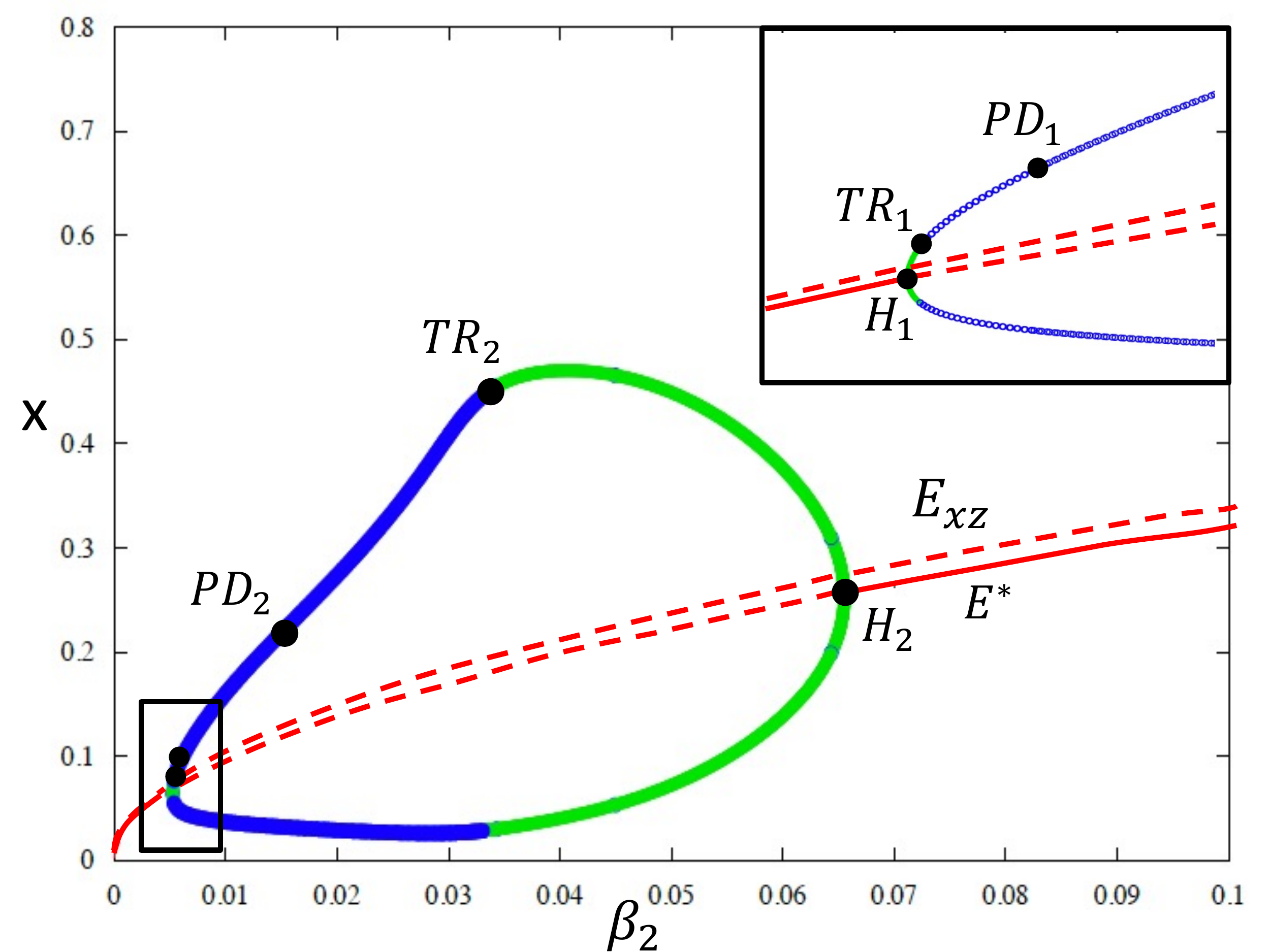}}
\caption{ One parameter bifurcation diagram of (\ref{inter}) with respect to $\beta_2$ with other parameters as in (\ref{parvalues}) and $\alpha=0.75$. Filled green (open blue) circles represent maximum and minimum values of $x$ in stable (unstable) limit cycles. H - Hopf bifurcation, PD- period-doubling bifurcation, TR - torus bifurcation}
\label{one_par_bif}%
\end{figure}
Using XPPAUT, a one-parameter bifurcation diagram was computed as shown in Figure \ref{one_par_bif}, where $\beta_2$ is the continuation parameter and the maximum norm is considered along the vertical axis. We first note that the boundary equilibrium state $E_{xz}$ exists as an unstable node or focus for all $\beta_2 \in (0, 1)$. For $\beta_2$ sufficiently small, the coexistent equilibrium $E^*$ exists as a stable attractor. It loses its stability at a supercritical Hopf bifurcation $H_1 \approx 0.00524$, giving birth to a family of stable periodic orbits. This  family of orbits loses its stability at a torus bifurcation $TR_1 \approx 0.00536$ giving way to MMOs and bursting oscillations. The time profiles of the solutions emerging from $TR_1$ display amplitude-modulated oscillations as shown in figure \ref{torus_canards1}(A). These solutions are headless mixed-type torus canards  \cite{DBKK} as they follow the repelling branch of limit cycles created at a subcritical Hopf bifurcation of the fast subsystem (\ref{layer1}) (not shown here). At $\beta_2$= 0.00536, we note solutions with very different oscillation profiles elucidating the presence of multiple timescales as shown in figure \ref{torus_canards1}(B).  The long quiescent phase of the solution  in figure \ref{torus_canards1}(B) is organized by a homoclinic bifurcation of the fast subsystem (\ref{layer1}). The trajectory spends a long time near a homoclinic orbit of (\ref{layer1}) in the superslow timescale and then follows the unstable branch of the limit cycles of (\ref{layer1}) and eventually jumps to an attracting sheet of the slow manifold (c.f. figure \ref{delayed_hopf}(D)).  Orbits with similar patterns were referred to as mixed-type torus canards with head in \cite{DBKK}.  A detailed study of these solutions is left for a future study.

        \begin{figure}[h!]     
  \centering 
\subfloat[]{\includegraphics[width=7.67cm]{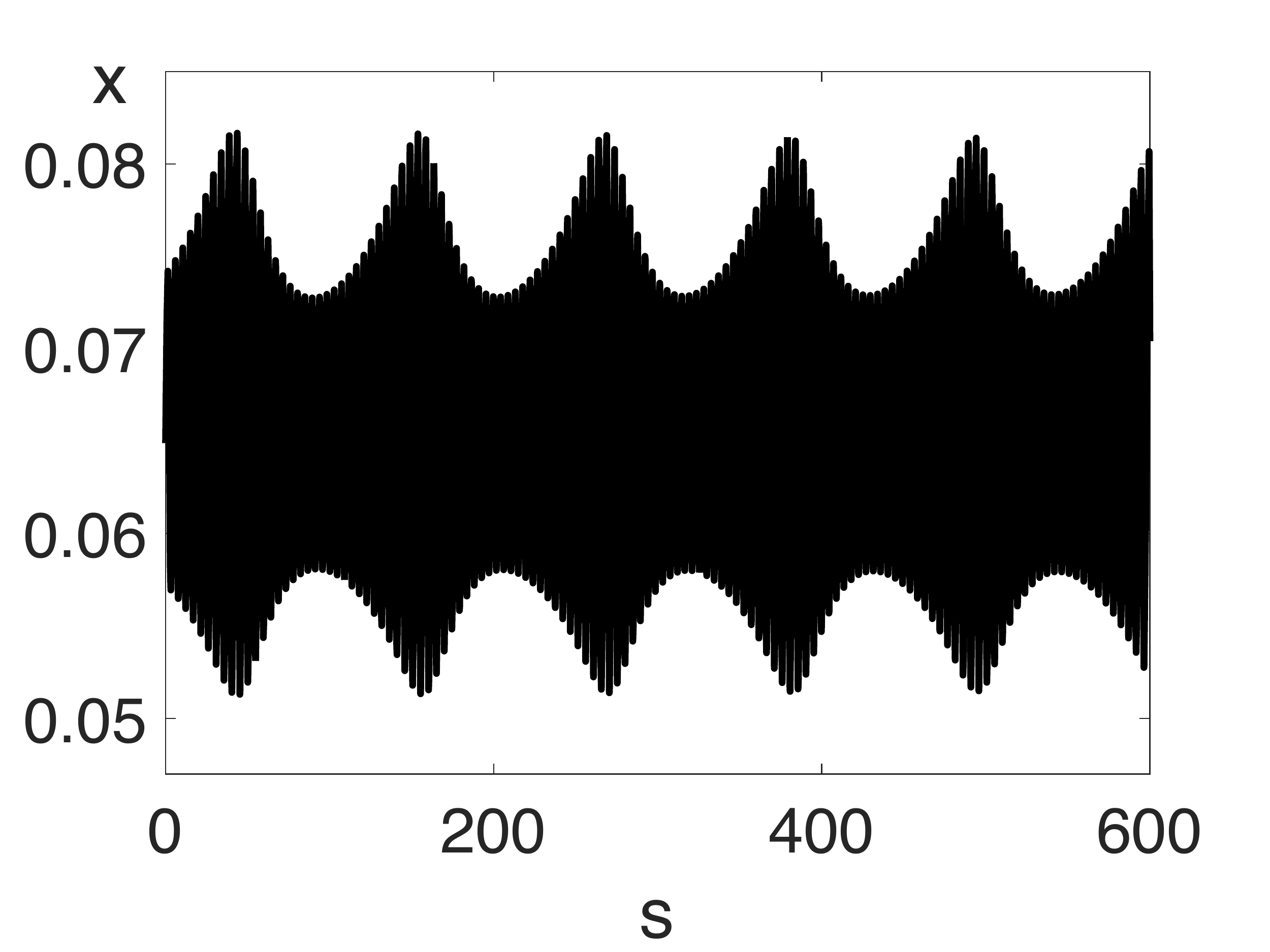}}
\quad
\subfloat[]{\includegraphics[width=7.67cm]{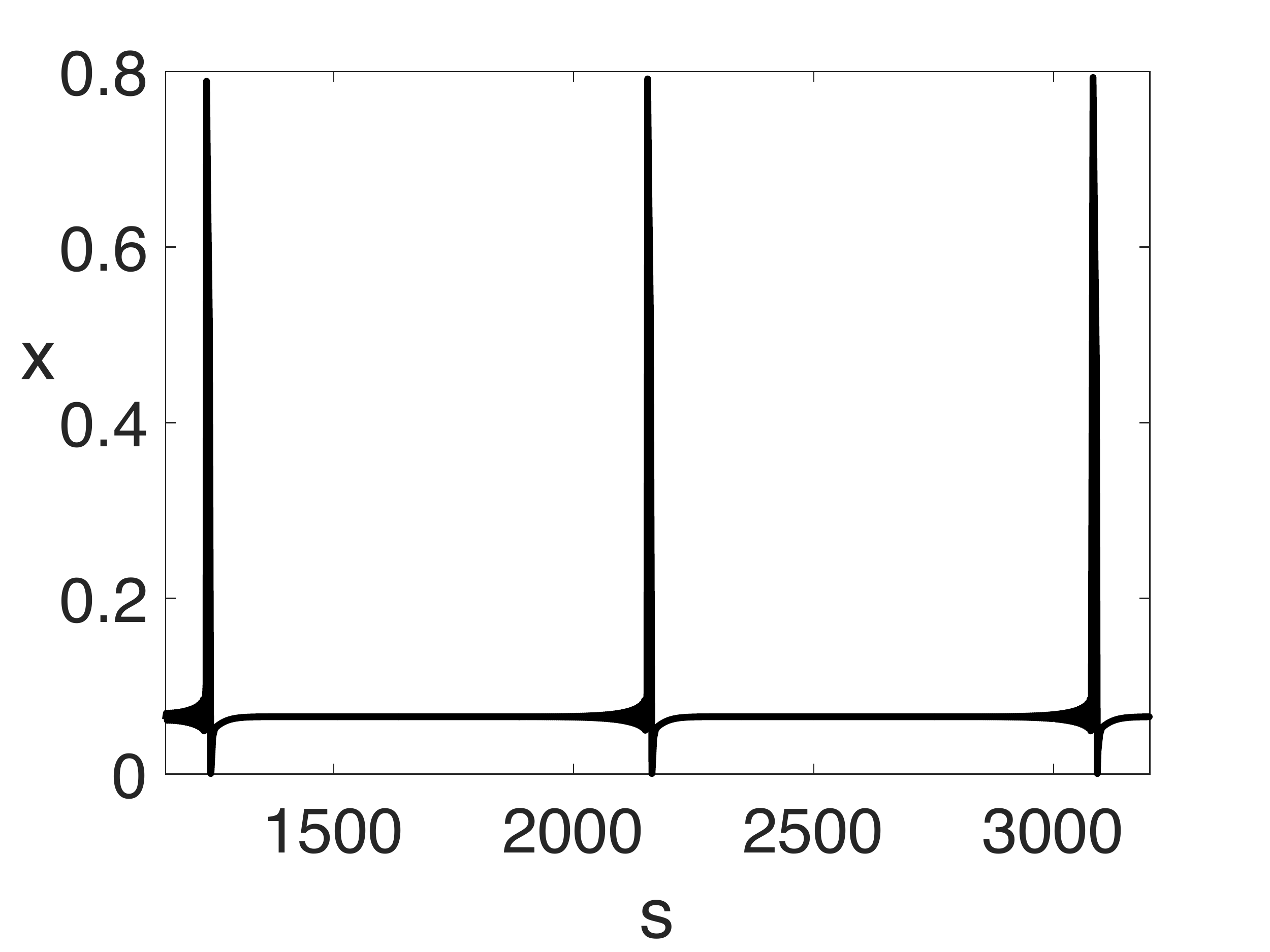}}

\caption{Time profiles of the $x$-components of solutions of system (\ref{inter}) for   $\alpha=0.75$ and other parameter values as in (\ref{parvalues}). (A) $\beta_2 =0.0053$. (B) $\beta_2 =0.00536$.}
\label{torus_canards1}
\end{figure}

 On further increasing $\beta_2$, the system undergoes a period-doubling bifurcation $PD_1 \approx 0.0064$, and MMOs are observed thereafter as shown in figure \ref{timeseries_example_1}. The MMOs persist until the system undergoes another period-doubling bifurcation $PD_2 \approx 0.0155$, after which bursting oscillations of sub-Hopf/fold cycle type ({\em subcritical elliptic bursting}) are observed (c.f figure \ref{bursting}(D)).   On further increasing $\beta_2$,  the system exhibits torus canards \cite{DBKK}, where the oscillations are qualitatively sub-Hopf/fold cycle type, except that the oscillations do not terminate at a saddle-node bifurcation of the periodics of the fast-subsystem, but instead continue along the branch of unstable limit cycles. The system undergoes another torus bifurcation $TR_2 \approx 0.0331$  giving rise to amplitude-modulated spiking orbits. The torus canards exist in a very small parameter regime that lies in a vicinity of $TR_2 \approx 0.0331$. 
After $TR_2$, the system exhibits spiking. The branch of spiking orbits persist until another  supercritcial Hopf bifurcation ($H_2\approx 0.066$) occurs, after which the coexistent state $E^*$ gains its stability and exists as a stable attractor.
   
\subsection{Two-parameter bifurcation}

We are interested in transitions between different dynamical regimes consisting of  spiking, bursting  and other types of oscillatory patterns in system (\ref{inter})  as the predation efficiency of the generalist predator and the fraction of the generalist predator's diet that consists of $x$ are varied. Figure \ref{two_par_bif_full} shows the two-parameter bifurcation structure of (\ref{inter})  in $(\beta_2, \alpha)$ space.  
Numerical continuation of Hopf bifurcations $H_1$ and $H_2$ in figure \ref{one_par_bif} generates a curve, denoted by HB, that divides the parameter space into two regions based on the stability of the steady state solution $E^*$.  
The equilibrium $E^*$ is stable outside the region bounded by HB and unstable otherwise.   The region enclosed by HB is further divided into several sub-regions consisting of different kinds of oscillatory dynamics such as mixed-mode oscillations, sub-elliptic bursting, large-amplitude oscillations, and sub-threshold or Hopf cycles. The curves obtained by continuing the torus bifurcation $TR_2$, and the period-doubled bifurcation $PD_2$ in figure \ref{one_par_bif} mark the boundaries between transitions from one type of dynamics to another; $TR_2$ separates sub-elliptic bursting from large-amplitude oscillations while $PD_2$ separates MMOs with single spikes from sub-elliptic bursting.

The parameter regime bounded by the torus curve, TR, consists of MMOs, subcritical elliptic bursting, and classical and mixed-type torus canard solutions. The mixed-type torus canards as well as classical torus canards occur in a very close vicinity of the boundary TR and appear during transition from sub-elliptic bursting to spiking, similar to the dynamics exhibited by the model considered in \cite{BAD}. The former type of dynamics are observed for lower values of $\beta_2$ while the latter for relatively higher values of $\beta_2$. The red dashed curve  in the region bounded by TR separates the one-spike periodic solutions from subcritical elliptic bursting with multiple spikes.  Many spike adding bifurcations occur in a very small parameter regime, and each time a spike is added, the periodic orbit transforms from a subcritical elliptic bursting orbit  from $n$ spikes to $n+1$ spikes. The precise mechanism for spike adding is left for future study. The SAOs in these periodic orbits occur near the middle branch,  $\mathcal{F}^0$,  of the fold curve. 
 In the regime with a single spike, the quiescent phase of the dynamics may persist for a prolonged time (see second - last panels of  figure \ref{timeseries_varying_alpha}). This regime is ecologically significant as it reveals the role of a generalist predator in regulating the population of prey.  A highly efficient generalist predator (i.e. for smaller values of $\beta_2$) can keep the focal prey density at a low level if the prey consists of a major part of its diet. On the other hand, in the regime consisting of elliptic bursting, we note that multiple spikes can occur in the prey dynamics even if it consists of a major part of the generalist predator's diet as shown in the last panel of figure \ref{timeseries_varying_alpha_1}.

As the efficiency of the generalist predator decreases (i.e. $\beta_2$ increases), the system could either exhibit MMOs with long epochs of SAOs, where the SAOs in the MMO orbits occur near $\mathcal{F}^+$ (see first panel of figure \ref{timeseries_varying_alpha_1}),  large-amplitude oscillations (see second and third panels of figure \ref{timeseries_varying_alpha_1}) or subcritical elliptic bursting patterns  (as shown in the last panel of figure \ref{timeseries_varying_alpha_1}) on varying $\alpha$. The MMO orbits with SAOs near $\mathcal{F}^+$ occur in a very narrow parameter regime close to the HB curve. 
Due to stiffness issues,  $PD_1$ and $TR_1$ could not be numerically continued. Hence, to obtain the boundary of the parameter regime separating the MMO orbits from the large-amplitude oscillations, similar to figure \ref{one_par_bif}, another one-parameter bifurcation diagram of system (\ref{inter}) with $\alpha=0.45$ is considered (not shown here). In this scenario as  $\beta_2$ is decreased,  $E^*$ experiences a supercritical Hopf bifurcation and gives birth to a family of stable limit cycles. The family loses its  stability at a torus bifurcation and the system exhibits MMOs of the type seen in the first panel of  figure  \ref{timeseries_varying_alpha_1}. The system then undergoes a PD bifurcation and transitions to relaxation oscillations. A numerical continuation of the PD bifurcation gave rise to the boundary 
separating the MMO orbits from the large-amplitude oscillations  in figure  \ref{two_par_bif_full}.
The remaining portion of the region enclosed by HB in figure  \ref{two_par_bif_full}  consists of either sub-threshold/Hopf cycles or relaxation type oscillations. The sub-threshold oscillations occur in a very close vicinity of HB, whereas  relaxation oscillations occur in the transition regime from MMOs to sub-elliptic bursting.

 \subsection{Analysis of three-timescale solutions.}
 
   \begin{figure}[h!]     
  \centering 
{\includegraphics[width=12.67cm]{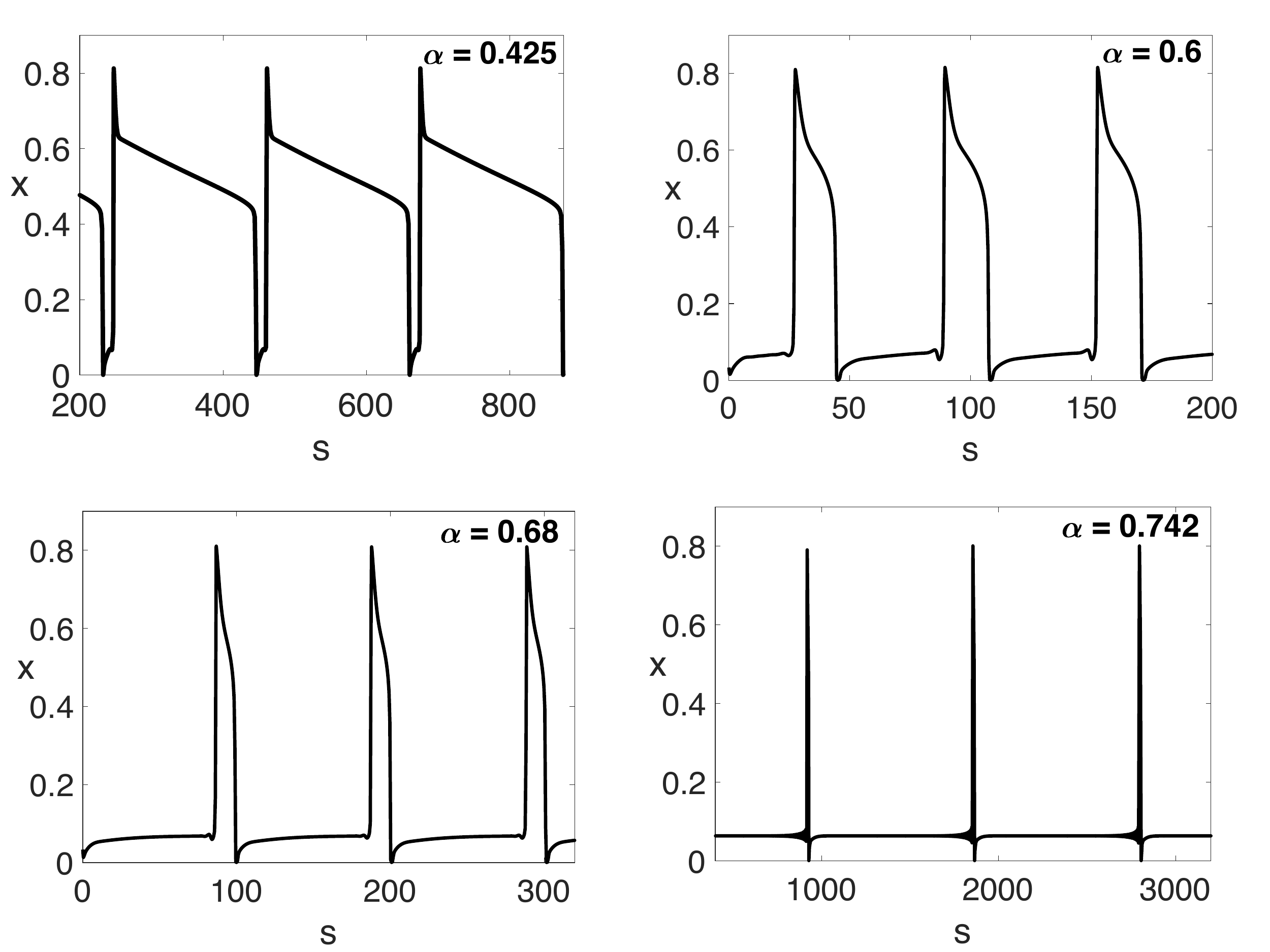}}
\caption{Time series for the $x$-coordinate of the orbits of system (\ref{inter}) for $\beta_2=0.005$ and varying $\alpha$. The other parameter values are as in (\ref{parvalues}). Note the change in profile of the MMO patterns.}
\label{timeseries_varying_alpha}%
\end{figure}

 In this subsection we focus on the roles of the critical manifold $\mathcal{M}_1$ and the superslow manifold $\mathcal{M}_2$ in organizing the flow of system (\ref{inter}).  We will see that the flows associated with the reduced problem (\ref{reduced}) and the fast subsystem (\ref{layer1}) give us good insight of the mechanisms responsible for organizing MMOs and bursting dynamics in the full system. We consider parameter regimes in figure \ref{two_par_bif_full} in which system (\ref{inter})  either exhibits MMOs with varying epochs of SAOs as shown in figure \ref{timeseries_varying_alpha} or transitions from MMOs to subcritical elliptic bursting patterns as shown in figure \ref{timeseries_varying_alpha_1}. 
  
   \begin{figure}[h!]     
  \centering 
\subfloat[]{\includegraphics[width=7.85cm]{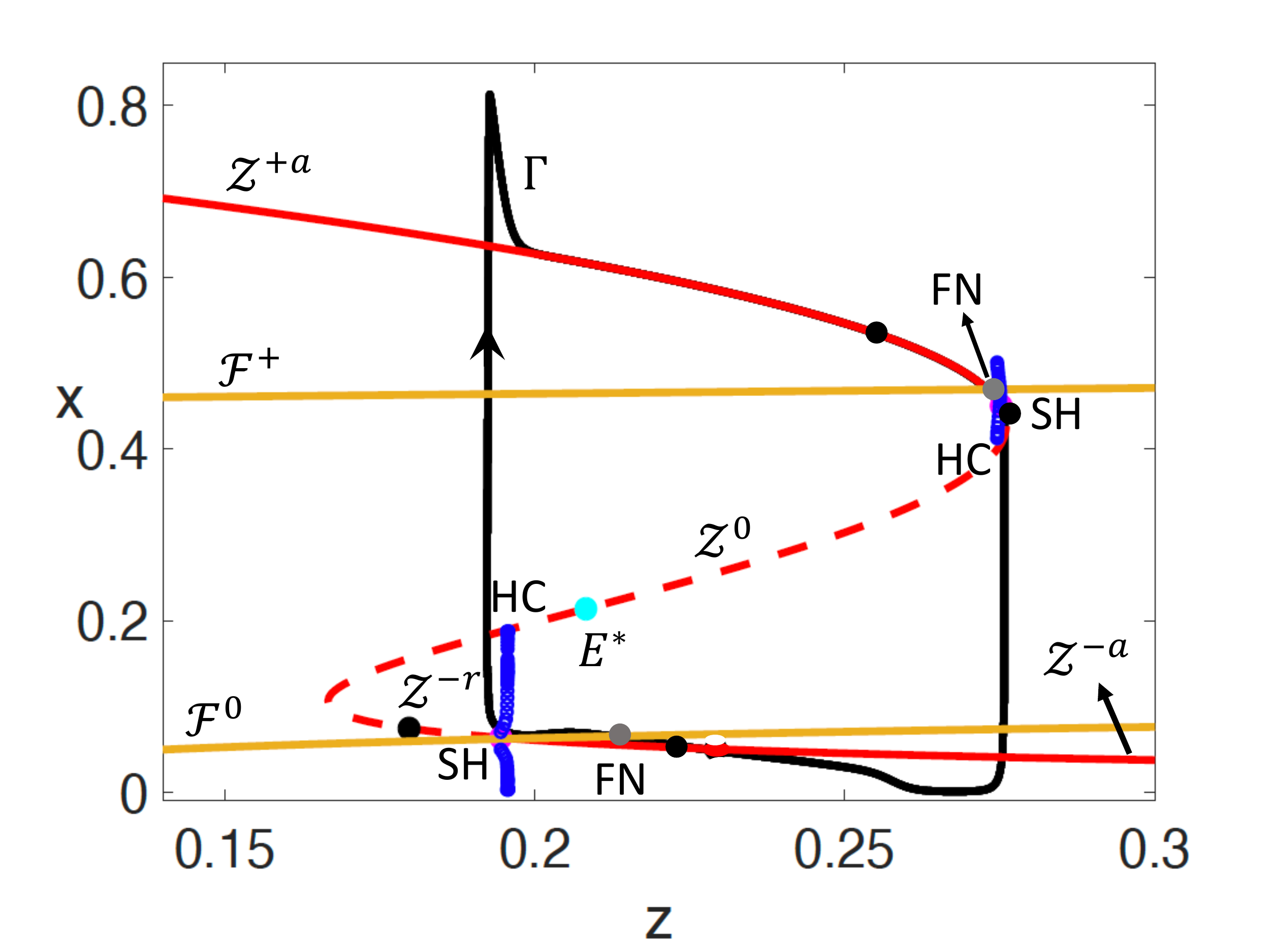}}
\quad
\subfloat[]{\includegraphics[width=7.85cm]{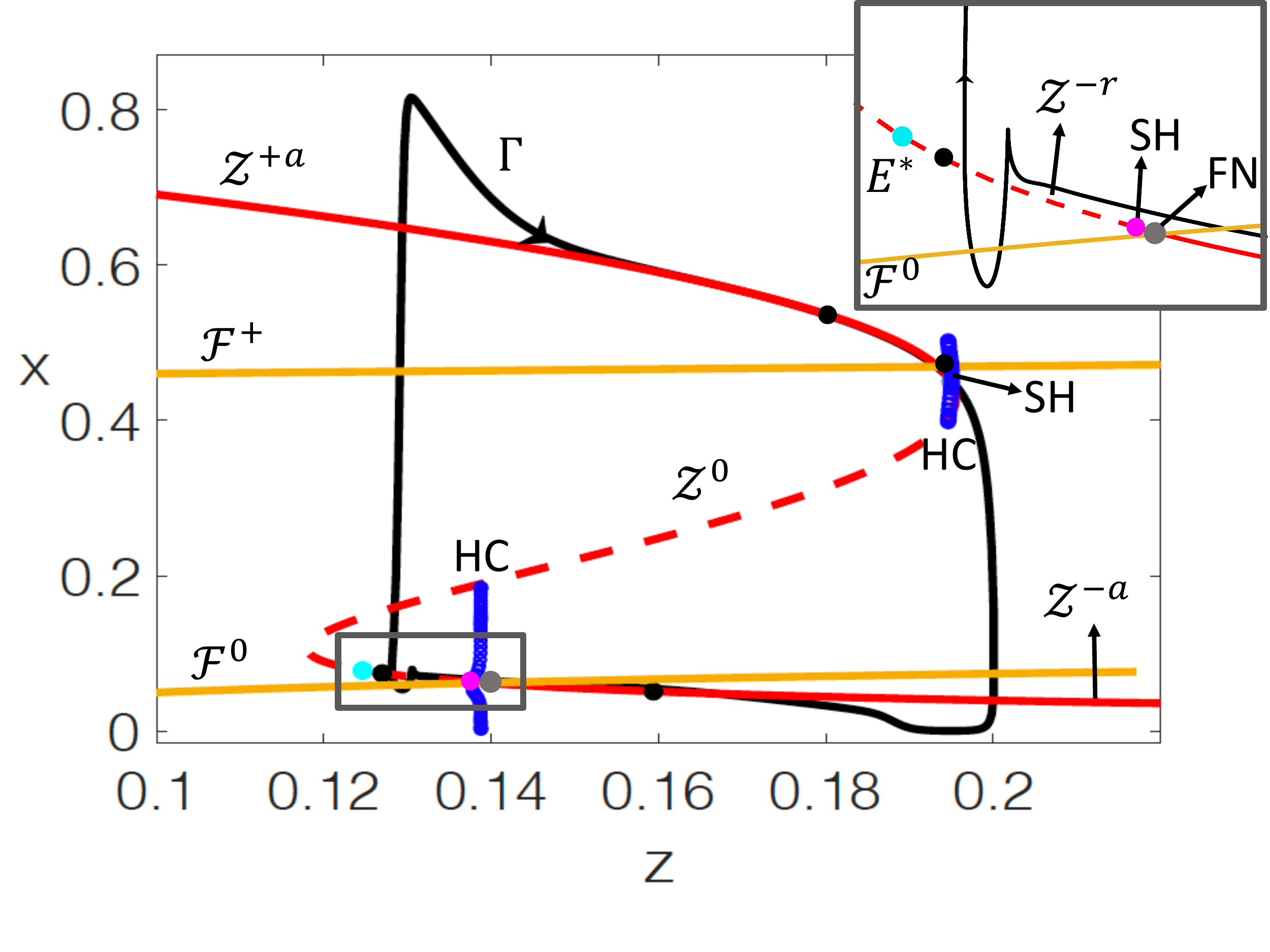}}
\quad
\subfloat[]{\includegraphics[width=7.85cm]{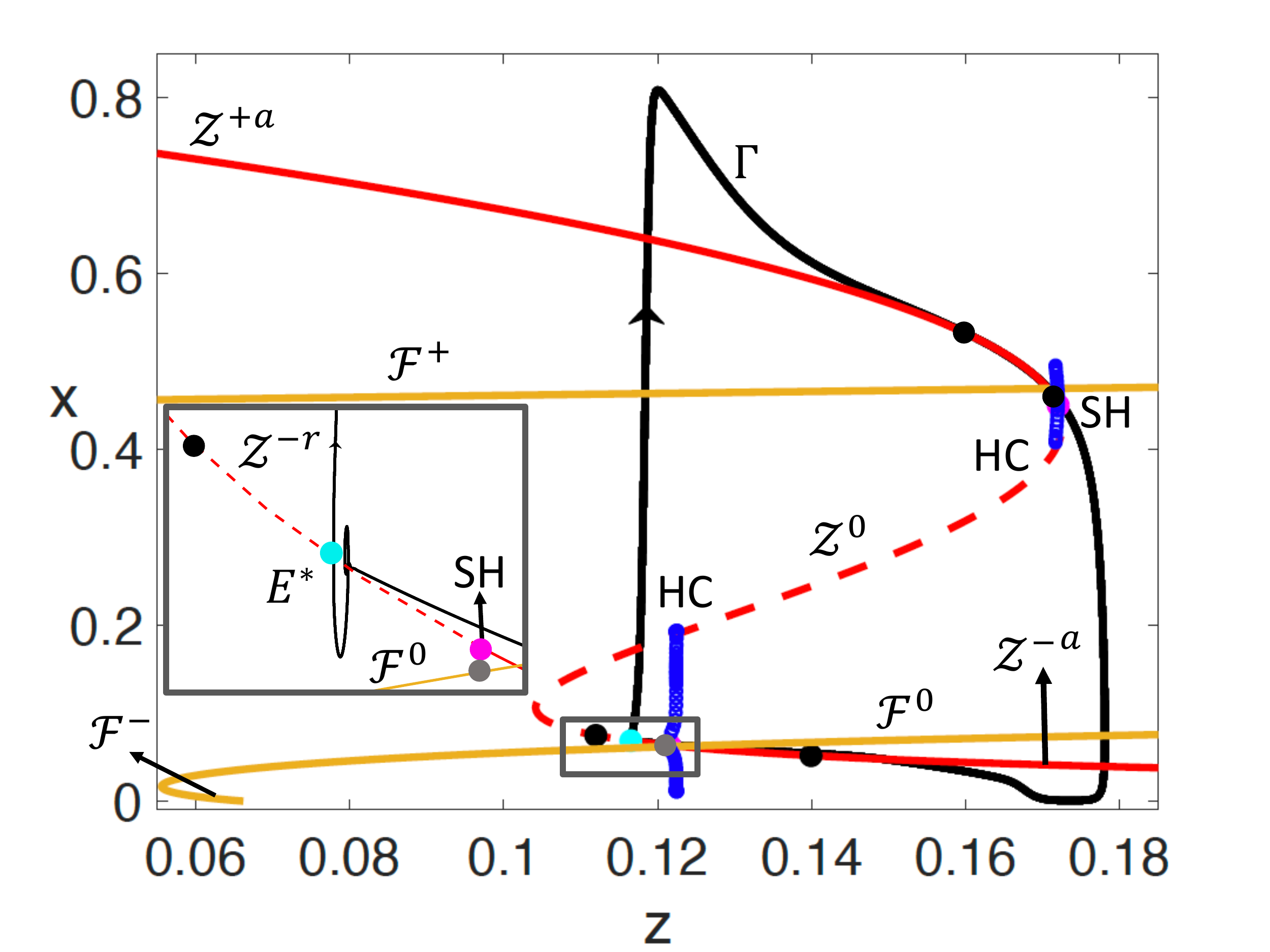}}
\quad
\subfloat[]{\includegraphics[width=7.85cm]{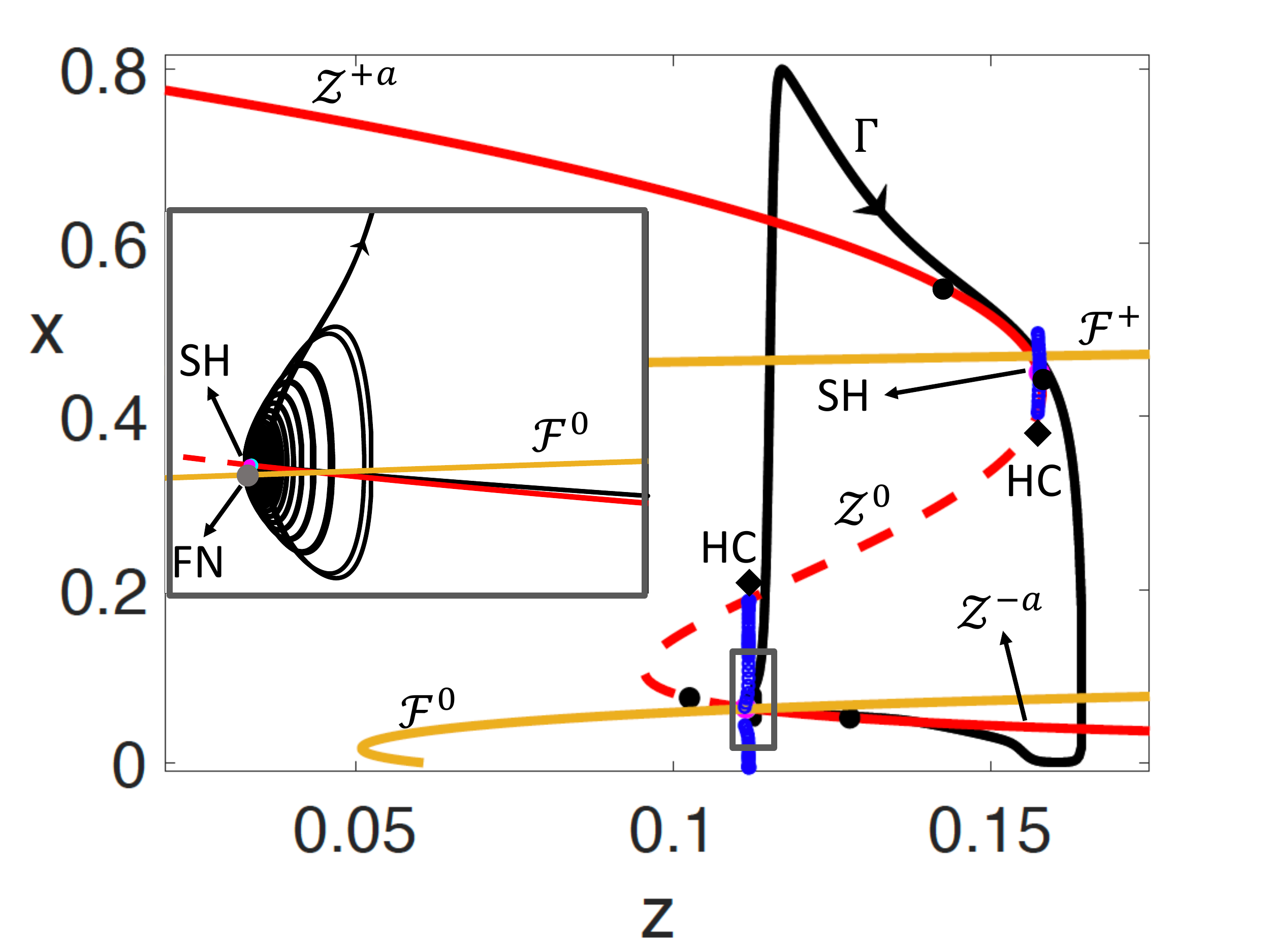}}
\caption{(A) Overlay of the projection of trajectory $\Gamma$ of (\ref{inter}) onto the $(z, x)$ - plane with the bifurcation diagram of the fast subsystem (\ref{layer1})  for $\beta_2=0.005$. (A) $\alpha=0.425$  (B) $\alpha=0.6$ (C) $\alpha=0.68$ (D) $\alpha=0.742$.  Shown are the fold curves $\mathcal{F}^{\pm}, \mathcal{F}^0$ (yellow), the z-curves $\mathcal{Z}^{\pm}, \mathcal{Z}^{0}$ (red), folded node (grey dot),  equilibrium point $E^*$ (cyan dot), delayed Hopf points (magenta dots) and degenerate nodes (black dots) of (\ref{layer1}).  Note that the trajectory exhibits delayed Hopf bifurcation on the lower branch of $\mathcal{Z}$ in panels (B)-(D).  FN - folded node,  SH - subcritical Hopf, SNP - saddle-node bifurcation of periodic orbit, HC -  homoclinic bifurcation, $\mathcal{Z}$ -  the superslow curve. Dashed lines denote instability. Open blue circles represent maximum and minimum values of $x$ in unstable limit cycles.}
\label{delayed_hopf}%
\end{figure}

 To this end, we superimpose the trajectories corresponding to the time profiles in figure \ref{timeseries_varying_alpha} projected on the $(z, x)$ - phase plane  on the bifurcation diagram of the fast subsystem  (\ref{layer1}) as shown in figure  \ref{delayed_hopf}. We also include the projection of the fold curve $\mathcal{F}$. We note that for the choice of the parameter values  in figure  \ref{delayed_hopf}, $\mathcal{F} = \mathcal{F}^+ \cup \mathcal{F}^0\cup \mathcal{F}^-$ is piecewise continuous  (belongs to case 3 in figure \ref{bif_foldcurve}), 
and the critical manifold $S$ may have up to four sheets. However, the fast dynamics of the MMO orbits are directed towards or away from that part of $S$  which consists of exactly three normally hyperbolic sheets. The superslow manifold consists of two attracting branches, two repelling branches, and a saddle branch, denoted by ${{\mathcal{Z}}^{\pm}}^a$, ${{\mathcal{Z}}^{\pm}}^r$ and ${{\mathcal{Z}}^{0}}$ respectively. For $\varepsilon$ and $\delta $ sufficiently small, the fast component of a trajectory of  (\ref{inter}), which is a perturbation of the dynamics governed by (\ref{layersub}), first
  brings the trajectory  towards the slow manifold $S^{+a}_{\varepsilon, \delta}$ where it initially overshoots ${\mathcal{Z}}^{+ a}$. The intermediate flow which is now a perturbation of (\ref{intersub}) governs the dynamics on  $S^{+a}_{\varepsilon, \delta}$ and brings the trajectory to the perturbed superslow manifold ${\mathcal{Z}^{+a}_{\varepsilon, \delta}}$ where the superslow flow takes over. The superslow flow, which is a perturbation of the dynamics  governed by (\ref{reduced2}), slowly takes the trajectory past $\mathcal{F}^+$  while reaching a vicinity of a homoclinic bifurcation point HC of (\ref{layer1}) on ${\mathcal{Z}}^{+ a}$. It eventually reaches a neighborhood of the fold $\mathcal{F}_{\mathcal{Z}}$,  and jumps to $\Pi^a_{\varepsilon, \delta}$, where its flow is governed by (\ref{intersub}). As the orbit descends along this manifold, it goes past the transcritical bifurcation $TC$ and stays on  $\Pi^r_{\varepsilon, \delta}$ for a while before it concatenates with a fast fiber, resulting in Pontryagin's delay of  loss of stability \cite{AS, kpnew, Sadhudcds}. The orbit then jumps to  the  attracting branch $S^{-a}_{\varepsilon, \delta}$ of the slow manifold and gets attracted towards ${{\mathcal{Z}^{-a}_{\varepsilon, \delta}}}$. As it follows ${{\mathcal{Z}^{-a}_{\varepsilon, \delta}}}$, it slowly passes through a neighborhood of a canard point until it reaches a neighborhood of the Hopf point  $\mathcal{Z}_{DH}^{\varepsilon}$ (denoted by SH in  in figure  \ref{delayed_hopf}) on ${{\mathcal{Z}}^{- a}}$.  After reaching $\mathcal{Z}_{DH}^{\varepsilon}$, the trajectory does not immediately jump to the  opposite attracting branch ${{\mathcal{Z}^{+a}_{\varepsilon, \delta}}}$ of the superslow manifold, but continues to drift close to ${\mathcal{Z}^{- r}}$, tracing $\mathcal{Z}^{-r}_{\varepsilon, \delta}$. 
The trajectory may remain close to ${{\mathcal{Z}}^{- r}}$ for an $O(\delta)$ distance past the Hopf bifurcation as seen in the insets in figure \ref{delayed_hopf}(B)-(C). The small amplitude oscillations  are below a visible threshold and indistinguishable from the superslow manifold $\mathcal{Z}$. However, the size of the oscillations grow  as the equilibrium of the full system approaches the delayed Hopf point. In figure \ref{delayed_hopf}(D), the trajectory reaches a vicinity of a homoclinic orbit of (\ref{layer1}) and spends a prolonged time near $E^*$ and SH before it jumps to $S^a_{\varepsilon, \delta}$. The trajectory exhibits MMOs, where 
the large amplitude oscillation in the MMO orbits can be viewed as a hysteresis loop that alternately jumps 
between the two subcritical Hopf points on the two branches of $\mathcal{Z}^{\pm}$, and the small amplitude oscillations are guided by a slow passage through SH, further influenced by the unstable manifold of the equilibrium $E^*$. According to the classification in \cite{izh}, the dynamics in figure \ref{delayed_hopf} can be referred to as subHopf/subHopf and can be mapped to the region labeled as \fbox{2} in figure \ref{two_par_fast}. Note the degenerate nodes and the subHopf points of the fast subsystem corresponding to figure  \ref{delayed_hopf}(b) were also shown in figure \ref{realeigen}.

      \begin{figure}[h!]     
  \centering 
\subfloat[]{\includegraphics[width=7.7cm]{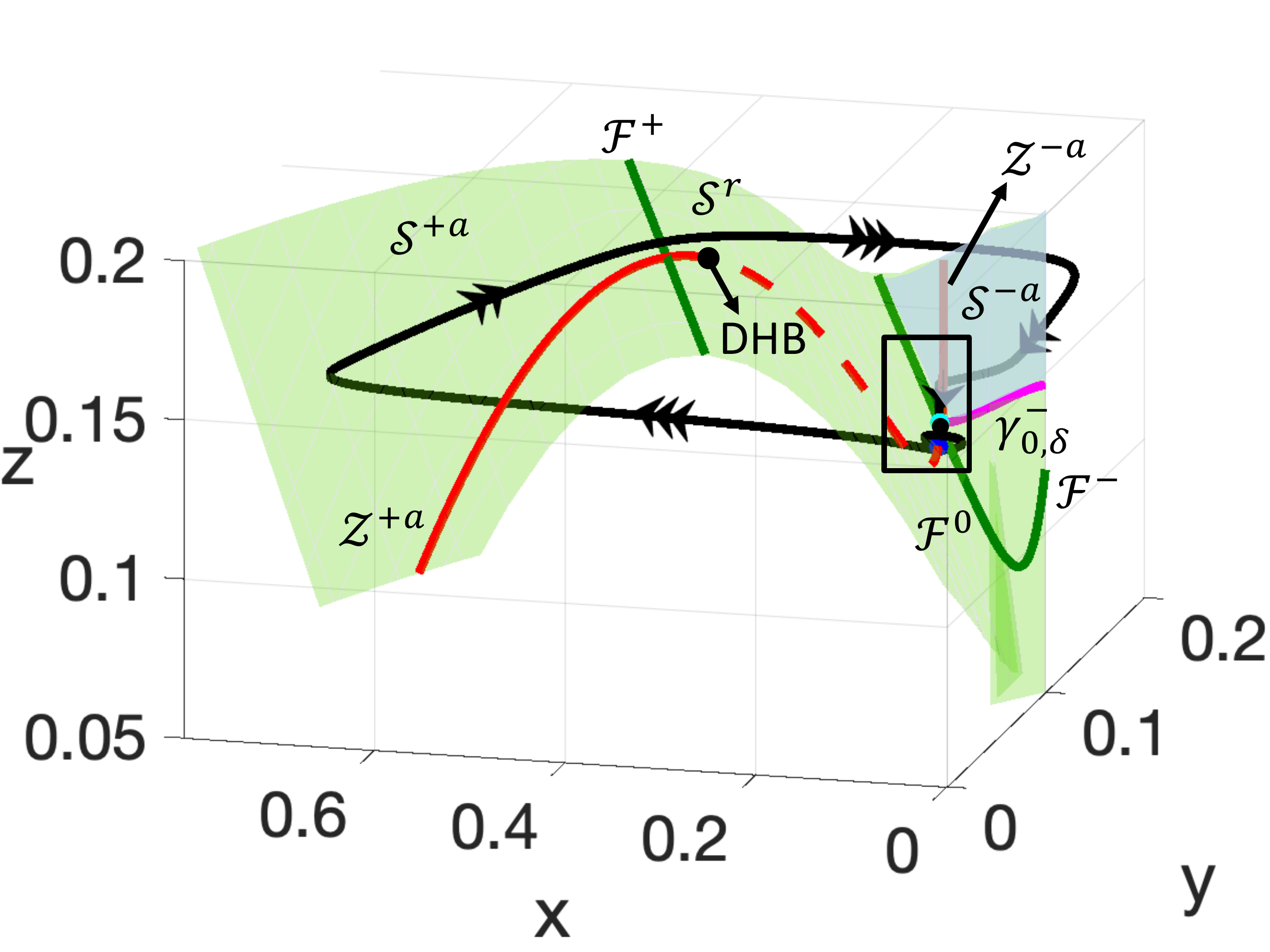}}
\quad
\subfloat[]{\includegraphics[width=7.7cm]{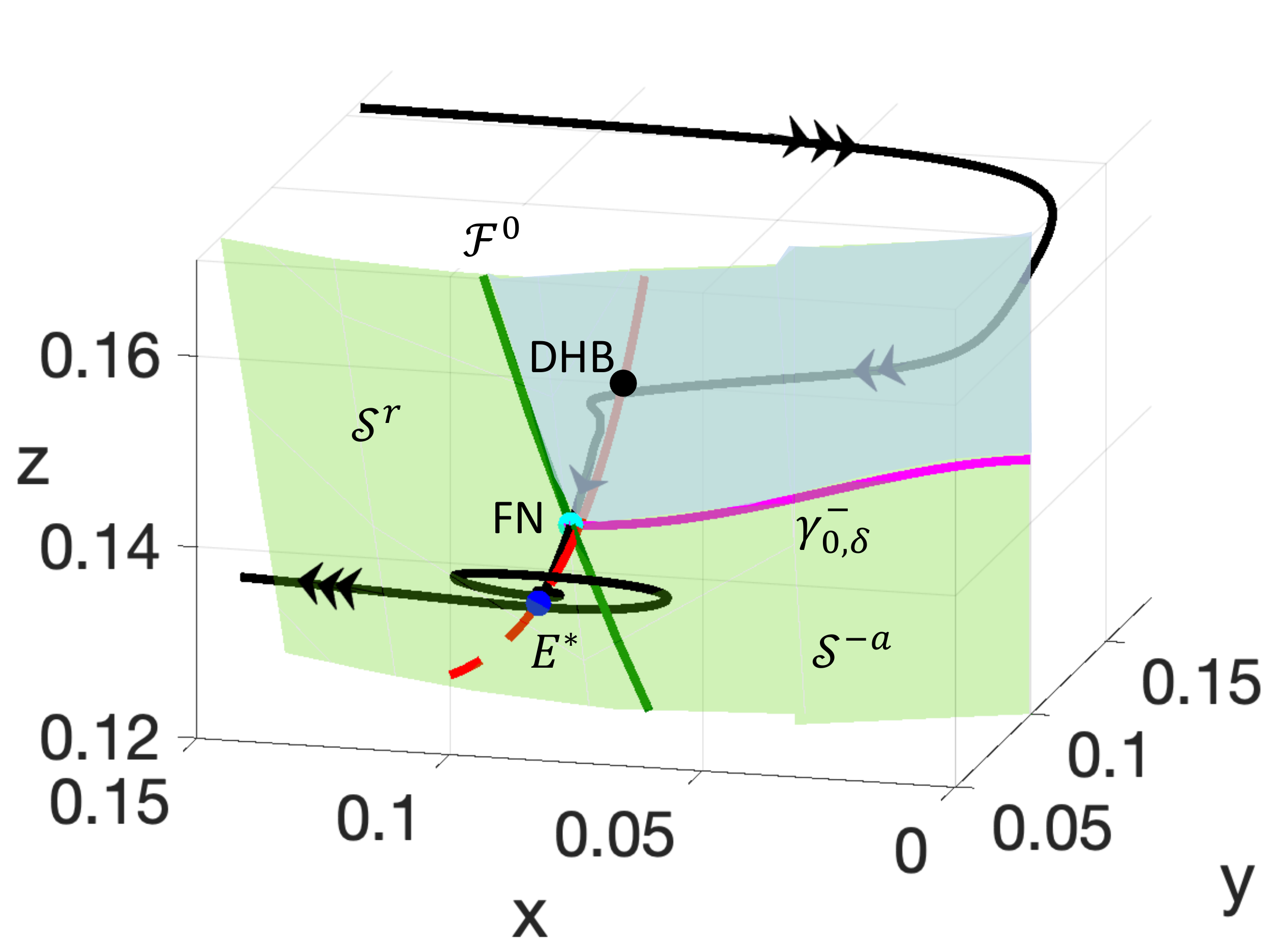}}
\caption{(A) Mixed-mode oscillations exhibited by system (\ref{inter}) for parameter values $\beta_2=0.01$ and $\alpha=0.75$. (B) Zoomed view of the dynamics near the lower fold $\mathcal{F}^0$ on the $xz$-plane. The trajectory enters into the singular funnel on $\mathcal{S}^{-a}$ shown by the shaded region and filters through the the folded-node singularity (FN) while passing close to the delayed-Hopf bifurcation point (DHB). The local vector field of the equilibrium $E^*$ further influences its dynamics before it jumps to $S^{+a}$ .}
\label{singular_funnel_proj_xz}%
\end{figure}

The parameter values chosen in figure \ref{delayed_hopf} belong to regions B, C or D in figure \ref{two_par_desing}, thus indicating the existence of at least one folded node singularity in the system. We note that in each panel of figure \ref{delayed_hopf}, the delayed Hopf point lies in a close vicinity of the folded node singularity on the lower branch $\mathcal{Z}^{-a}$ of the superslow manifold. It is not clear how the folded node singularity influences the dynamics of an MMO orbit while it makes a slow passage through the Hopf point, $\mathcal{Z}_{DH}^{\varepsilon}$. The folded node  approaches the Hopf point as $\alpha$ increases, and they both approach the equilibrium point $E^*$. The equilibrium $E^*$ is a saddle-focus with a two-dimensional unstable manifold. The vector field around $E^*$ also plays an important role in generating the SAOs in the MMO orbits as the trajectories pass closely to $E^*$ as seen in figure \ref{delayed_hopf}(B)-(D). 

To gain a better perspective, we consider the phase portrait of an MMO orbit and examine its relative position with respect to $E^*$, the canard point, the delayed Hopf point, the critical manifold $S$ and the superslow manifold $\mathcal{Z}$ as shown in figure \ref{singular_funnel_proj_xz}. The phase space in figure \ref{singular_funnel_proj_xz} qualitatively represents the dynamics of the orbits in figure \ref{delayed_hopf}(B)-(D) (cf. figure \ref{crit_mfld}(A)). 
The SAOs in these MMO orbits are observed  near the branch $\mathcal{F}^0$ of the fold curve (also see figure \ref{timeseries_example_1}).  By the canard theory, for local oscillations to occur, the trajectory must land in a neighborhood of the singular funnel  in one of the attracting sheets of the slow manifold and rotate around the primary weak canard during its passage through a folded node singularity (see \cite{BKW, TTVWB, Wesc}). In figure \ref{singular_funnel_proj_xz}, the trajectory lands in $S^{-a}_{\varepsilon, \delta}$ to one side of the strong canard $\gamma^{-}_{\varepsilon, \delta}$ and filters through the folded node while staying close to ${\mathcal{Z}^{-a}}$ during its passage. The primary weak canard (not shown here) lies close to the superslow manifold ${\mathcal{Z}^{-a}_{\varepsilon, \delta}}$, and in the singular limit merges with ${\mathcal{Z}^{-a}}$. In principle, one would need to draw the slow manifolds ${{S^r_{\varepsilon, \delta}}}$ and ${{S^{-a}_{\varepsilon, \delta}}}$  to study the locally twisted geometry of the intersection of these manifolds, which is beyond the scope of this paper. However, we note that  the SAOs generated are below a visible threshold with exponentially small amplitudes, and thus do not seem to be canard induced. 
In fact, the delayed Hopf point (DHB) lying in a close vicinity of the folded node singularity (see figure \ref{singular_funnel_proj_xz}(B)) is playing a crucial role in shaping the dynamics. Furthermore, the oscillations are initiated when the trajectory approaches $E^*$, which suggests that the unstable manifold of $E^*$ is also playing a role in  organizing the dynamics.

  \begin{figure}[h!]     
  \centering 
{\includegraphics[width=12.67cm]{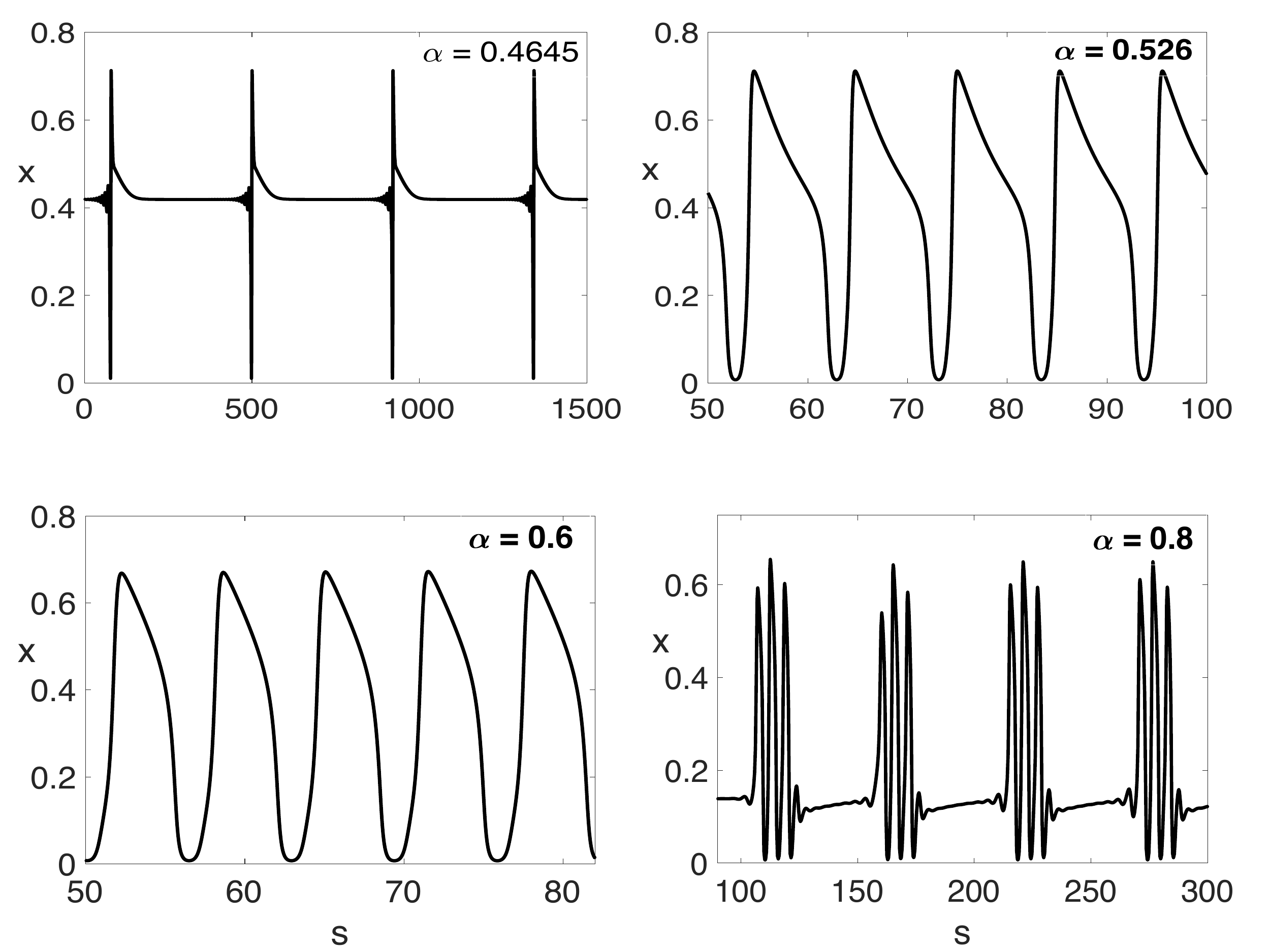}}
\caption{ Time series for the $x$-coordinate of the orbits of system (\ref{inter}) for $\beta_2=0.0245$ and varying $\alpha$. Note the transition from MMO patterns to bursting oscillations.}
\label{timeseries_varying_alpha_1}%
\end{figure}

Next, we consider the  parameter regime in which system (\ref{inter})  transitions from MMOs to bursting dynamics as shown in figure \ref{timeseries_varying_alpha_1}. In this regime, $\mathcal{F} = \mathcal{F}^+ \cup \mathcal{F}^0\cup \mathcal{F}^-$ is also piecewise continuous  (belongs to case 2(ii) in figure \ref{bif_foldcurve}),  
and the critical manifold $S$ may have up to four sheets.  However, similar to the parameter regime considered in figure \ref{delayed_hopf}, it turns out that the fast fibers of the orbits are directed towards or away from the part of $S$  consisting of exactly three normally hyperbolic sheets. The superslow manifold in this scenario also consists of two attracting branches ${{\mathcal{Z}}^{\pm}}^a$, two repelling branches ${{\mathcal{Z}}^{\pm}}^r$, and a saddle branch ${{\mathcal{Z}}^{0}}$.  

Figure \ref{bursting} includes the bifurcation diagrams of the the fast subsystem (\ref{layer1}) for varying $\alpha$ and a fixed $\beta_2$, superimposed with trajectories of system (\ref{inter}) corresponding to the time series shown in figure \ref{timeseries_varying_alpha_1}, projected on the $(z, x)$-phase space. The bifurcation diagram is similar to  figure \ref{delayed_hopf}, except that the unstable branch of periodic orbits born at the subcritical Hopf bifurcation (SH) of the lower branch ${\mathcal{Z}}^-$  gains stability at a saddle-node bifurcation ($SN_p$), remains stable for decreasing $z$ until it loses stability at another $SN_p$, and thereafter terminates in a homoclinic bifurcation (BHC) at a saddle point of (\ref{layer1}).  The equilibrium on ${\mathcal{Z}}^-$ which is now an unstable focus/node, while the equilibrium  on ${\mathcal{Z}}^+$, a stable focus, are both enclosed by the homoclinic loop.  This loop has been referred to as  the ``big homoclinic loop" in Section 3. The unstable branch of periodic orbits born at SH on ${\mathcal{Z}}^+$ also terminates in HC with a nearby saddle point.  

  \begin{figure}[h!]     
  \centering 
\subfloat[]{\includegraphics[width=7.67cm]{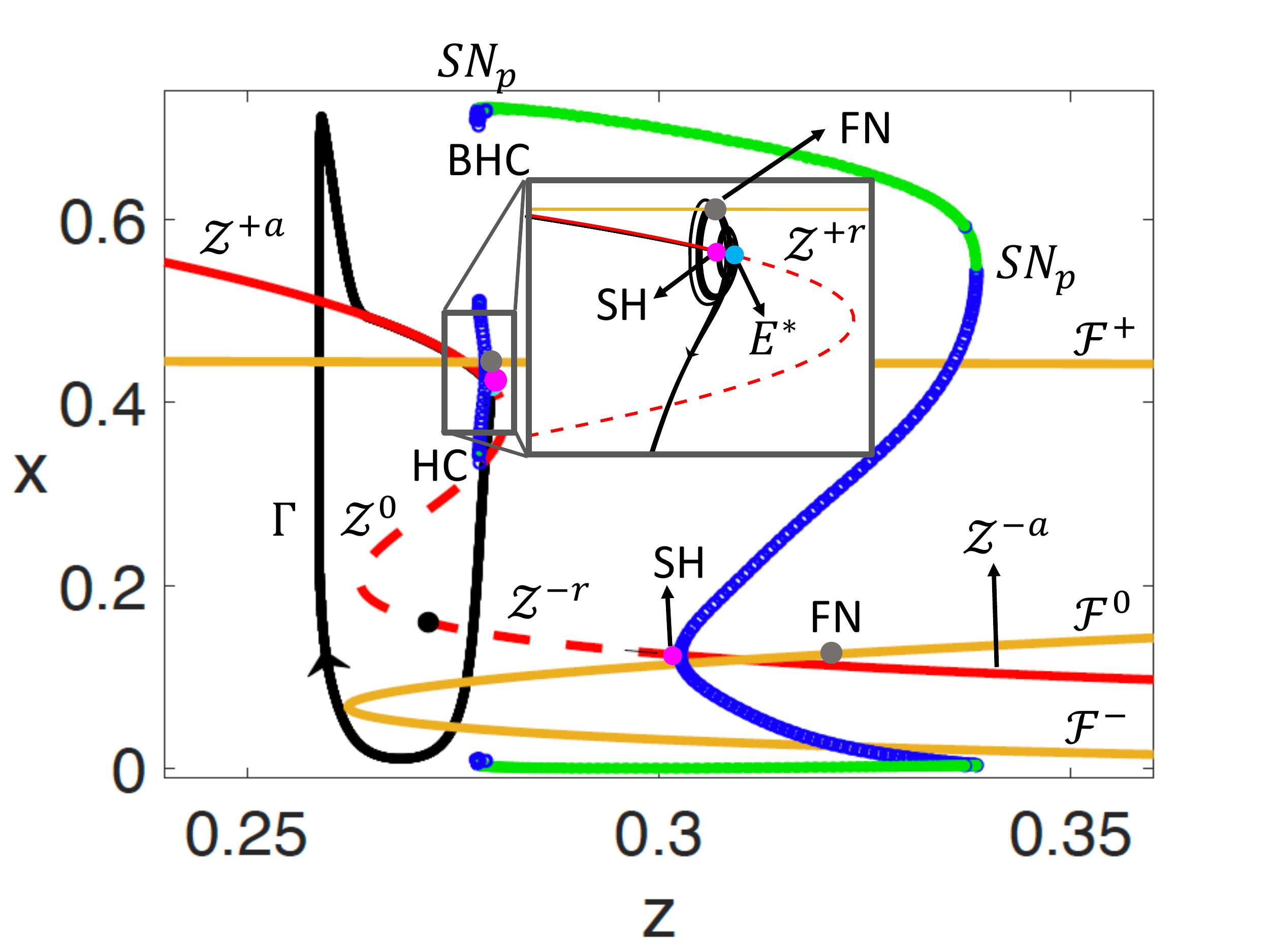}}
\quad
\subfloat[]{\includegraphics[width=7.67cm]{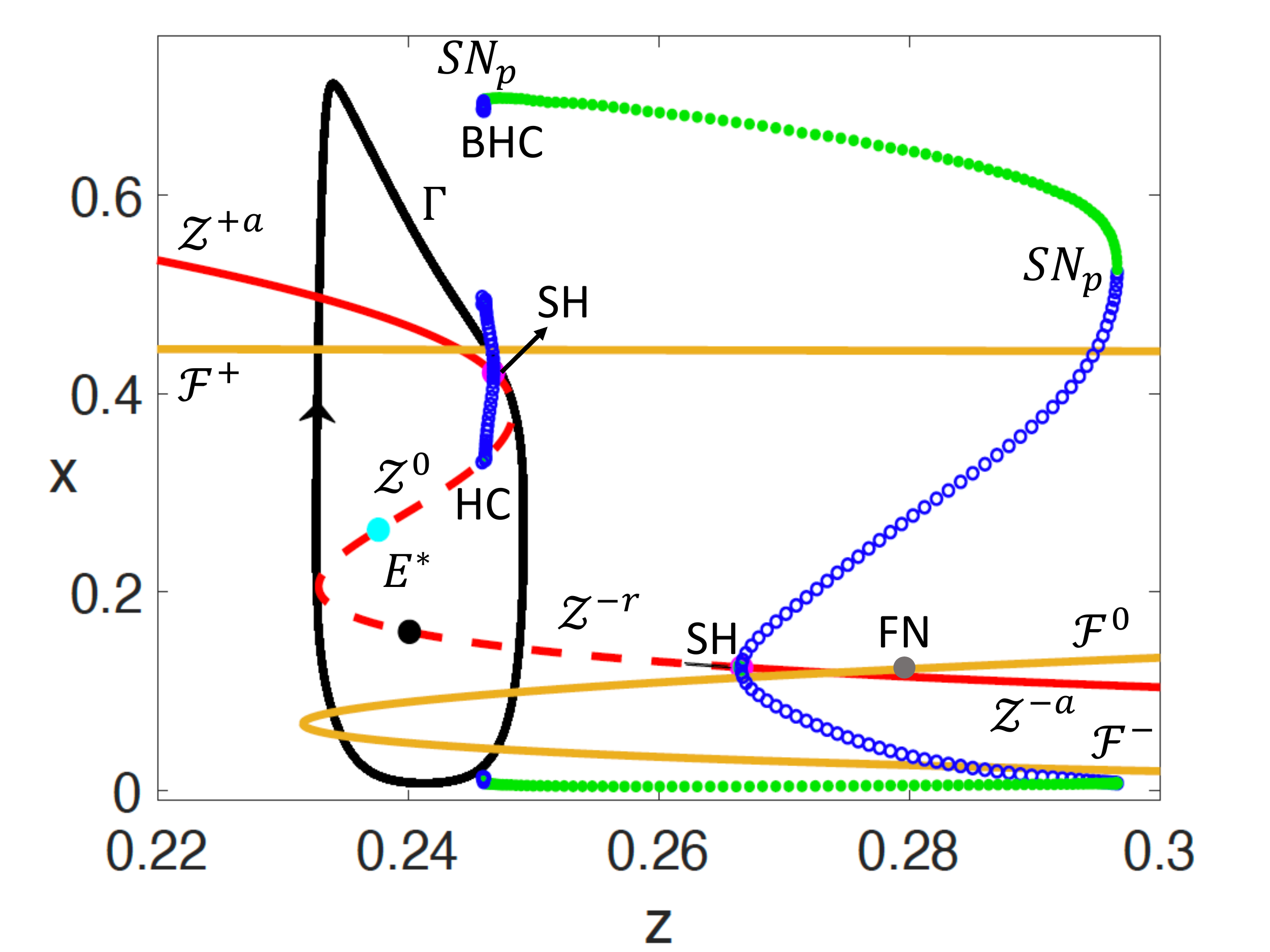}}
\quad
\subfloat[]{\includegraphics[width=7.7cm]{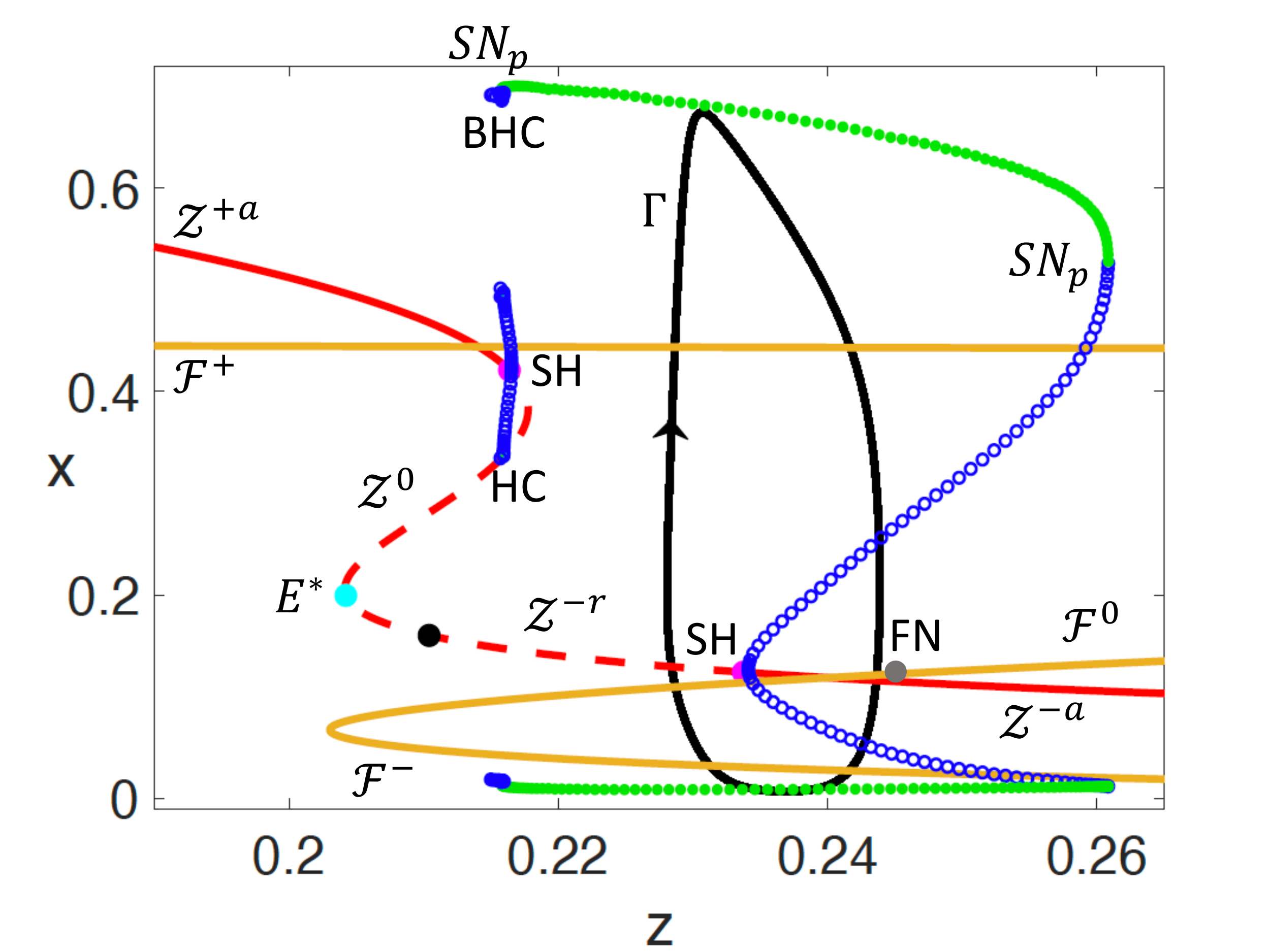}}
\quad
\subfloat[]{\includegraphics[width=7.7cm]{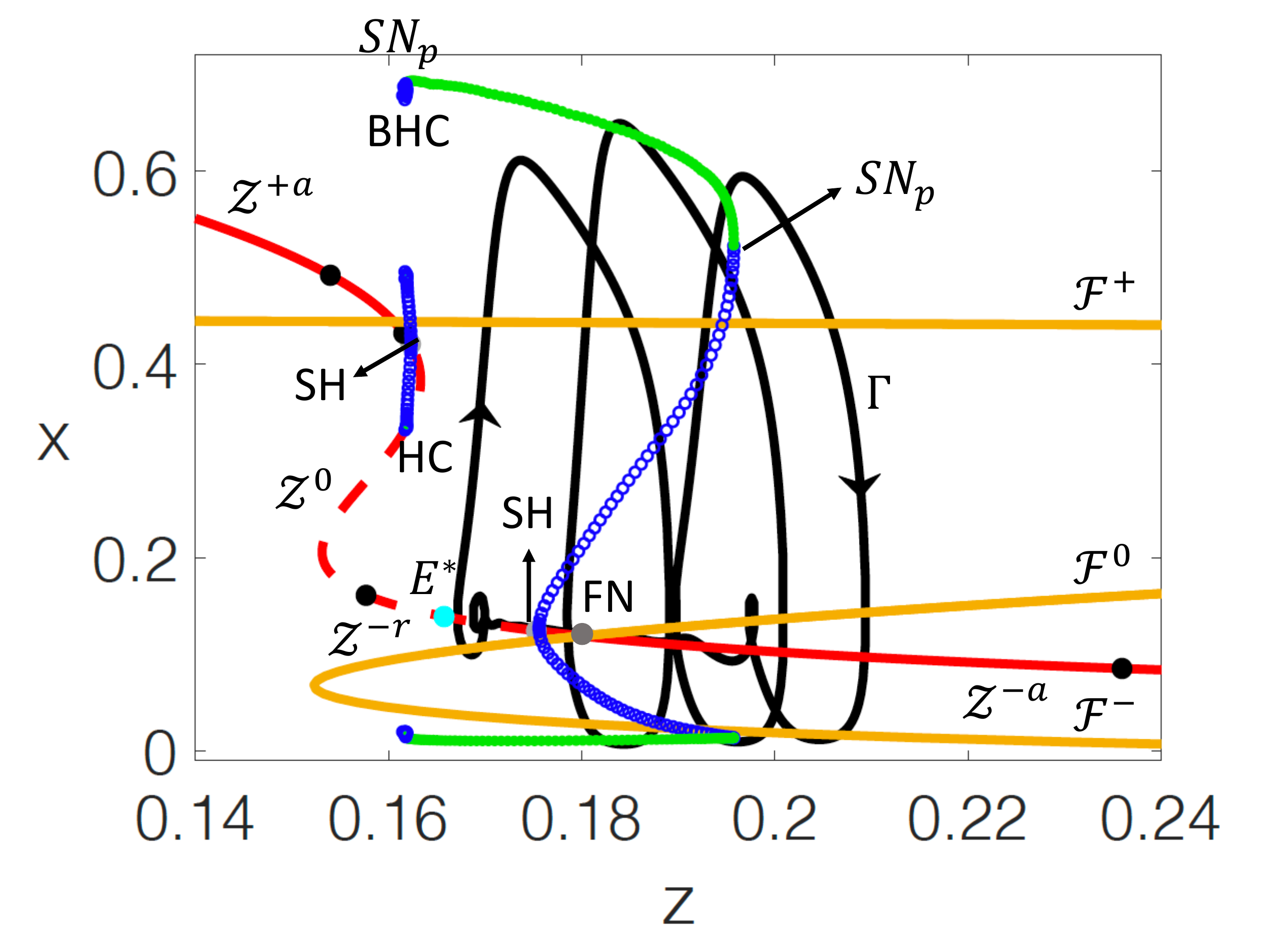}}
\caption{(A) Overlay of the projection of the trajectory $\Gamma$ of (\ref{inter}) onto the $(z, x)$ - plane with the bifurcation diagram of the fast subsystem (\ref{layer1}) for  $\beta_2 =0.0245$ and varying $\alpha$. Note the transition from spiking to   bursting patterns in the system. (A) $\alpha=0.4645$  (B) $\alpha=0.526$ (C) $\alpha=0.6$ (D) $\alpha =0.8$. $SN_p$ - saddle-node of periodics, BHC - big homoclinic loop, remaining labels and curves are as in figure \ref{delayed_hopf}. Filled green (open blue) circles represent maximum and minimum values of $x$ in stable (unstable) limit cycles.}
\label{bursting}%
\end{figure}

 The trajectory in figure \ref{bursting}(A) follows a close neighborhood of the upper branch $\mathcal{Z}^{+a}$ and exhibits small amplitude oscillations during its slow passage through the Hopf point SH $\in \mathcal{Z}_{DH}^{\varepsilon}$. As in  figure \ref{delayed_hopf}(A), in this case the orbit passes very close to a folded node singularity while making its way to the Hopf point, where its flow is governed by  (\ref{reduced2}). The local vector field of the nearby saddle-focus equilibrium $E^*$  initiates the small amplitude oscillations as shown in the inset of figure \ref{bursting}(A). The orbit then weakly follows the fast subsystem bifurcation diagram, and jumps to ${\Pi}^a_{\varepsilon, \delta}$ where it follows the intermediate flow and experiences a delayed loss of stability. Following a fast fiber, it jumps back to $S^{+a}_{\varepsilon, \delta}$ where the intermediate flow brings it to ${\mathcal Z}^+_{\varepsilon, \delta}$ and the cycle repeats.
  The phase portrait of the orbit along with the critical and superslow manifolds are shown in figure \ref{singular_funnel_proj_xz_1}. Similar to the dynamics in figure \ref{singular_funnel_proj_xz}, the orbit enters into the singular funnel on the attracting sheet $S^{+a}$ bounded by the singular strong canard $\gamma^+_{0, \delta}$ and the upper branch $\mathcal{F}^+$ of the fold curve. The SAOs in this case are organized by the unstable manifold of $E^*$.  However, note the difference in dynamics near the plane $\Pi$ in figure  \ref{singular_funnel_proj_xz} 
  and figure \ref{singular_funnel_proj_xz_1}. The difference arises due to the geometric structure of $S$ determined by $\mathcal{F}$, recalling that $\mathcal{F}$ belongs to case 3 in figure  \ref{singular_funnel_proj_xz} and to case 2 (ii) in figure \ref{singular_funnel_proj_xz_1}.

       \begin{figure}[h!]     
  \centering 
\subfloat[]{\includegraphics[width=7.7cm]{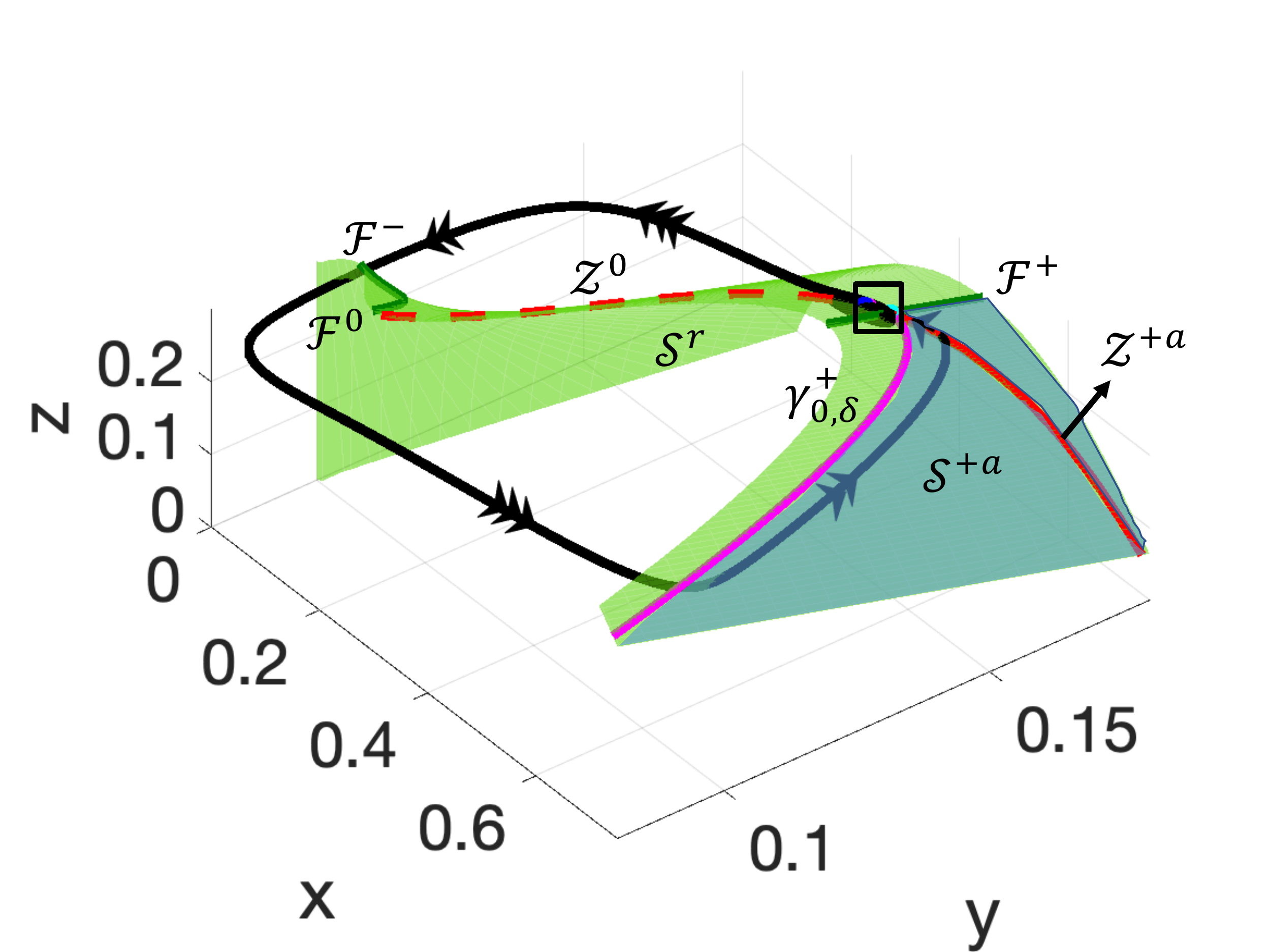}}
\quad
\subfloat[]{\includegraphics[width=7.7cm]{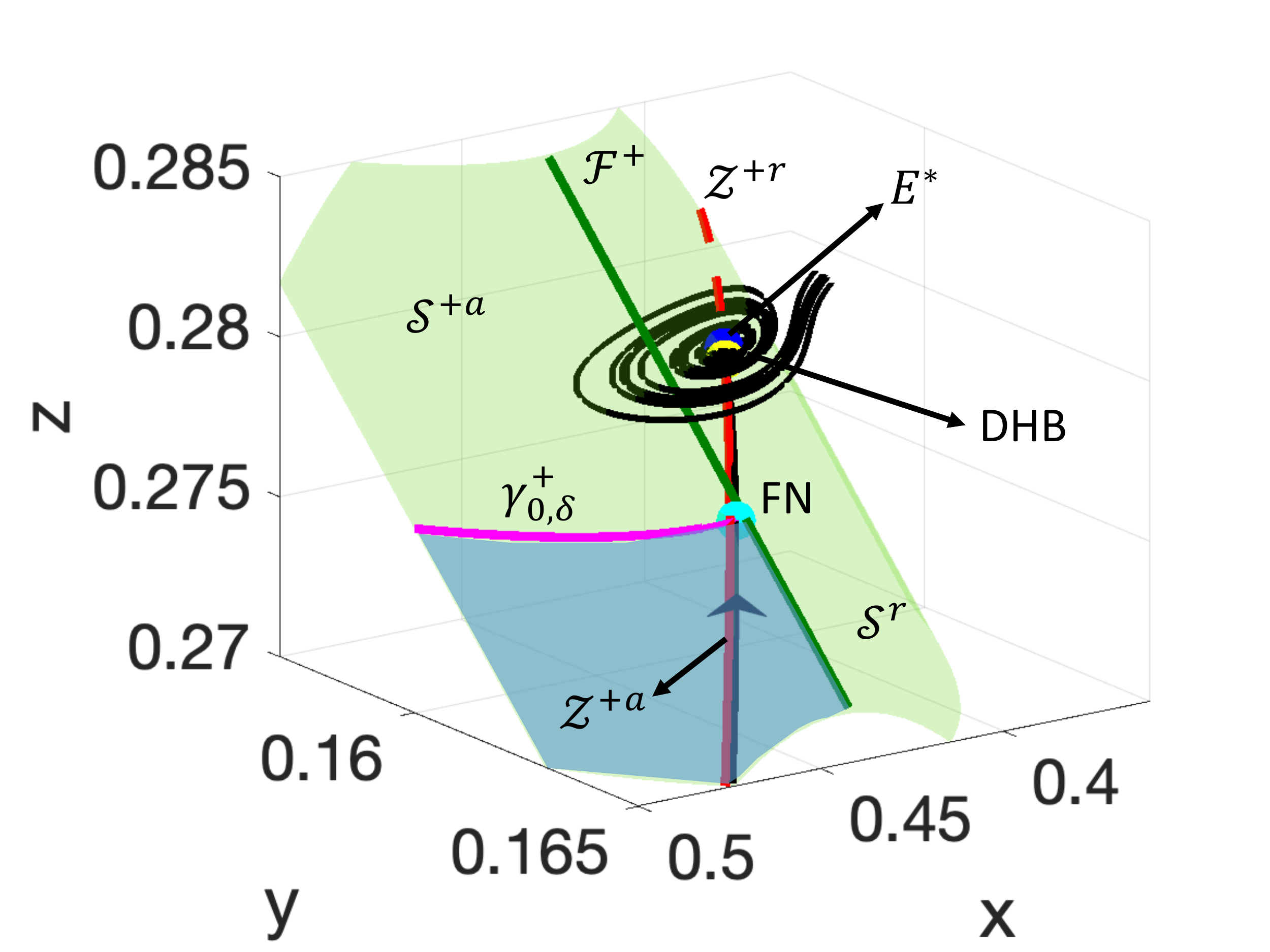}}
\caption{(A) Mixed-mode oscillations exhibited by system (\ref{inter}) for parameter values $\beta_2=0.0245$ and $\alpha=0.4645$. (B) Zoomed view of the dynamics near the upper fold $\mathcal{F}^+$ . The trajectory enters into the singular funnel on $\mathcal{S}^{+a}$ shown by the shaded region and filters through the the folded-node singularity (FN) while making its way to the delayed-Hopf bifurcation point (DHB). The local vector field of the equilibrium $E^*$ further influences its dynamics before it jumps to $S^{-a}$ .} 
\label{singular_funnel_proj_xz_1}%
\end{figure}

  The trajectories in figure \ref{bursting}(B)-(C) exhibit spiking behavior/relaxation oscillation type dynamics. The relaxation oscillation cycles in these figures involve only two timescales as they alternate between the fast timescale and the intermediate timescale on $S \setminus \mathcal{Z}$ and $\Pi$.  The dynamics of these orbits can be described as follows. Assuming that the orbit starts at a point on $S^{+a}_{\varepsilon, \delta}$, it follows the intermediate flow on $S^{+a}_{\varepsilon, \delta}$ until it reaches a vicinity of $\mathcal{F}^+$ where it gets connected by a fast fiber, and lands on  ${\Pi}^a_{\varepsilon, \delta}$.  As in  figure \ref{bursting}(A), the orbit undergoes the phenomenon of Pontryagin's delay of  loss of stability on the slow manifold ${\Pi}_{\varepsilon, \delta}$.   A fast fiber then concatenates with it and brings it to $S^{-a}_{\varepsilon, \delta}$, where it again follows the intermediate flow until it reaches a neighborhood of $\mathcal{F}^0$ and jumps to the opposite attracting branch $S^{+a}_{\varepsilon, \delta}$.  The cycle starts anew giving rise to relaxation-type dynamics.  We observe that in figure \ref{bursting}(b)-(c), the solutions do not get attracted to $\mathcal{Z}^{+a}$ in contrast to the solution shown in figure  \ref{bursting}(a). The difference in the behavior of the solutions can be attributed to the location of the folded singularities in system (\ref{nondim2}).  Note that the parameter values considered in  figure \ref{bursting}(a) and figures \ref{bursting}(b)-(c) lie in regions B and C of figure \ref{two_par_desing} respectively.  In region B, the reduced problem (\ref{reduced})  restricted to $S$ has two folded nodes, whereas it has only one folded node in region C. In  figure \ref{bursting}(a), the folded nodes lie on ${\mathcal{F}}^{\pm}$, whereas in figures \ref{bursting}(b)-(c), the folded node lies on ${\mathcal{F}}^{-}$. The  absence of folded node singularities on ${\mathcal{F}}^{+}$ leads to the absence of the singular funnel and the primary weak canard (which when exists, lies close to $\mathcal{Z}^{+a}$) on $S^{+a}$ for the parameter values considered in figure \ref{bursting}(b)-(c).  As a result the intermediate flow on $S^{+a}$ in such a scenario brings  an orbit directly  to $\mathcal{F}^+$ without necessarily approaching $\mathcal{Z}^{+a}$. Consequently,  the solutions shown in figure \ref{bursting}(b)-(c) do not  approach $\mathcal{Z}^{+a}$ and thus do not evolve along $\mathcal{Z}^{+a}$.

  In figure \ref{bursting}(D), the trajectory  exhibits a subHopf/fold cycle bursting pattern \cite{izh}, where the orbit closely follows  $\mathcal{Z}^{-a}$ while approaching the subcritical Hopf point SH $\in \mathcal{Z}_{DH}^{\varepsilon}$ located on that branch and experiences a slow passage through that point. It then spirals out to the attracting branch of the periodic orbits of the fast subsystem (\ref{layer1}) shadowing it, while slowly drifting towards larger values of $z$. This phase terminates when the trajectory approaches the $SN_p$ point and falls off the attracting branch of the periodic orbits. It then spirals in towards  $\mathcal{Z}^{-a}$ and the cycle repeats. Note the presence of bistability between $\mathcal{Z}^{-a}$ and the attracting branch of the periodic orbits of the fast subsystem for sustaining  such bursting patterns. The time series and the phase portrait of the bursting orbit are shown in figures \ref{timeseries_example_2} and \ref{crit_mfld}(B) respectively. According to the classification in \cite{izh}, such dynamics is referred to as a subcritical elliptic bursting as it is characterized by a subcritical Hopf bifurcation (SH) of the fast subsystem which triggers the onset of the spikes and a saddle-node ($SN_p$) of limit cycles which terminates the spikes. These dynamics can be mapped to the region labeled as \fbox{1} in figure \ref{two_par_fast}. The solution in figure \ref{bursting}(d) features three timescales as it alternates between the fast timescale, the intermediate timescale on $S^{\pm a}\setminus \mathcal{Z}^{\pm a}$, and the slow timescale on $\mathcal{Z}^{-a}$.  Assuming that the orbit starts on $S^{+a}_{\varepsilon, \delta}$, the intermediate flow brings it to a vicinity of $\mathcal{F}^+$ where it concatenates with a fast fiber and reaches $S^{-a}_{\varepsilon, \delta}$. The orbit now follows the intermediate flow on $S^{-a}_{\varepsilon, \delta}$ and may either reach a vicinity of $\mathcal{F}^0$  or get attracted to $\mathcal{Z}^{- a}_{\varepsilon, \delta}$. In the former case, it jumps back to $S^{+a}_{\varepsilon, \delta}$ and the cycle repeats until the orbit  gets attracted to  $\mathcal{Z}^{- a}_{\varepsilon, \delta}$.  In the latter case, the orbit slowly evolves along $\mathcal{Z}^{- a}_{\varepsilon, \delta}$ and experiences a slow passage through a Hopf point before it spirals out to the attracting branch of the periodic orbits of the fast subsystem as described above. In a recent work \cite{JKK}, the blow-up method was applied to analyze the global oscillatory transition near a regular folded limit cycle manifold in a class of three time-scale systems with two small parameters. It will be interesting to see if the method in \cite{JKK} can be extended to analyze the different kinds of bursting phenomena in this model.


\section{Discussion and future outlook}

Explanation of large population variations have ranged from availability of resources to predators, disease, seasonal reproductive cycles, precipitation patterns, temperature changes, human activities and so forth. In some species,  evidence of seasonal changes in their population abundance can be correlated to precipitation patterns and temperature changes. However, in ecosystems where seasonality is not very pronounced, trophic interactions, more specifically top-down regulations can play a central role in organizing population cycles, as evidenced by the studies in \cite{grotan,molleman}. Furthermore, the abundance of the type of predator characterized by their feeding preferences can play a major role in influencing the dynamics of prey  \cite{Hanski et al}. For instance,  it was suggested in \cite{Hanski et al} that  generalist predators (foxes, cats, common buzzards) seem to stabilize the populations of microtine rodents in the southern Fennoscandia, whereas specialist mammalian predators (small mustelids) seem to significantly contribute to their regular multiannual cycle  in northern Fennoscandia.  Hence, understanding the roles of generalist and specialist predators in regulating and shaping population dynamics is crucial in any ecosystem and has been an intriguing subject of interest in ecology. To that end, in this paper, we studied the dynamics between three interacting species, namely two classes of predators (specialist and generalist) competing for a common prey with the assumption that the predators do not prey upon each other.

Taking into account that each species operates on a different timescale, we introduced separation of timescales in the model, and obtained a slow-fast system featuring three timescales. Grouping the timescales by using $\varepsilon$ as the singular parameter and $\delta$ fixed, we partitioned the system into a one-dimensional fast subsystem described by the $x$ dynamics and a two-dimensional slow subsystem described by the $(y,z)$ dynamics. Similarly, treating $\delta$ as the singular parameter and $\varepsilon$  fixed, we obtained a family of two-dimensional $(x,y)$ fast subsystem parameterized by $z$.  We studied the role of critical and superslow manifolds in shaping the dynamics arising in this model. In contrast to other commonly studied slow-fast models that are motivated by applications in biological sciences, chemistry and ecology (see \cite{BKR, DGKKOW, K11, KC, LRV, MRfoodchain, VBW} and the references therein), we note that in this model, the component  of the critical manifold formed by the nontrivial nullcline of the fast subsystem, namely the surface $S$, is not uniformly ``S-shaped" and may contain up to two cusp points. The number of normally hyperbolic sheets of $S$ could vary between two to four, which gave a rich geometric structure to $S$. Moreover, the self-intersecting feature of the critical manifold gave rise to a bifurcation delay which was manifested in orbits during their passage past the invariant plane $\Pi$. Such dynamics have been studied in relatively few higher-dimensional models; some examples include \cite{kpnew, Sadhudcds}.

We applied slow-fast analysis techniques from GSPT and used two complementary geometric methods to examine the dynamics. The fast-slow decomposition methods for the two-timescale systems naturally extended to the three-timescale setting \cite{CT, kpk, LRV}. We noted that in both one-fast/two-slow and two-slow/one-fast analyses, the underlying geometry influences the different oscillatory patterns through a combination of local and global mechanisms. The one-fast/two-slow decomposition led to key structures such as the critical manifold and its singularities and gave insight to the mechanism of canard dynamics which could play a role in organizing the small amplitude oscillations in MMO patterns in the system. The two-fast/one-slow decomposition led to another set of key structures such as the superslow manifold and its folds, Hopf bifurcation points in the layer problem, and gave insight to the mechanism of dynamic Hopf bifurcation which also contributes to organizing the small amplitude oscillations in MMOs in this system. In addition, the two-fast/one-slow  analysis also provided a guide to locating transitions between spiking and slient phase in bursting patterns. From the viewpoint of this analysis, we noted that  the small oscillations in the bursting dynamics are generated from a slow passage through a dynamic Hopf bifurcation, and the large amplitude oscillations are hysteresis loops that alternately jump between two subcritical Hopf bifurcations (subHopf/subHopf bursts) or  at a subcritical Hopf and  a cyclic fold (subHopf/cyclic fold).  A two-parameter bifurcation analysis of the fast subsystem parametrized by the slow variable $z$ revealed  several interesting codimension-one  bifurcations such as subcritical Hopf, homoclinic, saddle-node of periodic orbits and codimension-two  bifurcations such as generalized Hopf, saddle-node separatrix loop and cusp. Transitions between the spiking and quiescent dynamics in the bursting oscillations were associated with these global bifurcations.

 Treating the efficiency  of the generalist predator, $\beta_2$, as the primary varying parameter, and $\alpha$, the proportion of diet of the generalist predator that consists of $x$, as the secondary parameter,  we note that the system progressed through different oscillatory regimes such as canard or delayed-Hopf induced MMOs, relaxation oscillation cycles featuring three timescales,  and subcritical elliptic bursting. Such oscillatory patterns have been studied in  neurological models, chemical kinetics and other prototypical three-timescale models (see \cite{DKP, kpk, LRV} and the references therein) but are novel in an ecological setting. In addition, oscillatory dynamics featuring a slow passage near the plane $\Pi$ before getting attracted to the adjacent attracting sheet of $S$ as   shown in figures \ref{crit_mfld}(a) and  \ref{singular_funnel_proj_xz} are novel in the three-timescale setting. These dynamics can be associated  with different types of cyclic patterns of population densities seen in various ecosystems \cite{nelson, singleton,stenseth} and perhaps can be  attributed to the role of a generalist predator in regulating the cycles.  In particular,  we found that a highly efficient generalist predator (i.e. for lower values of $\beta_2$) can keep the prey population at a very low density for a prolonged time, until its density slowly decreases below a certain threshold allowing the prey density to rise sharply, giving rise to cycles in the form of MMOs (figures \ref{crit_mfld}(a) and  \ref{singular_funnel_proj_xz}).  As the efficiency of the generalist predator decreases (i.e. as $\beta_2$ increases), we obtained bursts of high-frequency oscillations in the prey density as it exhibits a series of multiple outbreaks, giving rise to bursting oscillations (figure \ref{crit_mfld}(b)).  At a further reduced efficiency of the generalist predator, regular boom and bust cycles or relaxation oscillation dynamics are observed. We remark that in a three-timescale food chain model studied in \cite{MRfoodchain}, it was interpreted that high efficiency of the top predator implies cycles in a food chain. 
 In this paper, we obtain a general result in a similar spirit, namely the predation efficiency of the generalist predator influences the type of oscillatory pattern in a food-web, thus underscoring the role of a generalist predator in regulating population dynamics.

 The system also exhibited other types of interesting dynamics such as amplitude-modulated spiking, classical  torus canards and mixed-type torus canards. Torus canards and mixed-type torus canards have been well studied in two-timescale neuronal and chemical models (see \cite{BAD, burketal, DBKK, straube, Vo} and the references therein), but are novel in three-timescale settings. In this model, these dynamics are observed during transition to and from subHopf/fold cycle bursting (see figures \ref{two_par_bif_full},  \ref{one_par_bif}, and \ref{torus_canards1}). The mixed-type torus canards separate the regimes of Hopf cycles and MMOs featuring single spike accompanied with a long quiescent phase. These solutions occur in a very narrow parameter range and are preceded by amplitude-modulated spiking dynamics in the parameter space c.f. \cite{BAD, DBKK}. The long quiescent phase in the MMOs is organized by a homoclinic bifurcation in the fast subsystem.  The classical torus canards, on the other hand, mediate the transition between subcritical elliptic bursting  to relaxation oscillations  c.f. \cite{burketal}.  A detailed analysis of these dynamics is beyond the scope of this paper and is left for the  future.  We also remark that the transitions from MMOs exhibiting SAOs along $\mathcal{F}^0$ to MMOs exhibiting SAOs along $\mathcal{F}^+$ can be explained by the singular geometry and the reduced flows as in  \cite{KPK}. However, in contrast to \cite{KPK}, in this model, the folded singularity does not necessarily lie on $\mathcal{Z}$, the singular geometry is non-symmetric, and therefore, one may need additional work to classify the MMO dynamics using the singular flows. We leave this subject for  future study as well.

The presence of a folded node singularity in vicinity of a delayed-Hopf point and an unstable equilibrium of saddle-focus makes the dynamics near the folds all the more interesting.  In this model, the local vector field of the equilibrium plays a crucial role in organizing the local oscillations near the fold, while the delayed-Hopf point and the folded node singularity are instrumental in guiding the trajectory to the equilibrium. It will be interesting to investigate the precise dynamical mechanism inducing the jump from $S$ towards the critical manfold $\Pi$, and use the analysis to detect early warning signs of  an outbreak as has been carried out in  \cite{Sadhunew} for a two-timescale predator-prey model.  We leave this subject for  future study.

We finally remark that though the model considered in this paper is generic, yet it produces a host of interesting  oscillatory patterns, some of which qualitatively represent population patterns observed in  small mammals or insects. For instance, the time profile of the MMO orbit in figure \ref{timeseries_example_1}  showing patterns of long epochs of small amplitude oscillations near the middle branch of the fold curve, periodically  interspersed with large-amplitude fluctuations, qualitatively resembles population densities of microtine rodent populations \cite{singleton,stenseth}, though the large-amplitude variations in such populations are more sporadic in nature. Another interesting pattern is the time profile of a subcritical elliptic bursting characterized by a sequence of recurrent high-amplitude fluctuations separated by a transient low density state shown  in figure \ref{timeseries_example_2}(A). Such a pattern qualitatively resembles population densities of multivoltine insects such as smaller tea-tortrix, a pest on tea leaves in Japan, which may sporadically exhibit multiple outbreaks annually as has been studied in \cite{nelson}. It will be interesting to study the effect of stochasticity on this system in parameter regimes associated with MMOs and bursting oscillations. We also leave this subject for future study.


\section{Acknowledgements}  The author would like to thank the referees for their careful reading of the original manuscript and many valuable comments and suggestions that greatly improved the presentation of this paper.

\section*{Appendix:  Classification of the fold curve $\mathcal{F}$}

Recalling that $\mathcal{F}=  \left\{(x, y, z) \in S:   y=\mu(x),\  z= \nu(x), \ x \in [0, 1]\setminus \{x_d\} \right\}$, where $\mu(x)$ and $\nu(x)$ are defined by (\ref{foldzcomp}), we note that $\nu(0)>0$, $\nu'(0)<0$, and $\nu(x)$ has a unique root at $x=(1-\beta_1)/2$ and an infinite discontinuity at $x=x_d$. Similarly $\mu(0)>0$, $\mu'(0)>0$ and $\mu(x)$ also has an infinite discontinuity at $x=x_d$. Depending on the value of $\beta_2$, $\mu(x)$ has either no roots or two roots (repeated or distinct) in $(0, 1)\setminus \{x_d\}$ if $x_d<1$ and in $(0, 1)$ if $x_d>1$. It is clear from (\ref{foldzcomp}) that $\nu(x)\mu(x)<0$ in a neighborhood of $x_d$. Hence the fold curve $\mathcal{F} \not\subset {\mathbb{R}^3}^+$ if $x$  lies in  a neighborhood of $x_d$. To this end, we define  the points $x_1, x_2$ by
 \bess x_1= \textnormal {min} \left\{\frac{1-\beta_1}{2}, x_d \right\} ,\  x_2= \textnormal {max} \left\{\frac{1-\beta_1}{2}, x_d \right\},
\eess
and consider the  following cases determined by the number of extreme values of  $\nu(x)$ in the interval $[0, x_1]$:

 Case 1: $x_1=\frac{1-\beta_1}{2}$ and $\nu(x)$ has no relative extrema in $[0, x_1]$.\\
 In this case, $\nu(x)>0$ on $[0, x_1)\cup (x_d, 1]$ if $x_d<1$, and $\nu(x)>0$ on $[0, x_1)$  if $x_d\geq 1$. In either case, $\mu (x)>0$ only when $ x \in [0, x_d)$. Hence, $\mathcal{F} \subset {\mathbb{R}^3}^+$ if $x \in [0, x_1]$. The curve $\mathcal{F}$ is monotonic (does not have any folds) in ${\mathbb{R}^3}^+$.  The surface $S$ will be uniformly divided into an attracting sheet and a repelling sheet that meet at $\mathcal{F}$.
 
  Case 2: $x_1=\frac{1-\beta_1}{2}$ and $\nu(x)$ has two relative extreme points in $[0, x_1]$.  Here we can have two sub-cases:\\
  (i) $\mu(x)$ has no zeros in $[0, 1]\setminus \{x_d\}$ if $x_d<1$  or $[0, 1]$ if $x_d>1$. \\
  In this case, $\mu(x)>0$ if $x \in [0, x_d)$ and $\nu(x)>0$ if $x\in [0, x_1)$. 
  Similar to Case 1,  $\mathcal{F}\subset {\mathbb{R}^3}^+ $ if $x\in [0, x_1)$. However, in this situation,  $\mathcal{F}$ is cubic-shaped and divides the surface into four different regions (see figure \ref{fold_curves}(A)). Denoting the locations of the relative minimum and maximum of $\nu(x)$ by $x_m$ and $x_M$ respectively, where $x_m<x_M$, we can write $ \mathcal{F}  = \mathcal{F}^- \cup \mathcal{F}^0  \cup \mathcal{F}^+$ with
\bess  \mathcal{F}^- & = &\left\{(\eta, F(\eta, z(\eta)), z(\eta)) \in {\mathbb{R}^3}^+ :  z(\eta)=\nu(\eta), \  0\leq \eta < x_m\right\},\\
 \mathcal{F}^0 & = &\left\{(\eta, F(\eta, z(\eta)), z(\eta))  \in {\mathbb{R}^3}^+ :  z(\eta)=\nu(\eta), \ x_m\leq \eta \leq  x_M\right\},\  \textnormal {and} \\
\mathcal{F}^+  &= &\left\{(\eta, F(\eta, z(\eta)), z(\eta))  \in {\mathbb{R}^3}^+ :  z(t)= \nu(t),\  x_M< \eta \leq x_1 \right\}.
\eess
 We then have that for $0\leq y\leq \mu(x_M)$, $S$ is uniformly attracting. For $\mu(x_M)<y \leq \beta_1$, $S$ has two attracting sheets, ${S_{\pm}^a}$, and one repelling sheet, ${S^r}$, joined along the two branches $\mathcal{F}^+$ and $\mathcal{F}^0$  of the fold curve.  For $\beta_1< y\leq \mu(x_m)$, $S$ has two attracting sheets ${S_{\pm}^a}$ and two repelling sheet ${S_{\pm}^r}$ separated by the three branches of $\mathcal{F}$. For $y> \mu(x_m)$, $S$ has an attracting sheet, $S^a$, and a repelling sheet $S^r$.  Figure \ref{fold_curves}(A) shows the variation in the number of attracting and repelling branches of $S$ with $y$.

  (ii) $\mu(x)$ has two repeated or distinct roots in $[0, x_1)$.\\
   In this case, $\mathcal{F}$ is defined piecewise, dividing the surface into three different regions. We may write $ \mathcal{F}  = \mathcal{F}^- \cup \mathcal{F}^0  \cup \mathcal{F}^+$, where 
\bess  \mathcal{F}^- & = &\left\{(\eta, F(\eta, z(\eta)), z(\eta))  \in {\mathbb{R}^3}^+ :  z(\eta)=\nu(\eta), \ 0\leq \eta < x_m\right\},\\
 \mathcal{F}^0 & = &\left\{(\eta, F(\eta, z(\eta)), z(\eta))  \in {\mathbb{R}^3}^+ :  z(\eta)=\nu(\eta), \ x_m\leq \eta \leq  {\mu}_1\right\},\  \textnormal {and} \\
\mathcal{F}^+  &= &\left\{(\eta, F(\eta, z(\eta)), z(\eta))  \in {\mathbb{R}^3}^+ : z(\eta)=\nu(\eta), \ {\mu}_2\leq \eta \leq x_1 \right\},
\eess
where ${\mu}_1\leq {\mu}_2$ are roots of $\mu(x)=0$  that lie in the interval $(0, x_1)$. The surface $S$ has two attracting  branches and one repelling branch, ${S_{\pm}^a}$ and  ${S^r}$ respectively for $0\leq y\leq \beta_1$, two attracting and two repelling branches, ${S_{\pm}^a}$ and ${S_{\pm}^r}$ respectively for $\beta_1<y<\mu(x_m)$, and an attracting branch and a repelling branch, $S^a$ and $S^r$ respectively for $y>\mu(x_m)$ (see figures \ref{crit_mfld}(B) and \ref{fold_curves}(B)).

Case 3.  $x_1 =x_d$ and $\nu(x)$ has exactly one relative extreme point in $(0, x_d)$.\\
Since $\nu(0)>0$, $\nu'(0)<0$ and  $\nu(x) \to \infty$ as $x \to { x_d}^-$, $\nu(x)$ attains its local minimum  at $x_m \in (0, x_d)$. In this case, $\nu(x)>0$ on $[0, x_d)$ and $(x_2, 1]$, and $\mu(x)>0$ on $[0, {\mu}_1) \cup (x_d, {\mu}_2)$, where  ${\mu}_1$, ${\mu}_2 \in (0, 1)$ are positive roots of $\mu$ such that ${\mu}_1<x_d$ and ${\mu}_2>x_2$ (Note from (\ref{foldzcomp})  that $\mu(x_2)>0$, hence ${\mu}_2>x_2$.) It then follows that the fold curve $\mathcal{F}\subset {\mathbb{R}^3}^+$ if $x \in [0, {\mu}_1) \cup (x_2, {\mu}_2)$, and is therefore piecewise continuous.  We may write  $ \mathcal{F}  = \mathcal{F}^- \cup \mathcal{F}^0  \cup \mathcal{F}^+$, where 
\bess  \mathcal{F}^- & = &\left\{(\eta, F(\eta, z(\eta)), z(\eta))  \in {\mathbb{R}^3}^+ :  z(\eta)=\nu(\eta), \ 0\leq \eta <  x_m\right\} \\
 \mathcal{F}^0 & = &\left\{(\eta, F(\eta, z(\eta)), z(\eta))  \in {\mathbb{R}^3}^+ :  z(\eta)=\nu(\eta),\  x_m\leq \eta \leq  {\mu}_1\right\},\  \textnormal {and} \\
\mathcal{F}^+  &= &\left\{(\eta, F(\eta, z(\eta)), z(\eta))  \in {\mathbb{R}^3}^+ :  z(\eta)=\nu(\eta),\  x_2\leq \eta \leq {\mu}_2 \right\}.
\eess
The number of attracting and repelling sheets that $S$ possesses is similar to Case 2 (ii) (see figure \ref{crit_mfld}(A)).


\medskip
\medskip


\begin{thebibliography}{100}

\bibitem{AS} S. Ai and S. Sadhu, {\em The entry-exit theorem and relaxation oscillations in slow-fast planar systems}, Journal of Diff. Eq., 268, (2020) 7220 - 7249.

\bibitem{AY} S. Ai and Y. Yingfei,  {\em Relaxation Oscillations in Predator–Prey Systems}, Journal of Dynamics and Differential Equations, (2021)1-28. 

\bibitem{BE} S. M. Baer and T. Erneux, {\em Singular Hopf Bifurcation to Relaxation Oscillations}, SIAM J. Appl. Math., 46, (1986), 721 39. 

\bibitem{BAD} E. Baspinar, D. Avitabile, and M. Desroches, {\em Canonical models for torus canards in elliptic bursters}, Chaos, 31 (2021) 063129.

\bibitem{bazykin} A.D.  Bazykin,  {\em Nonlinear Dynamics of Interacting Populations}, Series A: Monographs and Treatise; World Scientific Series on Nonlinear Science: Singapore, 1998.

\bibitem{BLL}  N. Bolohan, V.  LeBlanc and F. Lutscher, {\em Seasonal dynamics of a generalist and a specialist predator on a single prey}, Math. Appl. Sc. and Engr.  2, (2021) 72 - 148.

\bibitem{BKK} H. Broer, T.J. kaper and M. Krupa {\em Geometric desingularization of a cusp singularity in slow-fast systems with applications to Zeeman's examples}, Journal of Dynamics and Differential Equations, Springer Verlag, 2013.

\bibitem{BKR} M. Br${\mathrm{\o}}$ns, T.J. Kaper, and H. G. Rotstein, {\em Introduction to Focus Issue: Mixed Mode Oscillations: Experiment, Computation, and Analysis}, Chaos 18, 015101 (2008).

\bibitem{BK} M. Br${\mathrm{\o}}$ns and R. Kaasen, {\em Canards and mixed-mode oscillations in a forest pest model}, Theoretical Population Biology 77 (2010) 238-242. 

\bibitem{BKW} B.M. Br${\mathrm{\o}}$ns,  M. Krupa, M. Wechselberger, {\em Mixed Mode Oscillations Due to the Generalized Canard Phenomenon}, Fields Institute Communications 49 (2006) 39-63.

\bibitem{burketal} J. Burke,  M. Desroches,  A.M. Barry, 
T. J. Kaper, and M. A Kramer, {\em A showcase of torus canards in neuronal bursters}, J. Math.  Neurosci. 2: 3  (2012)  2-30.

\bibitem{CT} P.T. Cardin and M. A.  Teixeira, {\em Fenichel Theory for Multiple Time Scale Singular Perturbation Problems}, SIAM J. Appld. Dyn. Syst. 16 (2017) 1425- 1452.

\bibitem{DKP} P. De Maesschalck, E. Kutafina, and N. Popovic, {\em Sector-delayed-Hopf-type  
mixed-mode oscillations in a prototypical three-time-scale model}, Appl. Math.
Comput. 273,  (2016) 337 - 352. 

\bibitem{BD1} B. Deng, {\em Food chain chaos due to junction-fold point}, Chaos 11 (2001) 514-525.

\bibitem{DH} B. Deng and G. Hines, {\em Food chain chaos due to Shilnikov's orbit}, Chaos 12 (2002) 533-538.

\bibitem{DH1}  B. Deng and G. Hines, {\em Food chain chaos due to transcritical point} Chaos 13 ( 2003 ) 578 - 585

\bibitem{BD2} B. Deng, {\em Food chain chaos with canard explosion},  Chaos 14 ( 2004) 1083 - 1092.

\bibitem{DBKK} M. Desroches, J. Burke, T.J. Kaper, and M.A. Kramer, {\em Canards of mixed type in a neural burster}, Phy. Review E 85 021920 (2012).

\bibitem{DGKKOW} M. Desroches, J. Guckenheimer, B. Krauskopf, C. Kuehn, H.M. Osinga, M. Wechselberger, {\em Mixed-Mode Oscillations with Multiple Time Scales}, SIAM Review 54 (2012) 211-288.

\bibitem{DZL} L. Duan, D. Zhai, and Q. Lu, {\em Bifurcation and bursting in Morris-Lecar model for class I and class II excitability}, Disc.\& Cont. Dynm. Syst (S) (2011) 391-399.

\bibitem{ELS} A. Erbach, F. Lutscher, G. Seo, {\em Bistability and limit cycles in generalist predator–prey dynamics}, Ecol. Complexity, 14 (2013) 48-55.

\bibitem{F} N. Fenichel, {\em Geometric singular perturbation theory for ordinary differential equations},  J. Diff. Eq., 31(1979) 53-98.

\bibitem{GTW} J. Guckenheimer, J.H. Tien and A.R. Willms, {\em Bifurcations in the Fast Dynamics of Neurons: Implications for Bursting}, in Bursting: The Genesis of Rhythm in the Nervous System, World
Scientific, (2005) 89 - 122.

\bibitem{grotan} V. Grotan, R. Lande, S. Engen, B.E Saether, P.J. DeVries, {\em Seasonal cycles of species diversity and similarity in a tropical butterfly community}, J. of Animal Ecol. 81 (2012) 714 - 723.

\bibitem{Hanski et al} I. Hanski, L. Hansson and H. Henttonen, {\em Specialist Predators, Generalist Predators, and the Microtine Rodent Cycle}, J. Animal Ecol. 60 (1991) 353 - 367.

\bibitem{Hassell et al} M. P. Hassell and R. M. May,  {\em Generalist and Specialist Natural Enemies in Insect Predator-Prey Interactions}, J. of Animal Ecol., 55, No. 3 (1986)  923 - 940.


 \bibitem{izh} E. M. Izhikevich, {\em Neural excitability, spiking and bursting}, Internat. J. Bifur. Chaos  Appl. Sci.  Engrg., 10 (2000), pp. 1171 - 1266.
 
 \bibitem{JKK} S. Jelbert, S.-V. Kuntz and C. Kuehn, {\em Geometric Blow-up for Folded Limit Cycle Manifolds in Three Time-Scale Systems}, 	arXiv:2208.01361.

\bibitem{kpnew} P. Kaklamanos and  N. Popovi\'c, {\em Complex oscillatory dynamics in a three-timescale El Ni\~no Southern Oscillation model}, arXiv:2207.03230.

\bibitem{kpk} P. Kaklamanos, N. Popovi\'c, and K. U. Kristiansen, {\em Bifurcations of mixed-mode oscillations in three-timescale systems: an extended prototypical example}, Chaos 32 013108 (2022).
 
\bibitem{KW}  Y. Kang and L. Wedekin, {\em  Dynamics of a intraguild predation model with generalist or specialist predator}. J. Math. Biol. 67, 1227 - 1259 (2013). 
 

\bibitem{K11} C. Kuehn, {\em Multiple Time Scale Dynamics} Springer (2015).

\bibitem{KPK} M. Krupa, N. Popovi\'c and N. Kopell, {\em Mixed-mode oscillations in three time-scale systems: a prototypical example}, SIAM J, Appl. Dyn. Syst. 7 (2008), 361-420.

\bibitem{KC} M. Kuwamura and H. Chiba, {\em Mixed-mode oscillations and chaos in a prey-predator system with dormancy of predators}, Chaos 19 (2009) 1-10.


\bibitem{KR}  Y. A. Kuznetsov and S. Rinaldi, {\em Remarks on food chain dynamics} Math. Biosci. 133 (1996) 1 - 33.

\bibitem{LRV} B. Letson, J. Rubin and T. Vo, {\em Analysis of interacting local oscillation
mechanisms in three-timescale systems}, SIAM J. Appl. Math. 77 3 (2017) 1020-1946. 

\bibitem{LXY}  W. Liu, D. Xiao, Y. Yi,  {\em Relaxation oscillations in a class of predator-prey systems}, J. Diff. Equ.188 (2003), 30-331. 


\bibitem{molleman}  F. Molleman, {\em Moving beyond phenology: new directions in the study of temporal dynamics of tropical insect communities}, Current Science 114 (2018). 


\bibitem{MR} S. Muratori and S. Rinaldi, {\em Remarks on competitive coexistence}, SIAM J.  Applied Math. 49 (1989) 1462-1472.

\bibitem{MRfoodchain} S. Muratori and S. Rinaldi, {\em Low- and high-frequency oscillations in three-dimensional food chain system}, J. Appl. Math. 52 (1992) 1688 - 1706.


\bibitem{neish1} A. I. Neishtadt, {\em Persistence of stability loss for dynamical bifurcations I}, Differ. Equ. 24 (1987) 23, 1385 -1391.

\bibitem{neish2} A. I. Neishtadt, {\em Persistence of stability loss for dynamical bifurcations II}, Differ. Equ. 24 (1988) 171 - 176 .

\bibitem{nelson} W.A. Nelson, O.N.  Bjornstad and T. Yamanaka, {\em Recurrent Insect Outbreaks Caused by Temperature-Driven Changes in System Stability}, Science, 314 (2013) 796-799.

\bibitem{PAGK} J.C. Poggiale, C. Aldebert, B. Girardot, and  B.W. Kooi,  {\em Analysis of a predator-prey model with specific time scales: a geometrical approach proving the occurrence of canard solutions}, J. of Math. Bio.  80, (2020) 39 - 60.

 
\bibitem{RMpopbio} S. Rinaldi, and S. Muratori,  {\em Limit cycles in slow-fast forest-pest models}, Theor. Popul. Biol. 41, (1992) 26 - 43.


 \bibitem{SCT} S. Sadhu and S. Chakraborty Thakur, {\em Uncertainty and Predictability in Population Dynamics of a Two-trophic Ecological Model: Mixed-mode Oscillations, Bistability and Sensitivity to Parameters}, Ecological Complexity 32 (2017) 196-208.
 
 
\bibitem{Sadhudcds} S. Sadhu, {\em Complex oscillatory patterns near singular Hopf bifurcation in a two-timescale ecosystem}, Discrete \& Continuous Dynamical Systems - B, 26 (2021) 5251 - 
5279.
 
  \bibitem{Sadhunew} S. Sadhu, {\em Analysis of the onset of a regime shift and detecting early warning signs of major population changes in a two-trophic three species predator-prey model with long-term transients}, J. Math. Biol. (in print).

\bibitem {seoetal} G. Seo and G. Wolkowicz, {\em Pest control by generalist parasitoids: a bifurcation theory approach}, DCDS-S, 13 (11) (2020) 3157 - 3187.

\bibitem{singleton} G. R. Singleton, P.R.  Brown, R.P. Pech, J. Jacob, G.J. Mutze and C. J. Krebs, {\em One hundred years of eruptions of house mice in Australia - a natural biological curio}, Biol. J. of  Linnean Soc., 84, (3) (2005) 617 - 627.


\bibitem{stenseth} N. C. Stenseth, {\em Population Cycles in Voles and Lemmings: Density Dependence and Phase Dependence in a Stochastic World}, Oikos, 87 (3) (1999) 427-461.

\bibitem{straube} R. Straube, D. Flockerzi and M.J.B. Hauser, {\em Sub-Hopf/fold-cycle bursting and its relation to (quasi-)periodic oscillations}, J. Phys.: Conf. Ser. 55 020 (2006).

\bibitem{TTVWB}  W. Teka,  J. Tabak,  T. Vo, M. Wechselberger,  R.  Bertram {\em The dynamics underlying pseudo-plateau bursting in a pituitary cell model}, J. Math. Neuro. Sc. 1: 12 (2011).

\bibitem{TL} R. Tyson and F. Lutscher, {\em Seasonally varying predation behaviour and climate shifts are predicted to affect predator-prey cycles}, Am. Nat. (2016) 188 (5) 539 - 553.

\bibitem{VBW} T. Vo, R. Bertram, and M. Wechselberger, {Multiple geometric viewpoints of
mixed mode dynamics associated with pseudo-plateau bursting}, SIAM J. Appl.
Dyn. Syst. 12 (2013) 789 - 830.

\bibitem{Vo} T. Vo, {\em Generic torus canards},  Physica D: Nonlinear Phenomena, 356 (2017) 37 -64.

\bibitem{WZ}  C. Wang  and X. Zhang, {\em Canards, heteroclinic and homoclinic orbits for a slow-fast predator-prey model of generalized Holling type III}, J. Diff. Eqns. 267, 6 (2019)  3397 - 3441.

\bibitem{Wesc} M. Wechselberger, {\em Existence and bifurcation of canards in $\mathbb{R}^3$ in the case of a folded node}, SIAM J. Appl. Dyn. Syst. 4 (1) (2005) 101 - 139. 



\end{thebibliography}
\end{document}